\documentclass[11pt]{amsart}
\usepackage{amsmath,amscd,amssymb,amsthm}

\def\cal{\mathcal}

\def\C{{\cal C}}

\def\I{{\cal I}}

\def\K{{\cal K}}

\def\S{{\cal S}}

\def\V{{\cal V}}

\def\frak{\mathfrak}

\def\da{{\frak a}}

\def\dF{{\frak F}}

\def\dg{{\frak g}}

\def\dh{{\frak h}}

\def\dl{{\frak l}}

\def\dgo{{\frak o}}
\def\dgp{{\frak p}}

\def\ds{{\frak s}}
\def\dsl{\ds\dl}

\def\dt{{\frak t}}
\def\du{{\frak u}}
\def\dU{{\frak U}}

\def\dZ{{\frak Z}}
\def\osp{\dgo\ds\dgp}

\def\Bbb{\mathbb}

\def\bC{\Bbb C}
\def\bD{\Bbb D}

\def\bF{\Bbb F}
\def\bH{\Bbb H}
\def\bJ{\Bbb J}
\def\bL{\Bbb L}
\def\bN{\Bbb N}

\def\bP{\Bbb P}
\def\bQ{\Bbb Q}
\def\bR{\Bbb R}

\def\bX{\Bbb X}
\def\bZ{\Bbb Z}

\def\ep{\epsilon}

\def\h{\hat}
\def\la{\langle}
\def\lb{\left[}

\def\lf{\lfloor}
\def\o{\overline}
\def\ot{\otimes}

\def\ra{\rangle}
\def\rb{\right]} 

\def\rf{\rfloor}
\def\scong{\buildrel\ds\over\cong}

\def\t{\tilde}
\def\tcong{\buildrel\dt\over\cong}

\def\ad{\mathop{\rm ad}\nolimits}
\def\and{\mbox{\rm \ and\ }}

\def\Ber{\mathop{\rm Ber}\nolimits}

\def\Comp{\mathop{\rm Comp}\nolimits}
\def\CQ{\mathop{\rm CQ}\nolimits}

\def\CS{\mathop{\rm CS}\nolimits}
\def\Diff{\mathop{\rm Diff}\nolimits}

\def\dim{\mathop{\rm dim}\nolimits}

\def\End{\mathop{\rm End}\nolimits}

\def\even{{\mathop{\rm even}\nolimits}}
\def\Ext{\mathop{\rm Ext}\nolimits}

\def\Hom{\mathop{\rm Hom}\nolimits}

\def\Im{\mathop{\rm Im}\nolimits}

\def\mod{\mathop{\rm mod}\nolimits}

\def\NS{\mathop{\rm NS}\nolimits}

\def\odd{{\mathop{\rm odd}\nolimits}}
\def\odU{\o{{\frak U}}{}}

\def\Poly{\mathop{\rm Poly}\nolimits}
\def\px{\partial_x}
\def\pxi{\partial_\xi}

\def\reg{\mathop{\rm reg}\nolimits}

\def\Roo{{\bR^{1|1}}}

\def\SBol{\mathop{\rm SBol}\nolimits}

\def\scrm{\scriptsize\rm}

\def\SNCR{\mathop{\rm SNCR}\nolimits}

\def\Soo{{S^{1|1}}}
\def\SP{\mathop{\rm SP}\nolimits}

\def\Span{\mathop{\rm Span}\nolimits}
\def\SQ{\mathop{\rm SQ}\nolimits}

\def\Stabilizer{\mathop{\rm Stabilizer}\nolimits}

\def\Sym{\mathop{\rm Sym}\nolimits}
\def\Symb{\mathop{\rm Symb}\nolimits}

\def\triv{{\mbox{\scrm triv}}}
\def\ts{\textstyle}

\def\Vec{\mathop{\rm Vec}\nolimits}

\def\VR{\Vec(\Bbb R)}
\def\VRM{\Vec(\Bbb R^m)}
\def\VRoo{\Vec(\Roo)}

\def\cind{completely indecomposable}

\def\dog{differential operator}

\def\expro{extremal projector}

\def\ie{{\em i.e.,\/}}
\def\iff{if and only if}

\def\ind{indecomposable}

\def\irr{irreducible}

\def\lsa{Lie superalgebra}
\def\lwv{lowest weight vector}
\def\meno{\medbreak\noindent}

\def\PBW{Poincar\'e-Birkhoff-Witt}

\def\psdog{pseudo\dog}
\def\psidog{$\Psi$DO}

\def\r{representation}

\def\sq{subquotient}

\def\tdm{tensor density module}

\def\tfm{tensor field module}
\def\th{\thinspace}

\def\uea{universal enveloping algebra}

\def\vf{vector field}
\def\vs1{$\Vec(S^1)$}

\def\PR{\bC[x]}
\def\PRoo{\bC[x,\xi]}
\def\PSoo{{\mathbb C}[x, x^{-1}, \xi]}

\newtheorem{lemma}{Lemma}[section]

\newtheorem{thm}[lemma]{Theorem}
\newtheorem{prop}[lemma]{Proposition}
\newtheorem{cor}[lemma]{Corollary}

\title[Conformal symbols over $\Roo$]{Conformal symbols and the action of contact vector fields over the superline}

\author{Charles H.\ Conley}
\address{Department of Mathematics\\University of North Texas\\Denton TX 76203}
\email{conley@unt.edu}

\begin{document}

\begin{abstract}
Let $\K$ be the \lsa\ of contact \vf s on the supersymmetric line $\Roo$.  We compute the action of $\K$ on the modules of differential and \psdog s between spaces of tensor densities, in terms of their conformal symbols.  As applications we deduce the geometric subsymbols, 1-cohomology, and various uniserial subquotients of these modules.  We also outline the computation of the $\K$-equivalences and symmetries of their \sq s. \end{abstract}

\maketitle

\section{Introduction}  \label{Intro} 
\setcounter{lemma}{0}

Suppose that $\V$ is an infinite dimensional Lie algebra (or superalgebra) of \vf s on a manifold $M$, for example $\Vec(M)$ itself or its contact subalgebra.  Consider the adjoint action of $\V$ on the algebra $\Diff(M)$ of \dog s on $M$.  Note that this action preserves the degree filtration $\Diff^k(M)$.  Assume that $\V$ contains a distinguished maximal finite dimensional subalgebra $\da$, its {\em conformal\/} or {\em projective\/} subalgebra, under whose action the degree filtration has a unique invariant splitting.  This splitting is known as the {\em conformal\/} or {\em projective symbol.\/}

Over the past decade there have been many papers concerned with the problem of computing the action of $\V$ on $\Diff(M)$ in terms of the conformal symbol.  The analysis leads to the more general setting of modules of \dog s between spaces of {\em tensor fields,\/} because the symbol modules $\Diff^k(M)/\Diff^{k-1}(M)$ are \tfm s.

In the seminal paper \cite{CMZ97}, Cohen, Manin, and Zagier solved this problem completely for $M=\Bbb R$ and $\V=\VR$, where $\da\cong\ds\dl_2$ is the projective subalgebra (in fact, they went further and computed composition of \dog s in terms of the conformal symbol).  Lecomte and Ovsienko have studied the case that $M=\Bbb R^m$ and $\V=\VRM$, where $\da\cong\dsl_{m+1}$ is the projective subalgebra \cite{LO99}.  

Results on this problem have several applications, for example to modular forms \cite{CMZ97}, the $\Diff(M)$-valued cohomology of $\V$ \cite{FF80, LO00}, symmetries and equivalences of \sq s of \dog\ modules \cite{DO97, LO99}, and \ind\ bounded modules \cite{FF80, Co05}.

To our knowledge, the first examination of a case in which $\V$ is not all of $\Vec(M)$ was commenced in \cite{CMZ97}.  Toward the end of that paper the \lsa\ $\K$ of contact \vf s on the supersymmetric line $\Roo$ is considered: the conformal symbol is computed for $\Diff(\Roo)$.  This result has recently been generalized to \dog s between generic pairs of \tdm s \cite{GMO07}: only {\em resonant\/} pairs are excluded.  It is this generic ({\em non-resonant\/}) case which we study here.  We take the approach used in \cite{Co01} for $\VR$: we use the {\em step algebra\/} of the conformal subalgebra, a copy of $\osp_{1|2}$, to deduce formulas for the action of $\K$ in terms of the conformal symbol without actually using the explicit formula for the conformal symbol.  We work only with polynomials and our definitions are algebraic.  For a more geometric point of view, see \cite{GMO07}.

It should be noted that there is a large body of work concerning $\K$ and its analog over the supersymmetric circle $\Soo$; see for example the references in \cite{GMO07}.  Ovsienko has initiated the study of the conformal symbol for general contact manifolds, which promises to be a rich field \cite{Ov06}.

This paper is organized as follows.  In Section~\ref{DOG} we define $\K$ and the modules of tensor densities and \dog s on which it acts.  This section concludes with the important {\em fine filtration\/} discovered in \cite{GMO07}, a $\K$-invariant filtration of the \dog\ modules twice as fine as the obvious order filtration.  Its \sq s are \tdm s.  

Section~\ref{OSP12} gathers various results on $\osp_{1|2}$-modules.  Section \ref{GenCoho} contains some generalities on cohomology and \r s of finite length of \lsa s, along with preliminary results on the $\K$-cohomology of the \tdm s.  In particular, we make a conjecture generalizing to $\K$ Goncharova's description of the cohomology of the \tdm s of $\VR$ \cite{Gon73}.

In Section~\ref{CS} we state our main result, the formula for the action of $\K$ on non-resonant modules of \dog s in terms of the conformal symbol.  The formula is complicated and we have only simplified it fully on length~6 \sq s.  In Section~\ref{PSDOG} we generalize to \psdog s, and deduce certain symmetries of the $\K$-action with respect to the conformal symbol arising from the Berezinian and the super Adler trace.

Section~\ref{Apps} gives various applications of the main result.  First we show that under the action of $\K$, \psdog s have an invariant ``sesquisymbol'' extending their usual symbol.  Then in Section~\ref{ITs} we use the super Adler trace to relate our work to $\K$-invariant maps from the tensor product of two \tdm s to a third.  Coupled with our main result, this gives a classification of such maps in the non-resonant case, reproducing a result of Leites, Kochetkov, and Weintrob \cite{LKW91}.

In Section~\ref{KCoho} we compute the conformally relative $\Ext^1$ groups between \tdm s of $\K$, and give some information on the $\Ext^2$ groups.  In particular, we find analogs of the famous ``$\sqrt{19}$'' 1-cocycles of $\VR$ discovered by Feigin and Fuchs \cite{FF80}, and make a conjecture concerning the generalization to $\K$ of their description of the cohomology of the \dog\ modules of $\VR$.  In Section~\ref{Uniserial} we follow the approach of \cite{Co05} to construct various uniserial (\ie\ \cind) modules composed of tensor density modules, as \sq s of \psdog\ modules.

In Section~\ref{Equivs} we discuss equivalences between \sq s of \psdog\ modules, in the spirit of \cite{DO97} and \cite{LO99}.  Corresponding to each composition series there is a 2-parameter family of \sq s.  Within each family, subquotients of length~$<6$ are generically equivalent, but there are exceptional cases.  We expect that \sq s of length~$>6$ are equivalent only to their conjugates, but we have not yet proven this.  

A new phenomenon occurs in length~6: the equivalence classes within a given family are curves.  In appropriate coordinates they form the pencil of conics through four fixed points which depend only on the composition series.  However, the computations are quite long and there are related results over $\VR$ in length~5, so we will present these both together in a future paper.

We also compute the algebras of symmetries, \ie\ $\K$-invariant endomorphisms, of \sq s of \psdog\ modules.  In the non-resonant case these are always of the form $\bC^e$, where $e$ is the number of \ind\ summands of the \sq.

Section~\ref{Proofs} contains the proof of the main results.  Section~\ref{Resonant} discusses the resonant case and non-$\osp_{1|2}$-relative 1-cohomology.  Section~\ref{LSA} is an appendix on representations of \lsa s in general.

\section{Differential operators on $\Roo$}  \label{DOG}
\setcounter{lemma}{0}

Throughout this paper we work only with vector spaces over $\bC$.  We use $\bZ^+$, $\bZ^-$, and $\bN$ for the positive, negative, and non-negative integers, respectively, and we write $\lf x\rf$ and $\{x\}$ for the integer and fractional parts of $x\in\bR$, respectively.  We use the following notation for the {\em falling factorial\/} and the {\em generalized binomial symbol:\/} for $z\in\bC$ and $n\in\bN$,
$$ \ts \lb{z\atop n}\rb := z(z-1)\cdots(z-n+1), \quad
   {z\choose n} := \frac{1}{n!} \lb{z\atop n}\rb. $$

Given a superspace $V$, we write $V_\even$ and $V_\odd$ for its even and odd parts instead of the usual $V_0$ and $V_1$, in order to avoid confusion with other numerical subscripts.  We denote the parity endomorphism $v\mapsto (-1)^{|v|}v$ of $V$ by $\ep_V: V\to V$ (we often write simply $\ep$ when $V$ is specified by the context).  If $V$ is a \r\ of a \lsa\ $\dg$, we denote the subspace of $\dg$-invariant vectors by $V^\dg$.  We will need the parity-reversing functor $\Pi$, which we write as a superscript:
\begin{displaymath}
   V^\Pi_\even := V_\odd,\quad V^\Pi_\odd := V_\even.
\end{displaymath}
Define $V^{k\Pi}$ to be $V$ if $k$ is even and $V^\Pi$ if $k$ is odd.  Write $\dim V$ for the dimension of $V$.

In Section~\ref{LSA} we have given a review of some basic constructions involving \r s of \lsa s, such as duals, tensor products, supersymmetric and superalternating tensor algebras, and invariant inner products.  In particular, we discuss the interactions of these objects with $\Pi$.  We also discuss \uea s: the \PBW\ theorem, centers and supercenters, and automorphisms.  We will refer to this section regularly.

\subsection{The superline}
Let $\Roo$ be the superline, with even coordinate $x$ and odd coordinate $\xi$.  Here $\xi^2=0$, so $\PRoo$ has basis $\{1, \xi\}$ over $\PR$.

The space of polynomial \vf s on $\Roo$ is 
$$ \VRoo := \Span_{\PRoo} \bigl\{ \px, \pxi \bigr\}. $$  
It is a Lie superalgebra acting by superderivations on $\PRoo$:
$$ [X,Y] := XY - (-1)^{|X||Y|}YX, \quad
   X(FG) := X(F)G + \ep^{|X|}(F) X(G). $$
When we write an expression such as $XF$, the $F$ should be interpreted as a multiplicative operator: $XF(G):=X(FG)$.  Thus the superbracket $[X,F]$ is $X(F)$.  We will write simply $F'$ for $\px F$.

\subsection{Contact vector fields} \label{K}
Define elements $D$ and $\o D$ of $\VRoo$ by
$$ D := \pxi + \xi\px, \quad \o D := \pxi - \xi\px. $$
Note that $\pxi$, $D$, and $\o D$ are odd and satisfy the formulas $\pxi^2=0$,
\begin{equation} \label{DoD}
   [D,D] = 2D^2 = 2\px, \quad
   [D,\o D] =0, \quad
   [\o D,\o D] = 2\o D^2 = -2\px.
\end{equation}

We will be interested in the contact structure on $\Roo$ defined by the form
$$ \omega := dx+\xi d\xi. $$
The corresponding completely non-integrable distribution is generated by $\o D$.  Thus the space of vector fields tangent to the distribution is $\PRoo\o D$, the vector fields annihilated by $\omega$.

\meno {\bf Definition.} {\em
The Lie superalgebra of\/ {\em contact vector fields on $\Roo$\/} is\/}
$$ \K := \Stabilizer_{\VRoo}\bigl(\PRoo\o D\bigr) 
   = \Stabilizer_{\VRoo}\bigl(\PRoo\omega\bigr). $$

\medbreak
The reader may check that $\K_\odd=\PR D$, $\K_\odd$ generates $\K$, and $\K$ is the image of the even linear injection $\bX:\PRoo\to\VRoo$ defined as follows:
\begin{equation} \label{CH}
   \mbox{\rm For\ } f\in\PR,\ \ 
   \bX(f) := f\partial_x + \ts\frac{1}{2}f'\xi\partial_\xi
   \mbox{\rm\ \ and\ \ }
   \bX(\xi f)\ :=\ \ts\frac{1}{2}fD.
\end{equation}

Given $X\in\VRoo$, the function $\bX^{-1}(X)$ in $\PRoo$ is called the {\em contact Hamiltonian\/} of $X$.  We remark in passing that the {\em contact bracket\/} $\{F,G\} := \bX^{-1}[\bX(F),\bX(G)]$ on $\PRoo$ is given by
\begin{equation} \label{CB}
   \{F,G\} = FG' - F'G - {\ts\frac{1}{2}}(-1)^{|G|}D(F)D(G).
\end{equation}
Equations~(\ref{CH}) and~(\ref{CB}) match~(2.5) and~(2.7) of \cite{GMO07}, respectively.  

Note that $\VRoo=\K\oplus\PRoo\o D$, but $\K$ is not closed under multiplication by elements of $\PRoo$.  For reference, $f, g\in\PR$ give
\begin{equation} \label{XB}
\begin{array}{rcl}
   \bigl[\bX(f),\bX(g)\bigr] &=& \bX(fg'-f'g), \\[4pt]
   \bigl[\bX(f),\bX(\xi g)\bigr] &=& \bX(fg' - \ts\frac{1}{2}f'g), \\[4pt]
   \bigl[\bX(\xi f),\bX(\xi g)\bigr] &=& \bX(\ts\frac{1}{2}fg).
\end{array}
\end{equation}

We will use a Virasoro-like basis of $\K$: for $n\in \bN$,
$$ e_{n-1}:=\bX(x^n),\qquad e_{n-1/2}:=2\bX(\xi x^n). $$
Its brackets are
\begin{displaymath} \begin{array}{rcll}
   [e_n,e_m] &=& (m-n)e_{n+m} &\mbox{\rm\ if $n,m\in\bN-1$}, \\[4pt]
   [e_n,e_m] &=& (m-n/2)e_{n+m} &\mbox{\rm\ if $n\in\bN-1$ 
   and $m\in\bN-1/2$}, \\[4pt]
   [e_n,e_m] &=& 2e_{n+m} &\mbox{\rm\ if $n,m\in\bN-1/2$}.
\end{array} \end{displaymath}

\subsection{Tensor density modules} \label{Flambda}
By definition, $\PRoo\omega$ is a module of $\K$.  The \tdm s are the associated 1-parameter family of scalar powers $\bigl(\PRoo\omega\bigr)^{\ot\lambda}$.  (They may also be viewed as the 1-parameter family of deformations of $\PRoo$ associated to the 1-dimensional space $H^1(\K,\PRoo)$; see the remark at the end of Section~\ref{Coho} and Proposition~\ref{F012coho}bc.)  Thus:

\meno{\bf Definition.} {\em
For $\lambda\in\bC$, the\/ {\em \tdm\/} $\bigl(\pi_\lambda,\bF(\lambda)\bigr)$ of $\K$ is\/}
$$ \bF(\lambda) := \PRoo\omega^\lambda, \quad
   \pi_\lambda\bigl(\bX(F)\bigr) (\omega^\lambda G) := 
   \omega^\lambda \bigl(\bX(F)G + \lambda F'G \bigr). $$

\medbreak
Note that since $\omega$ is even, $\bF_\even(\lambda) = \PR\omega^\lambda$ and $\bF_\odd(\lambda) = \PR\xi\omega^\lambda$.  We give the operators $\pi_\lambda(e_n)$ for reference: for $n\in\bN$,
\begin{displaymath}
   \pi_\lambda(e_{n-1}) = x^n\partial_x + nx^{n-1}(\lambda+\ts\frac{1}{2}\xi\partial_\xi),
   \quad \pi_\lambda(e_{n-1/2}) = x^nD+2n\lambda\xi x^{n-1}.
\end{displaymath}

The following lemma redefines $\bX$, making the important point that as $\K$-modules, $\K$ and $\bF(-1)$ are equivalent.  Its proof is~(\ref{XB}).

\begin{lemma} \label{Xomega}
Henceforth regard $\bX$ as the map $\omega^{-1}F\mapsto \bX(F)$ from $\bF(-1)$ to $\K$.  As such, it is a $\K$-equivalence.
\end{lemma}

\subsection{Differential operators} \label{DOGs}
Our focus in this paper is on the action of $\K$ on the \dog s between \tdm s.  Its study will lead us to consider the action of $\K$ on arbitrary linear maps between \tdm s of both standard and reversed parity, so for brevity we make the following definition.

\meno {\bf Definition.} {\em
Set $\bH(\lambda,p) := \Hom\bigl(\bF(\lambda),\bF(\lambda+p)\bigr)$.  Let $\sigma_{\lambda,p}$ and $\sigma'_{\lambda,p}$ be the actions $\Hom(\pi_\lambda,\pi_{\lambda+p})$ and $\Hom(\pi_\lambda,\pi^\Pi_{\lambda+p})$ of $\K$ on $\bH(\lambda,p)$ and $\bH(\lambda,p)^\Pi$, respectively (see~(\ref{superhom})).}

\meno {\bf Remark.}
In the category of superspaces, $\Hom\bigl(\bF(\lambda)^\Pi,\bF(\lambda+p)^\Pi\bigr)$ is equal to $\bH(\lambda,p)$, and both $\Hom\bigl(\bF(\lambda),\bF(\lambda+p)^\Pi\bigr)$ and $\Hom\bigl(\bF(\lambda)^\Pi,\bF(\lambda+p)\bigr)$ are equal to $\bH(\lambda,p)^\Pi$ (see Section~\ref{Parity}).  However, care must be taken with the $\K$-actions: $\Hom(\pi_\lambda^\Pi,\pi_{\lambda+p}^\Pi)$ and $\Hom(\pi_\lambda^\Pi,\pi_{\lambda+p})$ are equal to $\sigma_{\lambda,p}$ and $\sigma'_{\lambda,p}$, respectively, but the actions $\sigma_{\lambda,p}^\Pi$ and $\sigma'_{\lambda,p}$ are not equal.  By Lemma~\ref{leprep}, the map $L(\ep_{\bF(\lambda+p)})=R(\ep_{\bF(\lambda)})$ is an odd $\K$-equivalence from $\sigma_{\lambda,p}$ to $\sigma'_{\lambda,p}$.

\medbreak
Now we turn to \dog s.  For $G\in\PRoo$, $X\in\VRoo$, $p\in\bC$, and $j\in \bN$, define
\begin{equation} \label{monodog}
   \omega^p G X^j:\bF(\lambda)\to\bF(\lambda+p), \quad
   (\omega^p G X^j) (\omega^\lambda F) := \omega^{\lambda+p} G X^j(F).
\end{equation}

\meno{\bf Definition.} {\em
The space $\bD(\lambda,p)$ of\/ {\em polynomial \dog s from $\bF(\lambda)$ to $\bF(\lambda+p)$\/} is the span of all operators as in~(\ref{monodog}).  In light of~(\ref{DoD}),}
$$ \bD(\lambda,p) = \Span\bigl\{\omega^p G \o D^j: 
   \th G\in\PRoo,\th j\in\bN\bigr\}. $$

\medbreak
Since $\pi_\lambda(X)$ is itself a polynomial \dog, $\bD(\lambda,p)$ is invariant under $\sigma_{\lambda,p}$.  Moreover, $\ep_{\bF(\lambda)}=1-2\xi\partial_\xi$ is also a polynomial \dog, so $\bD(\lambda,p)^\Pi$ is invariant under $\sigma'_{\lambda,p}$.  When we write $\bH(\lambda,p)^\Pi$ or $\bD(\lambda,p)^\Pi$, we are regarding them as $\K$-modules under the action $\sigma'_{\lambda,p}$.

Note that composition of \dog s makes $\bigoplus_{\lambda,p} \bD(\lambda,p)$ an algebra, on which the \r\ $\bigoplus_{\lambda,p} \sigma_{\lambda,p}$ acts by superderivations.  Therefore we have the following formula for the action of $\K_\odd$ on $\bD(\lambda,p)$: for $f\in\PR$ and $G\in\PRoo$, $\sigma_{\lambda,p}(fD)(\omega^p G \o D^j) = $
$$ \sigma_{\lambda,p}(fD)(\omega^p)G\o D^j +
   \omega^p \sigma_{\lambda,0}(fD)(G) \o D^j +
   \omega^p \ep(G) \sigma_{\lambda,0}(fD)(\o D^j). $$
Since $\pi_\lambda(fD) = fD + 2\lambda\xi f'$, this yields $\sigma_{\lambda,p}(fD)(\omega^p G \o D^j) = $
\begin{equation} \label{sigac1}
   \omega^p\Bigl\{\bigl(2p\xi f'G + fD(G)\bigr) \o D^j + \ep(G)
   [fD+2\lambda\xi f', \o D^j] \Bigr\}.
\end{equation}
The useful facts $[D,\o D]=0$, $D = \o D - 2\xi\o D^2$, and $fD = Df - \xi f'$ permit us to rewrite the last term of~(\ref{sigac1}):
\begin{equation} \label{sigac2}
   [fD+2\lambda\xi f', \o D^j] =
   - (\o D - 2\xi\o D^2)[\o D^j,f] - (-1)^j(2\lambda - 1) [\o D^j,\xi f'].
\end{equation}

\subsection{The fine filtration} \label{FF}
Let $\bD^k(\lambda,p)$ be the total order filtration on $\bD(\lambda,p)$, a $\VRoo$-invariant $\bN$-filtration.  It was discovered in \cite{GMO07} that it has an important refinement to a $\K$-invariant $\bN/2$-filtration, the {\em fine filtration.\/}  Geometrically, for $j\in\bZ^+$ the space $\bD^{j-1/2}(\lambda,p)$ consists of order~$j$ operators whose symbol is tangent to the contact distribution.  Algebraically,
$$ \bD^k(\lambda,p) := \Span\bigl\{\omega^p G \o D^j: G\in\PRoo,\ j\le 2k\bigr\},
   \quad k\in \bN/2. $$

\begin{prop} \label{fine} {\rm \cite{GMO07}}
The fine filtration is $\K$-invariant and compatible with composition: $\bD^{k'}(\lambda+p,p')\circ\bD^k(\lambda,p) = \bD^{k'+k}(\lambda,p'+p)$.  Composition is commutative up to fine symbol.  For $k\in\bN/2$ there is an even $\K$-equivalence from $\bD^k(\lambda,p)/\bD^{k-1/2}(\lambda,p)$ to $\bF(p-k)^{2k\Pi}$, defined by
$$ \omega^p G \o D^{2k} + \bD^{k-1/2}(\lambda,p) \mapsto \omega^{p-k} G. $$
\end{prop}

\meno{\it Proof.\/}
Composition compatibility and fine symbol commutativity are clear.  To prove $\K$-invariance, use fine symbol commutativity together with the facts that $\sigma_{\lambda,p}$ acts by superderivations and that the spaces $\PRoo\omega^p$ and $\PRoo\o D$ are $\K$-invariant.  To prove the fine symbol equivalence, verify that $G\o D\mapsto \omega^{-1/2} G$ is a $\K$-equivalence from $\PRoo\o D$ to $\bF(-1/2)$ and then use the argument in the last sentence.  $\Box$

\meno {\bf Remark.}
The action $\sigma'_{\lambda,p}$ does not preserve the filtration $\bD^k(\lambda,p)^\Pi$.  However, using the fact that $L(\ep_{\bF(\lambda+p)})$ is an odd equivalence from $\sigma_{\lambda,p}$ to $\sigma'_{\lambda,p}$ together with $\ep = 1-2\xi\o D$, one can show that $\sigma'_{\lambda,p}(\Omega)$ increases the degree by at most~1 for all $\Omega\in\dU(\dg)$.

\section{Representations of $\osp_{1|2}$}  \label{OSP12}
\setcounter{lemma}{0}

The {\em conformal\/} subalgebra $\ds$ of $\K$ consists of the infinitesimal linear fractional contact transformations of $\Roo$ (see \cite{GMO07} or \cite{CMZ97}):
\begin{equation} \label{sdefn}
   \ds := \Span\{e_{-1}, e_{-1/2}, e_0, e_{1/2}, e_1\},
\end{equation}
It is a maximal subalgebra of $\K$ and is isomorphic to $\osp_{1|2}$.  In this section we collect various properties of $\ds$-modules.  Their proofs are simple and for the most part well-known, so they are sketched or omitted.

\meno{\bf Definition.} {\em
The\/ {\em affine subalgebra\/} $\dt$ of\/ $\ds$ and the nilradical $\du$\/ of $\dt$ are
\begin{displaymath}
   \dt := \Span\{e_{-1}, e_{-1/2}, e_0\}, \quad
   \du := \Span\{e_{-1}, e_{-1/2}\}.
\end{displaymath}

In any $\ds$-module $V$, the eigenvalues and eigenspaces of $e_0$ are called {\em weights\/} and {\em weight spaces.\/}  We write $V_\mu$ for the $\mu$-weight space.  Weight vectors annihilated by $e_{-1/2}$ (respectively, $e_{1/2}$) are called {\em lowest\/} (respectively, {\em highest\/}) weight vectors.\/}

\subsection{Properties of $\bF(\lambda)$, $\bH(\lambda,p)$, and $\bD(\lambda,p)$ as $\ds$-modules}
\label{FDH} 
For reference, the action of $\ds$ on $\bF(\lambda)$ is given by
\begin{equation} \label{saction}
\begin{array}{c}
   \pi_\lambda(e_{-1/2}) = D, \quad
   \pi_\lambda(e_{-1}) = \partial_x, \quad
   \pi_\lambda(e_{1/2}) = xD+2\lambda\xi, \\[6pt]
   \pi_\lambda(e_1) = x(x\partial_x + \xi\partial_\xi +2\lambda), \quad
   \pi_\lambda(e_0) = x\partial_x+\ts\frac{1}{2}\xi\partial_\xi + \lambda.
\end{array}
\end{equation}
Note that $e_0$ acts as an Euler operator: $x$ has weight~$1$, $\xi$ has weight~$1/2$, $D$ and $\o D$ have weight~$-1/2$, $\omega^\lambda$ has weight~$\lambda$, {\em etc.\/}

The next lemma gives the $\ds$-structure of the $\bF(\lambda)$; the proof is an easy exercise (see Section~\ref{Parity} for the parity functor $\Pi$).  For $\lambda\in\bN/2$, define $\bL(\lambda)$ to be the following $4\lambda+1$-dimensional subspace of $\bF(-\lambda)$:
$$ \bL(\lambda) := \omega^{-\lambda}\Span\bigl\{1, \xi, x, \xi x, x^2, 
   \ldots, x^{2\lambda-1}, \xi x^{2\lambda-1}, x^{2\lambda}\bigr\}. $$

\begin{lemma} \label{F}
\begin{enumerate}

\item[(a)] $\bF(\lambda)^\du=\bC\omega^\lambda$.

\smallbreak
\item[(b)] $\bF(\lambda)^\dt= \bF(\lambda)^\ds= \bF(\lambda)^\K$.  It is $\bC 1$ if $\lambda=0$, and~$0$ otherwise.

\smallbreak
\item[(c)] $\bF(\lambda)^{e_{1/2}}$ is $\bC\omega^\lambda x^{-2\lambda}$ for $\lambda\in-\bN/2$, and~$0$ otherwise.

\smallbreak
\item[(d)] $\bF(\lambda)$ is $\ds$-\irr\ unless $\lambda\in -\bN/2$, when $\bL(-\lambda)$ is an \irr\ $\ds$-submodule and $\bF(\lambda)/\bL(-\lambda)$ is $\ds$-equivalent to $\bF(-\lambda+1/2)^\Pi$.

\end{enumerate}
\end{lemma}

Now we describe the $\ds$-structure of $\bD(\lambda,p)$.  The {\em super Bol operator\/} defined in \cite{Gi93, GT93} plays an important role:

\meno{\bf Definition.} {\em
For $p\in \bN/2$, the\/ {\em affine super Bol operator\/} of\/ $\bD(\lambda,p)$ is\/}
$$ \SBol_p(\lambda) := \omega^p\o D^{2p}. $$

\begin{lemma} \label{Bol}
\begin{enumerate}

\item[(a)] $\bH(\lambda,p)^\du= \omega^p \bC[[\o D]]$, and $\bD(\lambda,p)^\du = \omega^p \bC[\o D]$.

\smallbreak
\item[(b)] $\bH(\lambda,p)^\dt = \bD(\lambda,p)^\dt$.  It is $\bC \SBol_p$ if $p\in\bN/2$, and~$0$ otherwise.  For $p\in\bN/2$, $\bigl(\bH(\lambda,p)^{2p\Pi}\bigr)^\dt = \bC\omega^p D^{2p}$.

\smallbreak
\item[(c)] $\bH(\lambda,p)^\ds = \bD(\lambda,p)^\ds$.  It is $\bC \SBol_p$ if $p=0$ or if $p\in \frac{1}{2} + \bN$ and $\lambda = \frac{1}{4} - \frac{p}{2}$, and~$0$ otherwise.  When it is $\ds$-invariant, $\SBol_p(\lambda)$ has kernel\/ $\bL(-\lambda)$ and image\/ $\bF(-\lambda+1/2)$.

\smallbreak
\item[(d)] $\bH(\lambda,p)^\K = \bD(\lambda,p)^\K$ for all $(\lambda,p)$.  It is $\bC \SBol_p$ if $p=0$ or if $(\lambda,p)=(0,1/2)$, and~$0$ otherwise.

\end{enumerate}
\end{lemma}

\meno{\it Proof.\/}
For~(a), check that any $T$ in $\bH(\lambda,p)$ has a unique expression as a formal sum $\omega^p \sum_{i=0}^\infty T_i \o D^i$ for some $T_i\in\Poly(\Roo)$.  Hence $[D,T]= \omega^p \sum_i D(T_i) \o D^i$, which is zero \iff\ $T_i\in\bC$ for all~$i$.

The first two sentences of~(b) follow from~(a).  For the third, use Lemma~\ref{leprep} to see that the space in question is $\bC\SBol_p\circ\,\ep_{\bF(\lambda)}^{2p}$.  Then use $\ep_{\bF(\lambda)} = 1-2\xi\pxi$ to verify $\o D^{2p}\ep^{2p} = (-1)^{\lf p \rf} D^{2p}$.

For~(c), use~(b) and Lemma~\ref{F}d and note that when an $\ds$-map exists, it must be the unique $\dt$-map.  Or simply compute:
\begin{equation} \label{SBoleoh}
   \sigma_{\lambda,p}(e_{1/2}) \SBol_p =
   \bigl(p + (4\lambda-1)\{p\}\bigr) \omega^p \o D^{2p-1}.
\end{equation}

For~(d), use~(c) together with the fact that $e_{3/2}$ generates $\K$ under $\ds$ to come down to proving that $\sigma_{\lambda,p}(e_{3/2})\SBol_p(\lambda)=0$ \iff\ either $p=0$ or $p=1/2$ and $\lambda=0$.  Verify this by hand.  $\Box$

\subsection{The Casimir operator and its square root} \label{Casimir}
Let $\Lambda_\ds$ be the following even element of $\dU^2(\ds)$:
\begin{displaymath}
   \Lambda_\ds := \ts\frac{1}{2}(e_{-1/2}e_{1/2}-e_{1/2}e_{-1/2})
   = e_0 - e_{1/2}e_{-1/2}.
\end{displaymath}

\meno{\bf Definition.} {\em
In the terminology of \cite{ABF97}, the\/ {\em Scasimir operator\/} of $\ds$ is $T_\ds:=\Lambda_\ds-1/4$.  The\/ {\em Casimir operator\/} of $\ds$ is $Q_\ds:=T_\ds^2-1/16$.}

\begin{prop} \label{ghost}
The ghost and super centers of $\dU(\ds)$ are $\t\dZ(\ds)=\bC[T_\ds]$ and $\dZ(\ds)=\bC[T_\ds^2]=\bC[Q_\ds]$.  The actions of $T_\ds$ and $Q_\ds$ on the $\bF(\lambda)$ are
$$ \pi_\lambda(T_\ds) = (\lambda-{\ts\frac{1}{4}})\ep_{\bF(\lambda)}, \quad
   \pi_\lambda(Q_\ds) = \lambda^2-{\ts\frac{1}{2}}\lambda. $$
\end{prop}

The proof of this proposition is a relatively straightforward generalization of the proof that the Casimir operator $Q_{\ds_\even}=e_0^2-e_0-e_1e_{-1}$ of $\ds_\even$ (a copy of $\dsl_2$) generates the center of $\dU(\ds_\even)$.  Concerning $\pi_\lambda(T_\ds)$, recall that $\ep_{\bF(\lambda)}=1-2\xi\partial_\xi$.  Let us remark that $Q_{\ds_\even}= \Lambda_\ds^2-\Lambda_\ds$.  In particular, $\Lambda_\ds^2$ and hence also $Q_\ds$ are of degree~2.  

The existence of the Scasimir operator $T_\ds$, a square root of $Q_\ds+1/16$, is immediate from the relation between $Q_\ds$ and $Q_{\ds_\even}$ observed in \cite{Pi90}.  As far as we know, $T_\ds$ was first defined explicitly in \cite{Le95}.  A generalization to $\osp_{1|2n}$ was discovered independently in \cite{ABF97} and \cite{Mu97}.

\subsection{The extremal projector and the step algebra} \label{EPSA}
There is a large body of work on \expro s and step algebras.  To our knowledge the extremal projector of $\osp_{1|2}$ first appeared in \cite{BT81}, and step algebras go back to Mickelsson \cite{Mi73}.  The framework we give here is due to Zhelobenko; the results we state may be proven using the methods in \cite{Zh89}.  We have given parallel results for $\dsl_2\subset\VR$ in \cite{Co01}.

Let $M(\ds)$ be $\dU(\ds)/\dU(\ds)\du$, the {\em universal Verma module\/} of $\ds$.  Let $\dF(\ds)$ be the space of formal series in the $e_{1/2}^i e_{-1/2}^j$ with bounded $|i-j|$ and coefficients in $\bC[e_0]$.  Note that $\dF(\ds)$ is naturally an associative algebra containing $\dU(\ds)$ as a subalgebra and that $M(\ds)$ is an $\dF(\ds)$-module.

Given any free left $\bC[e_0]$ module $N$, let $\o N$ denote $\bC(e_0)\ot_{\bC[e_0]} N$, its extension by the fraction field $\bC(e_0)$.  Then $\o\dF(\ds)$ is an associative algebra and $\o M(\ds)$ is an $\o\dF(\ds)$-module.  The {\em \expro\/} $P(\ds)\in\o\dF(\ds)$ is
\begin{equation} \label{Psdfn}
   P(\ds) := \sum_{n=0}^\infty p_n(e_0) e_{1/2}^n e_{-1/2}^n, \mbox{\rm\ where\ }
   p_n := \prod_{j=1}^{\lf (n+1)/2 \rf} \frac{(j-2e_0)^{-1}}{\lf n/2 \rf!}.
\end{equation}

\begin{prop} \label{Ps}
\begin{enumerate}

\item[(a)]  $\o\dF(\ds)$ is naturally isomorphic to the algebra of right $\bC(e_0)$-linear endomorphisms of $\o M(\ds)$.

\smallbreak
\item[(b)]  $\o M(\ds) = \o M(\ds)^{\du} \oplus e_{1/2} \o M(\ds)$.

\smallbreak
\item[(c)]  $P(\ds)$ acts on $\o M(\ds)$ as projection to $\o M(\ds)^\du$ along $e_{1/2} \o M(\ds)$.  In particular, $P(\ds)^2=P(\ds)$ and $e_{-1/2}P(\ds)=0=P(\ds)e_{1/2}$.  The superautomorphism $j$ of Section~\ref{UEA} preserves $P(\ds)$.

\end{enumerate}
\end{prop}

Now suppose that $\dg$ is any \lsa\ containing $\ds$.  The {\em step algebra\/} of $\dg$ with respect to $\ds$ is $S(\dg,\ds) := \bigl(\dU(\dg)/\dU(\dg)\du\bigr)^\du$.  One checks that it is an associative algebra which acts on $V^\du$ for any $\dg$-module $V$.  Moreover, $\o\dF(\ds)$ acts on $\odU(\dg)/\odU(\dg)\du$, and $P(\ds)$ projects it to $\o S(\dg,\ds)$ along the image of $e_{1/2}$.  

We now state some results specific to $\K$.  For $\mu\in 1+\bZ^+/2$, use the coefficients $p_n$ of $P(\ds)$ to define elements $\t s_\mu$ of $\odU(\K)$ and $\t s_{\mu n}$ of $\bC(e_0)$:
\begin{equation} \label{ts}
   \t s_\mu := \sum_{n=0}^{2\mu-2} \t s_{\mu n}(e_0) e_{1/2}^n e_{\mu-n/2}
   := \sum_{n=0}^{2\mu} p_n(e_0) e_{1/2}^n \ad(e_{-1/2})^n e_\mu.
\end{equation}
To prove Lemma~\ref{step}, it is enough to note that $\t s_\mu\in\o S(\K,\ds)$ because $\t s_\mu \equiv P(\ds)e_\mu$ modulo $\odU(\K)\du$.  Lemma~\ref{345} is proven by direct computation.

\begin{lemma} \label{step}
Fix a \r\ $(\pi,V)$ of $\K$, and consider the space $V_\lambda^\du$ of $\lambda$-\lwv s.  If $p_{2\mu}$ has no pole at $\mu+\lambda$, then $\pi(\t s_\mu)$ is defined on $V_\lambda^\du$ and maps it to $V_{\lambda+\mu}^\du$.
\end{lemma}

\begin{lemma} \label{345}

\begin{enumerate}

\item[(a)]  
$\ad(e_{-1/2})^n e_\mu = C_{\mu n} e_{\mu-n/2}$, where 
\begin{displaymath}
   C_{\mu n} = 2^{-2\{n/2\}(-1)^{2\mu}}\lb {\lf \mu+1\rf \atop \lf(n+1)/2\rf - 2\{\mu n\}} \rb.
\end{displaymath}

\smallbreak \item[(b)]
$\t s_{\mu n} = C_{\mu n}\th p_n$ for $0\le n\le 2\mu-3$, and
$$ \t s_{\mu, 2\mu-2} = C_{\mu, 2\mu-2}\th p_{2\mu-2} + C_{\mu, 2\mu-1}\th p_{2\mu-1}
   + (e_0-\mu) C_{\mu, 2\mu}\th p_{2\mu}. $$

\smallbreak \item[(c)]
$\t s_{3/2}$, $\t s_2$, and $\t s_{5/2}$ are given by
$\t s_{3/2}\ =\ e_{3/2} - (e_0-1)^{-1} e_{1/2}^3$,
\begin{eqnarray*}
   \t s_2 &=& e_2 - {\ts\frac{3}{2}}(2e_0-1)^{-1} (e_{1/2} e_{3/2} + e_{1/2}^4), \\[4pt]
   \t s_{5/2} &=& e_{5/2} -2(2e_0-1)^{-1}\bigr(e_{1/2}e_2 + 
   {\ts\frac{3}{2}} e_{1/2}^2 e_{3/2} + 3(2e_0-3)^{-1} e_{1/2}^5 \bigr)
\end{eqnarray*}

\end{enumerate}
\end{lemma}

\subsection{\bf Decompositions of tensor products.} \label{TPs}
Here we recall the Clebsch-Gordan rules for $\osp_{1|2}$ giving the decompositions of the tensor products of the $\bL(\lambda)$, and state their analogs for the $\bF(\lambda)$.  We begin with two utilitarian lemmas; the proof of the first is left to the reader.

Let us say that a \r\ of $\ds$ is {\em good\/} if $e_0$ acts semisimply with finite dimensional weight spaces and weights bounded below, and $e_{-1/2}$ acts surjectively.

\begin{lemma} \label{LWVs}
Suppose that\/ $V$ is a \r\ of $\ds$ containing a $\bZ_2$-homogeneous \lwv\ $v$ of weight $\lambda$.  If $\lambda\not\in-\bN/2$, then $v$ generates a copy of\/ $\bF(\lambda)^{|v|\Pi}$.  If\/ $\dim V <\infty$, then $\lambda\in-\bN/2$ and $v$ generates a copy of\/ $\bL(-\lambda)^{|v|\Pi}$.  
\end{lemma}

\begin{lemma} \label{good}

\begin{enumerate}

\item[(a)]  Tensor products and $e_0$-split extensions (see Section~\ref{LengthN}) of good \r s are good.

\smallbreak
\item[(b)]  Up to (not necessarily even) equivalence, $\bF(\lambda)$ is the unique good \r\ with simple $e_0$ spectrum $\lambda+\bN$.

\smallbreak
\item[(c)]  $\bD^k(\lambda,p)$ is good for all $\lambda$, $p$, and~$k$.

\smallbreak
\item[(d)]  Fix\/ $T\in\bH(\lambda,p)_\mu$. If $\sigma_{\lambda,p}(e_{-1/2})(T)\in\bD^k(\lambda,p)$, then $T\in\bD^k(\lambda,p)_\mu$ unless $p-\mu \in k + \bZ^+/2$, when $T \in \bC \omega^p \o D^{2(p-\mu)} \oplus \bD^k(\lambda, p)_\mu$.

\smallbreak
\item[(e)]  Non-zero $\dt$-maps from $\bF(\mu)$ to $\bH(\lambda,p)$ exist only when $p-\mu\in\bN/2$, in which case they have image in $\bD^{p-\mu}(\lambda,p)$.

\end{enumerate}
\end{lemma}

\meno{\em Proof.\/}
For~(a), let $V$ and $W$ be good \r s.  Choose a flag $\{V_i\}$ for $V$ such that $V_0 = 0$ and $e_{-1/2}V_i\subset V_{i-1}$, and show that $V_i\ot W\subset e_{-1/2}(V\ot W)$ by induction on~$i$.  We leave the rest of ~(a) and~(b) to the reader.  For~(c), recall that $\bD^k$ is an extension of \tdm s and use~(a) and~(b).  For~(d), use~(c) to choose $S\in\bD^k(\lambda,p)_\mu$ such that $T-S$ is a \lwv\ and then apply Lemma~\ref{Bol}a.  Part~(e) follows from~(d).  $\Box$

\medbreak
To prove the following results on decompositions of tensor products under $\dt$ and $\ds$, it suffices to compute weight space dimensions and use Lemma~\ref{LWVs}, except for the $\ds$-decompositions at $\lambda=\mu\in-\bN/2$ in Proposition~\ref{Fdecomp}.  There one decomposes $\S_s^2\bF(\lambda)$ and $\Lambda_s^2\bF(\lambda)$ along the eigenspaces of the Casimir operator $Q_\ds$, and then applies Lemma~\ref{good}.

\begin{prop} \label{Ldecomp}
$\bL(\lambda)\ot\bL(\mu) \scong \bigoplus_{j=0}^{4\mu} \bL(\lambda+\mu-j/2)^{j\Pi}$ for all $\lambda\ge\mu$ in $\bN/2$.  When $\lambda=\mu$, the summands with $j\equiv 0$ or~$-1$ modulo~4 make up $\S^2_s\bF(\lambda)$, and the others make up $\Lambda^2_s\bF(\lambda)$.
\end{prop}

\begin{prop} \label{Fdecomp}
$\bF(\lambda)\ot\bF(\mu) \tcong \bigoplus_{j=0}^\infty \bF(\lambda+\mu+j/2)^{j\Pi}$ for all $\lambda$ and $\mu$ in $\bC$.  When $\lambda=\mu$, the summands with $j\equiv 0$ or~$-1$ modulo~4 make up $\S^2_s\bF(\lambda)$, and the others make up $\Lambda^2_s\bF(\lambda)$.

For $\lambda+\mu\not\in-\bN/2$ or $\lambda=\mu\in-\bN/2$, these are $\ds$-decompositions.
\end{prop}

\begin{lemma} \label{SymF}
$\S^n_s \bF(\lambda) \tcong \bigoplus_{j\in\bN}\, m^\S_j(n)\, \bF(n\lambda + j/2)^{j\Pi}$ for all $\lambda\in\bC$.  When $\lambda\not\in -\bN/2$, this is an $\ds$-decomposition.

The multiplicity $m^\S_j(n)$ is the number of ways to write\/~$j$ as a sum $j_1+\cdots+j_n$ such that $j_i\in\bN$, $j_1\le j_2\le \cdots\le j_n$, $j_i<j_{i+1}$ for $j_i$ odd, and $j_n-j_{n-1}$ is $0$ for $j_{n-1}$ even and $1$ for $j_{n-1}$ odd.  

In particular, for $n\ge 2$, $m^\S_0=m^\S_3=1$ and $m^\S_1=m^\S_2=0$.  For $n=3$,
$$ \S_s^3 \bF(\lambda)\ \tcong\ \ts\bigoplus_{i,\, j\in\bN,\ 
   b\in\{0,\, \frac{3}{2},\, \frac{5}{2},\, 4\}}\, 
   \bF(3\lambda + b+2j+3i)^{2b\Pi}. $$
\end{lemma}

\begin{lemma} \label{WedgeF}
$\Lambda^n_s\bF(\lambda) \tcong \bigoplus_{j\in\bN}\, m^\Lambda_j(n)\, \bF(n\lambda + j/2)^{j\Pi}$ for all $\lambda\in\bC$.  When $\lambda\not\in -\bN/2$, this is an $\ds$-decomposition.

The multiplicity $m^\Lambda_j(n)$ is the number of ways to write\/~$j$ as a sum $j_1+\cdots+j_n$ such that $j_i\in\bN$, $j_1\le j_2\le \cdots\le j_n$, $j_i<j_{i+1}$ for $j_i$ even, and $j_n-j_{n-1}$ is $0$ for $j_{n-1}$ odd and $1$ for $j_{n-1}$ even.  

In particular, for all $n$, $m^\Lambda_j=0$ for $j<n-1$ and $j=n+1$.  For $n\ge 3$, $m_{n-1}^\Lambda=m^\Lambda_n=m^\Lambda_{n+2}=1$, and $m^\Lambda_{n+3}=2$.
\end{lemma}

We remark that it is possible to classify the good \r s: they are direct sums of injective modules.  It is also possible to describe $\bF(\lambda)\ot \bF(\mu)$ for all $\lambda$ and $\mu$: injective modules arise as summands.  One approach is to generalize Appendix~A of \cite{CM07} from $\dsl_2$ to $\osp_{1|2}$.

The following proposition introduces the {\em supertransvectants\/} $\bJ^{\lambda,\mu}_k$.  It is a corollary of Proposition~\ref{Fdecomp}.  Explicit formulas for the supertransvectants may be found in \cite{GT93}, and in the notation of this paper, in \cite{GO07}.  In fact, the first paragraph of Proposition~\ref{ST} holds even if $\lambda+\mu\in -\bN/2$, provided that either $k\le -(\lambda+\mu)$ or $k\ge 1-2(\lambda+\mu)$.

\begin{prop} \label{ST}
For $\lambda+\mu\not\in-\bN/2$, $k\in\bN/2$, there is an $\ds$-surjection $\bJ^{\lambda,\mu}_k$ from $\bF(\lambda)\ot\bF(\mu)$ to $\bF(\lambda+\mu+k)$ of parity $(-1)^{2k}$, unique up to a scalar.

$\bJ^{\lambda,\mu}_k$ also exists for $\lambda=\mu\in-\bN/2$, unique up to a scalar if stipulated to be supersymmetric for $2k\equiv 0$ or $-1$ modulo~4 and superalternating otherwise.
\end{prop}

\section{Cohomology} \label{GenCoho}
\setcounter{lemma}{0}

In this section we derive some elementary results on $\K$-cohomology which will be used later.  We begin with a review of \lsa\ cohomology and its role in describing \r s of finite length.  It is defined by applying the rule of signs to Lie algebra cohomology.

\subsection{Lie superalgebra cohomology} \label{Coho}
Fix \lsa s $\dh\subset\dg$ and $\dg$-modules $U$, $V$, and $W$.  Write $\pi$ for the action of $\dg$ on $V$.  The space of {\em $\dh$-relative $n$-cochains of $\dg$ with values in $V$\/} is
$$ C^n(\dg,\dh,V) := \Hom_\dh\bigl(\Lambda_s^n(\dg/\dh), V \bigr). $$
If $\dh$ is omitted on the left hand side, it is assumed to be~$\{0\}$.

The {\em coboundary operator\/} $\partial:C^n(\dg,\dh,V)\to C^{n+1}(\dg,\dh,V)$ is an even map of square zero defining the cohomology complex: for $\phi\in C^n(\dg,\dh,V)$,
\begin{displaymath} \begin{array}{rcl}
   \partial\phi\bigl(\bigwedge_0^n X_i\bigr) &:=&
   \sum_{i=0}^n (-1)^{\{i + |X_i|(|\phi| + \sum_0^{i-1} |X_j|)\}}\,
   \pi(X_i) \phi\bigl\{\bigwedge_{\h i} X_j\bigr\} \\[6pt]
   && -\, (-1)^i \phi \bigl\{(\bigwedge_0^{i-1} X_j) \wedge
   \ad(X_i)(\bigwedge_{i+1}^n X_j)\bigr\}.
\end{array} \end{displaymath}

Write $\t\pi$ for the natural $\dg$-action on $C^n(\dg,V)$, and $i(X):C^n\to C^{n-1}$ for the contraction $i(X)\phi(\omega):=(-1)^{|\phi||X|}\phi(X\wedge\omega)$ for $\omega\in\Lambda^{n-1}_s\dg$.  Then
\begin{equation} \label{cohoprops}
   i(X)\circ\partial + \partial\circ i(X) = \t\pi(X), \quad
   \bigl[\t\pi(X), i(Y)\bigr] = i\bigl([X,Y]\bigr),
\end{equation}
and $\partial$ is a $\dg$-map.  Note that $C^n(\dg,\dh, V)$ may be viewed as the subspace of $C^n(\dg,V)$ annihilated by both $i(\dh)$ and $\t\pi(\dh)$.

The {\em $\dh$-relative Ext groups from $V$ to $W$\/} are the cohomology groups
$$ \Ext^n_{\dg,\dh}(V,W) := H^n\bigl(\dg,\dh,\Hom(V,W)\bigr). $$

The {\em cup product associated to composition,\/}
$$ \cup: C^n\bigl(\dg,\dh,\Hom(V,W)\bigr) \ot C^m\bigl(\dg,\dh,\Hom(U,V)\bigr)
   \to C^{n+m}\bigl(\dg,\dh,\Hom(U,W)\bigr), $$
is compatible with $\partial$ (and the $\dg$-actions when $\dh=0$) and drops to a cup product on the Ext groups.  On 1-cochains $\alpha$, $\beta$ it is
$$ \alpha \cup \beta(X\wedge Y) := (-1)^{|X||\beta|}\alpha(X)\circ\beta(Y) 
   - (-1)^{|Y|(|X|+|\beta|)} \alpha(Y)\circ\beta(X). $$

\meno {\bf Remark.}
Suppose that $A$ is a supercommutative superalgebra, $\pi$ is a \r\ of $\dg$ on $A$ acting by superderivations, and $\theta$ is an even $A$-valued 1-cocycle of $\dg$.  Then there is a 1-parameter family $\pi_\lambda$ of \r s of $\dg$ on $A$, defined by $\pi_\lambda(X)a := \pi(X)a + \lambda\theta(X)a$.

\subsection{Representations of finite length} \label{LengthN}
Fix \lsa s $\dh\subset\dg$ and $\dg$-modules $(\psi_r, V_r)$ for $1\le r\le n$.  As in \cite{Co05}, we define
$$ \mbox{\em an extension of\/\ } V_1\to V_2\to \cdots\to V_n $$
to be a $\dg$-module $(\psi,W)$ such that $W$ admits a $\dg$-invariant flag $W=W_1\supset W_2 \supset\cdots\supset W_{n+1}=0$, split in the category of vector spaces, such that $W_r/W_{r+1}$ is $\dg$-equivalent to $V_r$ for all~$r$.  In particular, when the $V_r$ are \irr\ $W$ is of Jordan-H\"older length~$n$.  The dual of an extension of $V_1\to\cdots\to V_n$ is an extension of $V_n^*\to\cdots\to V_1^*$.

Up to equivalence, an extension of $V_1\to\cdots\to V_n$ is a \r\ $\psi$ on $\bigoplus_1^n V_r$ whose block matrix $\psi_{rs}:\dg\to\Hom(V_s, V_r)$ is lower triangular with diagonal entries $\psi_{rr}=\psi_r$.  Regarded as 1-cochains, the subdiagonal entries satisfy the {\em cup equation:\/}
$$ \partial \psi_{rs} + \sum_{r>k>s} \psi_{rk} \cup \psi_{ks} = 0. $$
The invariant flag $\{W_i\}$ is $\dh$-split \iff\ all of the subdiagonal 1-cochains are $\dh$-relative.  The cup equation implies that the entries $\psi_{s+1,s}$ on the first subdiagonal are 1-cocycles.

Note that the $\psi_{rs}$ are necessarily even, because \r s are even by definition.  Odd $\Hom(V_s,V_r)$-valued 1-cochains arise in the construction of extensions in which the parities of some of the $V_i$ are reversed.  For example, if $\alpha$ is an odd $\Hom(V_1,V_2)$-valued 1-cocycle, then $\psi_{21} = L(\ep_{V_2})\circ\alpha$ defines an extension of $V_1\to V_2^\Pi$ (see Lemma~\ref{length2}).

A \r\ is said to be {\em uniserial,\/} or {\em \cind,\/} if all of its \sq s are \ind.  If $V_1,\ldots,V_n$ are \irr, then an extension $\psi$ of $V_1\to\cdots\to V_n$ is uniserial \iff\ each 1-cohomology class $[\psi_{s+1,s}]$ is non-trivial.  Dualization preserves both indecomposability and uniseriality.  In length~2 all \ind\ modules are uniserial and the situation is completely resolved by the following well-known lemma.  We also state a lemma concerning extensions of three modules, which follows from the cup equation.

\begin{lemma} \label{length2}
For $V_1$, $V_2$ \irr, there are bijections from the projective spaces $\bP\Ext^1_{\dg,\dh}(V_1,V_2)_\even$ and $\bP\Ext^1_{\dg,\dh}(V_1,V_2)_\odd$ to the sets of equivalence classes of uniserial $\dh$-split extensions of $V_1\to V_2$ and $V_1\to V_2^\Pi$, respectively.
\end{lemma}

\begin{lemma} \label{length3}
Suppose that $\alpha\in\Ext^1_{\dg,\dh}(V_1,V_2)_\even$ and $\beta\in\Ext^1_{\dg,\dh}(V_2,V_3)_\even$.  Then there exists an $\dh$-split extension $\psi$ of $V_1\to V_2\to V_3$ with $\psi_{21}=\alpha$ and $\psi_{32}=\beta$ \iff\ $[\beta\cup\alpha]$ is trivial in $\Ext^2_{\dg,\dh}(V_1,V_3)$.  In this case the extension is unique (up to equivalence) \iff\ $\Ext^1_{\dg,\dh}(V_1,V_3)_\even = 0$.
\end{lemma}

\subsection{Cohomology of $\K$} \label{PreKCoho}

Here we give some general results concerning the cohomology of $\K$.  Lemma~\ref{IC} reduces the 1-cohomology of $\Hom(V,W)$ to the $\ds$-relative case when $V$ and $W$ have distinct infinitesimal characters under $\ds$.  It may be proven by mimicking the proof of Lemma~3.1 of \cite{Co01}, so we will not give a proof here.  We remark that it is false for 2-cochains; see Proposition~\ref{F012coho}d.

Lemmas~\ref{FCL0} and~\ref{FCL1} are basic results permitting reduction to $\du$-relative $\dt$-invariant cohomology in many settings.  (We expect that the $\du$ statement in Lemma~\ref{FCL1} is true for all~$n$.)  Lemma~\ref{srelC1} describes $\ds$-relative 1-cochains.

\begin{lemma} \label{IC}
Let $V$ and $W$ be $\K$-modules on which $Q_\ds$ acts by distinct eigenvalues (see Section~\ref{Casimir}).  Then any $\Hom(V,W)$-valued 1-cochain with $\ds$-relative coboundary is cohomologous to a unique $\ds$-relative cochain.  Thus
$$ H^0\bigl(\ds, \Hom(V,W)\bigr) = 0, \quad
   H^1\bigl(\K, \Hom(V,W)\bigr) = H^1\bigl(\K, \ds, \Hom(V,W)\bigr). $$
\end{lemma}

\begin{lemma} \label{FCL0}
Let $(\pi,V)$ be a $\K$-module on which $e_0$ acts semisimply.  Suppose that $\phi$ is a $V$-valued $n$-cochain of $\K$ such that $\partial\phi$ is $e_0$-relative.  Then $\phi$ is cohomologous to a $\t\pi(e_0)$-invariant cochain $\phi_0$ whose contraction by $e_0$ is an $e_0$-relative cocycle.  If $\phi$ is $e_{-1}$- or $\du$-relative, then so is $\phi_0$.
\end{lemma}

\meno{\em Proof.\/}
Since $e_0$ acts semisimply on $\K$ and $V$, there exist unique (locally finite) $n$-cochains $\phi_\lambda$ of weight~$\lambda$ under $\t\pi(e_0)$ such that $\phi = \sum_\lambda \phi_\lambda$.  Verify that each $\partial\phi_\lambda$ is $e_0$-relative, and hence by~(\ref{cohoprops}) that $\partial i(e_0)\phi_\lambda = \lambda\phi_\lambda$.  

Define $T:=\sum_{\lambda\not=0} i(e_0)\phi_\lambda/\lambda$.  Then $\phi-\partial T = \phi_0$ is of weight~$0$, \ie\ $e_0$-invariant, and $i(e_0)\phi_0$ is an $e_0$-relative $(n-1)$-cocycle.  We leave the preservation of relativity to the reader.  $\Box$

\begin{lemma} \label{FCL1}
Let $(\pi,V)$ be a $\K$-module on which $e_{-1/2}$ acts surjectively.  Suppose that $\phi$ is a $V$-valued $n$-cochain of $\K$ such that $\partial\phi$ is $\du$-relative.  Then $\phi$ is cohomologous to an $e_{-1}$-relative cochain, and for $n\le 2$ it is cohomologous to a $\du$-relative cochain.
\end{lemma}

\meno{\em Proof.\/}
We will find an $(n-1)$-cochain $\delta$ such that $i(e_{-1})\partial\delta = i(e_{-1})\phi$.  Then by~(\ref{cohoprops}), $\phi-\partial\delta$ is $e_{-1}$-relative.  Assume that $i(e_{-1})\delta=0$; then we need $\t\pi(e_{-1})\delta = i(e_{-1})\phi$, \ie
$$ \pi(e_{-1}) \bigl(\delta(\omega)\bigr) = 
   \delta\bigl(\ad(e_{-1})\omega\bigr) + \phi(e_{-1} \wedge \omega)
   \mbox{\rm\ for all\ } \omega\in\Lambda^{n-1}_s\K. $$
Since $i(e_{-1})$ annihilates both sides and $\pi(e_{-1})=\pi(e_{-1/2})^2$ is surjective on $V$, we can solve for $\delta$ by induction upward on the weight of $\omega$.

For the $\du$-relativity statement we must find $\delta$ such that $i(e_{-1/2})\partial\delta = i(e_{-1/2})\phi$.  For $n=1$ it suffices to choose $\delta\in V$ such that $\pi(e_{-1/2})\delta = (-1)^{|\phi|}\phi(e_{-1/2})$.  

The $n = 2$ case is complicated by the fact that $i(e_{-1/2})^2\not=0$.  We may assume that $\phi$ is $e_{-1}$-relative, and we will also assume that $i(e_{-1/2})\delta = 0$ (the latter assumption does not seem to suffice for $n\ge 3$, although $i(e_{-1})\delta=0$ may).  The equation to be solved is then $\t\pi(e_{-1/2})\delta = i(e_{-1/2})\phi$, \ie
$$ (-1)^{|\phi|} \pi(e_{-1/2}) \bigl(\delta(e_m)\bigr) = \delta\bigl([e_{-1/2},e_m]\bigr)
   + \phi(e_{-1/2} \wedge e_m) \mbox{\rm\ for all\ } e_m\in\K. $$

At $m=-1$ this is $\pi(e_{-1/2})\delta(e_{-1})=0$, and at $m=-1/2$ it is $2\delta(e_{-1}) = -\phi(e_{-1/2} \wedge e_{-1/2})$.  In order to verify that these are compatible, we must use the fact that $\partial\phi$ is $\du$-relative: $\partial\phi(\bigwedge^3 e_{-1/2})=0$ yields $\pi(e_{-1/2})\phi(\bigwedge^2 e_{-1/2})=0$.  For $m\ge 0$, we can use the surjectivity of $\pi(e_{-1/2})$ to solve for $\delta(e_m)$ by induction upward on the weight of $\omega$.  $\Box$

\begin{lemma} \label{srelC1}
Fix a \r\ $\pi$ of $\K$ on a space $V$.

\begin{enumerate}

\item[(a)] There is an odd bijection $V^\du_{3/2}\to C^1(\K,\ds,V)$, $v\mapsto\phi$:
\begin{equation} \label{srelC1eqnA}
   \phi(e_n) := {\ts\frac{2^{-2\{n-1/2\}}}{\lf n-3/2 \rf!}} (-1)^{(2n-1)(|v|+1)}
   \pi(e_{1/2})^{2n-3} v \mbox{\rm\ \ for\ } n\ge 3/2.
\end{equation}

\smallbreak
\item[(b)] There is an odd bijection $\Hom_\ds\bigl(\bF(3/2),V\bigr)\to C^1(\K,\ds,V)$, $\Phi\mapsto\phi$:
\begin{equation} \label{srelC1eqnB}
   \phi := \Phi\circ\SBol_{5/2}(-1)\circ\bX^{-1}, \quad
   \phi\bigl(\bX(F)\bigr) = \Phi\bigl(\omega^{3/2}\o DF''\bigr).
\end{equation}

\end{enumerate}
\end{lemma}

\meno {\it Proof.\/}
The inverse of the bijection in~(a) is $\phi\mapsto\phi(e_{3/2})$.  The point is that an $\ds$-relative 1-cochain $\phi$ is determined by its value on $e_{3/2}$: for $n\ge 3/2$,
$$ \phi\bigl(\ad(e_{1/2})^{2n-3}(e_{3/2})\bigr) =
   (-1)^{(2n-1)|\phi|} \pi(e_{1/2})^{2n-3} \phi(e_{3/2}). $$
To finish the proof of~(a), verify that
\begin{equation} \label{Brsfact1}
   \ad(e_{1/2})^{2n-3} (e_{3/2}) = 2^{2\{n-1/2\}}\lf n-3/2 \rf! \th e_n.
\end{equation}

For~(b), use Lemmas~\ref{Xomega} and~\ref{Bol}c to prove that any $\ds$-relative 1-cochain must factor through $\SBol_{5/2}\circ\bX^{-1}$.  We leave the rest to the reader.  $\Box$

\subsection{$\bF(\lambda)$-valued cohomology} \label{FCoho}

In this section we compute the $\bF(\lambda)$-valued 0- and 1-cohomology of $\K$, and some of the 2-cohomology.  The results for $H^0$ and $H^1$ are known \cite{AB04}.  Define $\bF(0)$-, $\bF(1/2)$- and $\bF(3/2)$-valued 1-cochains by
$$ \theta := \SBol_1(-1)\circ\bX^{-1}, \
   \alpha := \SBol_{3/2}(-1)\circ\bX^{-1}, \
   \beta  := \SBol_{5/2}(-1)\circ\bX^{-1}. $$

\begin{prop} \label{F012coho}

\begin{enumerate}

\item[(a)] $H^0\bigl(\K,\bF(\lambda)\bigr)$ is $\bC 1$ for $\lambda=0$, and~$0$ otherwise.

\smallbreak
\item[(b)] The 1-cochains $\theta$, $\alpha$, and $\beta$ are all non-trivial 1-cocycles: $\theta$ is even, $\du$-relative, and $\dt$-invariant; $\alpha$ is odd and $\dt$-relative; and $\beta$ is odd and $\ds$-relative.  Moreover, $[\theta]$ is not $\dt$-relative and $[\alpha]$ is not $\ds$-relative.

\smallbreak
\item[(c)] $H^1\bigl(\K,\bF(\lambda)\bigr)$ is spanned by $[\theta]$ for $\lambda=0$, $[\alpha]$ for $\lambda=1/2$, $[\beta]$ for $\lambda=3/2$, and is~$0$ otherwise.  

\smallbreak
\item[(d)] The 2-cocycles $\theta\cup\alpha$, $\theta\cup\beta$, and $\beta\cup\beta$ are non-trivial, and $\theta\cup\theta$, $\alpha\cup\alpha$, and $\alpha\cup\beta$ are trivial.  Neither $[\theta\cup\alpha]$ nor $[\theta\cup\beta]$ is $\dt$-relative.

\smallbreak
\item[(e)] $H^2\bigl(\K,\bF(\lambda)\bigr)$ is spanned by $[\theta\cup\alpha]$ for $\lambda=1/2$, $[\theta\cup\beta]$ for $\lambda=3/2$, and $[\beta\cup\beta]$ for $\lambda=3$.  It is at least 1-dimensional for $\lambda=5$.  Otherwise it is~$0$, except possibly for $2\lambda\in 7+\bN$.

\smallbreak
\item[(f)] $\dim H^2\bigl(\K,\ds,\bF(\lambda)\bigr) = 1$ for $\lambda=3$ and~$5$ (at $\lambda=3$ it is spanned by $[\beta\cup\beta]$).  Otherwise it is~$0$, except possibly for $\lambda\in \{13/2, 7\} + 2\bN$.

\end{enumerate}
\end{prop}

\meno{\it Proof.\/}  Part~(a) is Lemma~\ref{F}b.  For the sequel we will need several facts.  First, by Lemma~\ref{Bol} there are even equivalences
\begin{equation} \label{Kmodsubalg}
   \K/\du\tcong\bF(0), \quad \K/\dt\tcong\bF(1/2)^\Pi, \quad \K/\ds\scong\bF(3/2)^\Pi.
\end{equation}
Hence by Lemmas~\ref{SymF}, \ref{WedgeF}, and~\ref{taep},
\begin{eqnarray}
   \label{wedge2u}
   \Lambda_s^2 (\K/\du) &\tcong& \bF(1/2)^\Pi \oplus \bF(1)
   \oplus \bF(5/2)^\Pi \oplus \bF(3) \oplus\cdots, \\[4pt]
   \label{wedge2t}
   \Lambda_s^2 (\K/\dt) &\tcong& \bF(1) \oplus \bF(5/2)^\Pi 
   \oplus \bF(3) \oplus \bF(9/2)^\Pi \oplus\cdots, \\[4pt]
   \label{wedge2s}
   \Lambda_s^2 (\K/\ds) &\scong& \bF(3) \oplus \bF(9/2)^\Pi 
   \oplus \bF(5) \oplus \bF(13/2)^\Pi \oplus\cdots, \\[4pt]
   \label{wedge3s}
   \Lambda_s^3 (\K/\ds) &\scong& \bF(9/2)^\Pi \oplus \bF(6)
   \oplus \bF(13/2)^\Pi \oplus\cdots.
\end{eqnarray}
Again by Lemma~\ref{Bol}, there exists a $\dt$-map from $\bF(\mu)$ to $\bF(\nu)$ \iff\ $\nu-\mu\in\bN/2$, in which case such maps are spanned by $\SBol_{\nu-\mu}(\mu)$.

For~(b), use Lemmas~\ref{Xomega} and~\ref{Bol} to verify that $\theta$ is $\du$-relative and $\dt$-invariant, $\alpha$ is $\dt$-relative, and $\beta$ is $\ds$-relative.  Therefore $\partial\theta$, $\partial\alpha$, and $\partial\beta$ have the same properties, respectively.  But by~(\ref{wedge2u}), (\ref{wedge2t}), and~(\ref{wedge2s}) there are no such 2-cochains, so they are zero.  By Lemma~\ref{F}, $B^1\bigl(\K,\bF(\lambda)\bigr)$ has no $\du$-relative $\dt$-invariant elements, implying the rest of~(b).

For~(c), suppose that $\phi\in Z^1\bigl(\K,\bF(\lambda)\bigr)$ is non-zero.  If it is $\ds$-relative, Lemma~\ref{srelC1}b gives $\lambda=3/2$ and $\phi\propto\beta$.  If it is $\dt$- but not $\ds$-relative, $\phi(e_{1/2})$ must be a non-zero \lwv\ of weight~$1/2$, so $\lambda=1/2$.  By Lemma~\ref{Bol}, $\phi$'s $\dt$-invariance implies $\phi\propto\alpha$.

If $\phi$ is not $\dt$-relative, we may use Lemmas~\ref{FCL0} and~\ref{FCL1} to assume that it is $\du$-relative, $\dt$-invariant, and satisfies $\phi(e_0)$ $\K$-invariant.  This implies $\lambda=0$ and $\phi\propto\theta$.

For~(d), $\theta\cup\theta=0$, and $[\alpha\cup\alpha]=[\alpha\cup\beta]=0$ follows from the proof of~(e).  If $\theta\cup\alpha$ or $\theta\cup\beta$ is a coboundary $\partial\delta$, we may assume by Lemmas~\ref{FCL0} and~\ref{FCL1} again that $\delta$ is $\dt$-invariant.  But this would force $\delta\propto\alpha$ or $\beta$, respectively, whence $\partial\delta=0$.  Thus $\theta\cup\alpha$ and $\theta\cup\beta$ are non-trivial and are their classes' only $\dt$-invariant representatives.  Since $i(e_0)(\theta\cup\alpha) = \alpha$ and $i(e_0)(\theta\cup\beta) = \beta$, neither class is $\dt$-relative.

For~(e), suppose that $\phi\in Z^2\bigl(\K,\bF(\lambda)\bigr)$ is non-zero.  Use Lemmas~\ref{FCL0} and~\ref{FCL1} one more time to assume that $\phi$ is $\du$-relative and $\dt$-invariant.  Then it factors through $\Lambda^2_s(\K/\du)$, so by~(\ref{wedge2u}) it is zero unless $\lambda\in\bZ^+/2$.  Suppose that $\phi$ is not $\dt$-relative.  Then $i(e_0)\phi$ is a non-zero $\dt$-relative 1-cocycle, necessarily non-trivial as there are no $\dt$-relative 1-coboundaries, so we must have $\lambda=1/2$ or~$3/2$.  Moreover, subtracting a multiple of $\theta\cup\alpha$ or $\theta\cup\beta$ from $\phi$ gives a $\dt$-relative 2-cocycle.

By~(\ref{wedge2t}), $\bF(\lambda)$-valued $\dt$-relative 2-cochains exist only for $\lambda\in 1+\bN/2$: $C^2\bigl(\K,\dt,\bF(\lambda)\bigr)$ is 1-dimensional for $\lambda = 1$, $3/2$, and~$2$; 2-dimensional for $\lambda=5/2$; and 3-dimensional for $\lambda=3$.  Except at $\lambda=3/2$, one of the dimensions is accounted for by $\partial\bigl(\SBol_{\lambda+1}(-1)\circ\bX^{-1}\bigr)$, the only $\dt$-invariant 2-coboundary.  The $\bF(3/2)$-valued $\dt$-relative 2-cochain $\phi$ is not a cocycle: $\phi(e_{1/2}\wedge e_1)\propto \omega^{3/2}$ implies $\partial\phi(\bigwedge^3 e_{1/2})\not=0$.

At $\lambda=5/2$, the second dimension is not a 2-cocycle: there is clearly a $\dt$-relative $\bF(5/2)$-valued 2-cochain $\phi$ mapping $e_{1/2}\wedge e_2$ to $\omega^{5/2}$ and $e_1\wedge e_{3/2}$ (and lower weights) to zero, and it is simple to verify $\partial\phi(e_{1/2}\wedge e_{1/2}\wedge e_{3/2})\not=0$.  At $\lambda=3$, a similar easy calculation shows that the second dimension is accounted for by a non-cocycle, and third is $\beta\cup\beta$.

For~(f), use~(\ref{wedge2s}) to verify $\dim C^2\bigl(\K,\ds,\bF(\lambda)\bigr)=1$ for $\lambda\in\{3,\, 9/2\} + 2\bN$.  By~(\ref{wedge3s}), the cochains at $\lambda=3$ and~$5$ are necessarily cocycles.  We must prove that the unique $\ds$-relative 2-cochain $\phi$ at $\lambda=9/2$ is not a cocycle.  Since $\phi$ kills weights below $9/2$, $\partial\phi(\bigwedge^3 e_{3/2})$ is a non-zero multiple of $\phi(e_{3/2}\wedge e_3)$.  To prove that this in turn is a non-zero multiple of $\omega^{9/2}$, it suffices to prove that $e_{3/2}\wedge e_3$ is not in the copy of $\bF(3)$, \ie\ that it is not a multiple of $\ad(e_{1/2})^3 (\bigwedge^2 e_{3/2})$.  This may be done directly.  $\Box$

\meno{\bf Remark.}  Let $F(\lambda)$ denote the $\VR$-module of $\lambda$-densities over $\bR$.  Goncharova proved the following important result: the affine-relative cohomology $H^n\bigl(\VR,\la \px, x\px \ra, F(\lambda)\bigr)$ is $\bC$ for $2\lambda = 3n^2\pm n$, and $0$ otherwise, the only non-projectively relative class being in $H^1 F(1)$ \cite{Gon73}.  (Note that this forces the ring structure on the cohomology of $\bigoplus_\lambda F(\lambda)$ induced by commutative multiplication of densities to be trivial.)

In light of Proposition~\ref{F012coho}, we make the following conjecture (which would imply that $\beta\cup\beta$ is the only non-trivial cup product on $\bigoplus_\lambda\bF(\lambda)$).  Preliminary calculations suggest that it is correct for $n\le 3$.

\meno {\bf Conjecture~1.}  $\dim H^n\bigl(\K,\dt,\bF(\lambda)\bigr) = 1$ for $\lambda = n^2\pm n/2$, and $0$ otherwise.  The only non-$\ds$-relative class is $[\alpha]$.

\section{The conformal symbol and the main result}  \label{CS}
\setcounter{lemma}{0}

We shall say that values of $p$ in $\bZ^+/2$ are\/ {\em resonant\/} and other values are\/ {\em non-resonant\/} (the terms\/ {\em regular\/} and\/ {\em singular\/} are used in the author's earlier papers).  The space $\bD^k(\lambda,p)$ is called\/ {\em resonant\/} if $p$ is resonant and $k\ge p$, and\/ {\em non-resonant\/} otherwise.  Define $k_p:=p-1/2$ if $p$ is a resonant value and $k_p:=\infty$ otherwise, so that $\bD^{k_p}(\lambda,p)$ is maximal non-resonant.  We begin by proving that the conformal symbol exists in the non-resonant case.

\begin{prop} \label{CQCS} {\rm \cite{GMO07}}
There exists a unique $\ds$-equivalence
$$ \CQ_{\lambda,p}: \bigoplus_{j=0}^{2k_p} \bF(p-j/2)^{j\Pi} \to \bD^{k_p}(\lambda,p) $$
such that $\CQ_{\lambda,p}(\omega^{p-j/2}G) \equiv \omega^p G\o D^j$ modulo $\bD^{(j-1)/2}$, the {\em conformal quantization.\/}  It is even and maps $\omega^{p-k}$ to $\omega^p \o D^{2k}$.  We define the {\em conformal symbol\/} $\CS_{\lambda,p}$ to be its inverse.
\end{prop}

\meno{\it Proof.\/}
By Proposition~\ref{fine}, $\bD^{k_p}$ has composition series $\{\bF(p-j/2)^{j\Pi}: 0\le j\le 2k_p\}$.  Hence the result follows from Proposition~\ref{ghost} by considering the eigenspaces of the Casimir operator $Q_\ds$, provided that its eigenvalues on the composition series are distinct.  This is the case \iff\ $p\not\in(2\bN+3)/4$.

For $p\in(2\bN+3)/4$ some of the eigenvalues of $Q_\ds$ are doubled: $\bF(p-k)^{2k\Pi}$ and $\bF(-p+k+1/2)^{2k\Pi}$ have the same eigenvalue for $k=0, 1/2, \ldots, p-3/4$.  Thus the eigenspaces corresponding to the doubled eigenvalues are \r s of $\ds$ with composition series $\{\bF(p-k)^{2k\Pi},\bF(-p+k+1/2)^{2k\Pi}\}$.  Here neither $p-k$ nor $-p+k+1/2$ is in $-\bN/2$, so it follows from Lemmas~\ref{F}d and~\ref{LWVs} that these \r s split uniquely. 

The image of $\omega^{p-k}$ can be computed using the fact that it is a \lwv\ together with Lemma~\ref{Bol}a.  $\Box$

\meno{\bf Remarks.}
The main results of \cite{GMO07} are explicit formulas for $\CS_{\lambda,p}$ and $\CQ_{\lambda,p}$ ($\CQ_{0,0}$ was given earlier in \cite{CMZ97}).  In this paper we avoid using these formulas by focusing on \lwv s.  However, they are interesting in their own right and their form is instructive.  Let us give $\CS_{\lambda,p}$ and point out some of its features.

Before we write the formula, check that the ``na\"\i ve symbol''
$$ \NS_{\lambda,p}: \bD(\lambda,p) \to \bigoplus_{j=0}^\infty \bF(p-j/2)^{j\Pi},
   \quad \omega^p G \o D^j\mapsto \omega^{p-j/2} G $$
is an even $\dt$-equivalence.  Therefore the map
$$ \omega^{p-j/2} G \mapsto \mbox{\rm\ the $\bF(p-i/2)^{i\Pi}$-summand of\ } 
   \CS_{\lambda,p}(\omega^p G \o D^j) $$
from $\bF(p-j/2)^{j\Pi}$ to $\bF(p-i/2)^{i\Pi}$ is in $\bigl(\bH(p- \frac{j}{2}, \frac{j-i}{2} )^{(j-i)\Pi}\bigr)^\dt$.  By Lemma~\ref{Bol}b, it must be in $\bC\omega^{(j-i)/2} D^{j-i}$.  This proves that $\CS_{\lambda,p}$ must take the form
\begin{displaymath}
   \CS_{\lambda,p}(\omega^p G \o D^j) = 
   \sum_{i=0}^j c_{ji}(\lambda,p)\, \omega^{p-i/2} D^{j-i}(G)
\end{displaymath}
for some constants $c_{ji}(\lambda,p)$.  The strategy of \cite{GMO07} is to compute the recursion relation induced on the $c_{ji}$ by $\ds$-invariance:

\begin{thm} \cite{GMO07}
$$ c_{ji} = (-1)^{\lf\frac{j-i+1}{2}\rf} 
   {\lf\frac{j}{2}\rf \choose \bigl\lf\frac{j-i}{2} + \{\frac{i}{2}\}\bigr\rf}
   {2\lambda + \lf\frac{j-1}{2}\rf \choose \bigl\lf\frac{j-i}{2} + \{\frac{i+1}{2}\}\bigr\rf}
   {2p-i-1 \choose \bigl\lf\frac{j-i+1}{2}\bigr\rf}^{-1}. $$
\end{thm}

We conclude these remarks with an example: (\ref{CH}) may be restated as
$$ \bX(\omega^{-1}F) = \CQ_{0,0}(\omega^{-1}F) = {\ts\frac{1}{2}}D(F)\o D - F\o D^2. $$

\meno{\bf Definition.} {\em
Let $\pi(\lambda,p)$ be the \r\ of $\K$ on $\bigoplus_{j=0}^{2k_p} \bF(p-j/2)^{j\Pi}$ obtained by transferring the canonical action $\sigma_{\lambda,p}$ of $\K$ on $\bD^{k_p}(\lambda,p)$ via $\CS_{\lambda,p}$:
$$ \pi(\lambda,p)(X) := \CS_{\lambda,p}\circ\sigma_{\lambda,p}(X)
   \circ\CS_{\lambda,p}^{-1}. $$

We will regard $\pi(\lambda,p)$ as a $(2k_p+1)\times(2k_p+1)$-matrix whose entries are 1-cochains of $\K$.  We choose the indices of the matrix entries to match the degrees of their range and domain, so they lie in $p-\bN/2$:\/}
$$ \pi_{p-i/2,p-j/2}(\lambda,p):\K\to   
   \bH\bigl(p-j/2,(j-i)/2\bigr)^{(j-i)\Pi}. $$

\begin{lemma} \label{pirs}
\begin{enumerate}

\item[(a)]  $\pi(\lambda,p)$ is lower triangular: $\pi_{rs}=0$ for $r<s$.

\smallbreak
\item[(b)]  $\pi_{rr}$ is the action $\pi_r^{2(p-r)\Pi}$ of $\K$ on $\bF(r)^{2(p-r)\Pi}$.

\smallbreak
\item[(c)]  $\pi_{rs}\in C^1\bigl(\K,\ds,\bH(s,r-s)^{2(r-s)\Pi})\bigr)$ for $r>s$.

\end{enumerate}
\end{lemma}

\meno{\it Proof.\/}
Part~(a) is due to the fact that $\sigma_{\lambda,p}$ and $\CS_{\lambda,p}$ preserve the relevant filtrations.  Part~(b) follows from Proposition~\ref{fine}.  Part~(c) follows from the $\ds$-covariance of $\CS_{\lambda,p}$.  $\Box$

\medbreak
As stated in the introduction, our main result is the action of $\K$ on $\bD^{k_p}(\lambda,p)$ in terms of the conformal symbol.  By this we mean the formula for the matrix entries $\pi_{rs}(\lambda,p)$.  It is relatively easy to compute these entries up to a scalar, and to prove that they are zero for $r-s=1/2$ or~$1$; this is done in Corollary~\ref{uptoscalar}.  However, the computation of the scalars is difficult.  In Theorem~\ref{main} we give an unsimplified formula for them, along with its simplification at $r-s=3/2$, $2$, and~$5/2$, the values needed in our applications.  First we define certain even $\ds$-relative 1-cochains of $\K$ and prove that the matrix entries $\pi_{rs}$ are multiples of them.

\begin{lemma} \label{beta}
For $2p\in 3+\bN$, $C^1\bigl(\K,\ds,\bH(\lambda,p)\bigr)$ is spanned by the $\ds$-relative 1-cochain $\o\beta_p(\lambda)$ of parity $(-1)^{2p}$ defined by~(\ref{srelC1eqnA}) with $v = 4\omega^p \o D^{2p-3}$.  Otherwise it is zero.

For $2p\in 3+\bN$, $C^1\bigl(\K,\ds,\bH(\lambda,p)^{2p\Pi}\bigr)$ is spanned by the even $\ds$-relative 1-cochain $\beta_p(\lambda)$ defined as $R(\ep_{\bF(\lambda)}^{2p}) \circ \o\beta_p(\lambda)$.  Otherwise it is zero.
\end{lemma}

\meno{\em Proof.\/}  The first paragraph is a corollary of Lemmas~\ref{Bol}a and~\ref{srelC1}a.  For the second, recall that $R(\ep)$ intertwines $\bH(\lambda,p)$ and $\bH(\lambda,p)^\Pi$.  $\Box$

\meno{\bf Remark.}
For $2p\in 3+\bN$, the restriction of $\CQ_{\lambda,p}\circ\ep^{2p+1}$ to $\bF(3/2)$ is an $\ds$-injection to $\bH(\lambda,p)$, and $\o\beta_p(\lambda)$ is the cochain associated to it by~(\ref{srelC1eqnB}).  In terms of the $\bF(3/2)$-valued 1-cocycle $\beta$, $\o\beta_p = \CQ_{\lambda,p}\circ\ep_{\bF(3/2)}^{2p+1}\circ\beta$.  Viewed as a map from $\K\ot\bF(\lambda)$ to $\bF(\lambda+p)$, $\o\beta_p \propto \bJ^{3/2,\lambda}_{p-3/2} \circ (\beta\ot 1)$.

The following formula for $\o\beta_p(\lambda)$ will be useful:
\begin{equation} \label{obeta}
   \o\beta_p(\lambda)\bigl(\ad(e_{1/2})^n e_{3/2}\bigr) = 4 (-1)^{2pn}
   \sigma_{\lambda,p}(e_{1/2})^n \bigl(\omega^p\o D^{2p-3}\bigr).
\end{equation}

The next corollary follows immediately from Lemmas~\ref{pirs}c and~\ref{beta}.  In the subsequent definition the superscript~$m$ is redundant until we consider \psdog s.  The final theorem of the section is the main result of the paper; it will be proven in Section~\ref{Proofs}.

Henceforth we suppress the argument $(\lambda,p)$ whenever convenient, and we use the notation
$$ c(\lambda,p) := \lambda+p/2-1/4. $$

\begin{cor} \label{uptoscalar}
For $r-s= 1/2$ or~$1$ (on the first two subdiagonals), $\pi_{rs}$ is zero.  For $r-s\ge 3/2$ (on the lower subdiagonals), it is a multiple of $\beta_{r-s}(s)$.
\end{cor}

\meno{\bf Definition.} {\em
For $r\in p-\bN/2$, $s\ge p-k_p$ (for non-resonance), $r-s\in 1+\bZ^+/2$, and $m\equiv 2(p-s)$ modulo $\bZ_2$, define scalars $b^m_{rs}(\lambda,p)$ by}
$$ \pi_{rs}(\lambda,p) := b^{2(p-s) \mod 2}_{rs}(\lambda,p) \beta_{r-s}(s). $$

\begin{thm} \label{main}
The scalars $b^m_{rs}$ are polynomial in~$\lambda$ and~$p$ and rational in~$s$.  They are given by the formula $b^m_{rs} = P^m_{rs}/B_{rs}$, where $B_{rs}$ is given by~(\ref{Brs}) and $P^m_{rs}$ is given by~(\ref{Pmrs}).  At $r-s= \frac{3}{2}$, $2$, and~$\frac{5}{2}$,
\begin{eqnarray*}
   b^0_{s+3/2,s} &=& -c{p-s\choose 2}\ \frac{1}{s+1/2}, \\
   b^1_{s+3/2,s} &=& \Bigl(p-s-1/2\Bigr)\ \frac{4ps+2p+1-16c^2}{16(s+1/2)}, \\
   b^0_{s+2,s} &=& {p-s\choose 2}\ \frac{(2s+3)(4p+2s+1)-48c^2}{32s(s+3/2)}, \\
   b^1_{s+2,s} &=& {p-s-1/2\choose 2}\ \frac{s(2p+s+1)-12c^2}{8s(s+3/2)}, \\
   b^0_{s+5/2,s} &=& -{p-s\choose 3} {s+2\choose 3}^{-1}\ 
      \frac{4(s+1)p+3-48c^2}{64}, \\
   b^1_{s+5/2,s} &=& -c{p-s-1/2\choose 2} {s+2\choose 3}^{-1}\ 
      \frac{4(s+1)p-s^2-2s+2-12c^2}{24}.
\end{eqnarray*}
\end{thm}

\section{Pseudodifferential operators and symmetries}  \label{PSDOG}
\setcounter{lemma}{0}

In this section we generalize our results to \psdog s (\psidog s) of arbitrary complex order, extending the 2-parameter family $\bD(\lambda,p)$ to a 3-parameter family which is useful in applications.  We also use conjugation of \dog s and the super Adler trace to prove certain important symmetries of the $b^m_{rs}(\lambda,p)$.  Our definition of \psidog s is based on the discussion at the end of Section~7 of \cite{CMZ97}:

\meno{\bf Definition.}  {\em
For $z\in\bC$, $m\in\bZ_2$, define formal symbols $\o D^z_m$ of parity~$m$:
$$ \o D^z_0 := e^{i\pi z/2} \partial_x^{z/2}, \quad
   \o D^z_1 := e^{i\pi (z-1)/2} \partial_x^{(z-1)/2} \o D. $$
For $k\in\bC$ and $\ell\in\bZ_2$, the space $\Psi^k_\ell(\lambda,p)$ of\/ {\em \psidog s of degree~$\le (k,\ell)$ from $\bF(\lambda)$ to $\bF(\lambda,p)$\/} consists of formal series:
$$ \Psi^k_\ell(\lambda,p) := \Bigl\{\omega^p \sum_{j\in \bN} G_j\o D^{2k-j}_{\ell-j}:
   \ G_j\in\PRoo\Bigr\}.  $$
Define $\Psi(\lambda,p)$ to be the space of all \psidog s from $\bF(\lambda)$ to $\bF(\lambda+p)$.\/}

\medbreak
Note that elements of $\Psi(\lambda,p)$ are not actually operators from $\bF(\lambda)$ to $\bF(\lambda+p)$.  However, composition of \dog s extends to \psidog s via $[\partial_x^z,\o D]=0$ and the {\em generalized Leibniz rule:\/}
$$ \partial_x^z \circ F := \sum_{j\in\bN} {\ts{z\choose j}} 
   F^{(j)}\circ \partial_x^{z-j} \mbox{\rm\ for all $z\in\bC$.} $$
This permits us to extend the actions $\sigma_{\lambda,p}$ and $\sigma'_{\lambda,p}$ of $\K$ on $\bD(\lambda,p)$ and $\bD(\lambda,p)^\Pi$ defined in Section~\ref{DOGs} to actions on $\Psi(\lambda,p)$ and $\Psi(\lambda,p)^\Pi$, respectively.  Observe that $\o D^{z'}_{m}\circ \o D^z_m = \o D^{z'+z}_{m'+m}$, and for $j\in\bN$, $\o D^j_j =\o D^j$.

Let us now extend some of the results of Sections~\ref{FF}, \ref{FDH}, and~\ref{CS} to \psidog s.  A modification of the proof of Lemma~\ref{Bol} yields Lemma~\ref{psiBol} (prove part~(c) by computing the action of $e_{1/2}$ on $\Psi^k_\ell(\lambda,p)^\dt$ directly).  The proof of Proposition~\ref{fine} gives Lemma~\ref{psifine}.  Then the idea of the proof of Proposition~\ref{CQCS} coupled with Lemma~\ref{psiBol} lead to Lemma~\ref{psiCQCS1}.

\begin{lemma} \label{psiBol}
\begin{enumerate}

\item[(a)] $\Psi^k_\ell(\lambda,p)^\du = \bigl\{\omega^p\sum_{j\in\bN} g_j \o D^{2k-j}_{\ell-j}: g_j\in\bC\bigr\}$.

\smallbreak
\item[(b)] $\Psi^k_\ell(\lambda,p)^\dt = \bC\omega^p\o D^{2p}_{\ell-2k+2p}$ for $k-p\in\bN/2$, and~$0$ otherwise.

\smallbreak
\item[(c)] $\Psi^k_\ell(\lambda,p)^\ds = \bC\omega^p\o D^{2p}_1$ for $k-p\in\bN/2$, $\ell-2k+2p$ odd, and $c(\lambda,p)=0$.  It is $\bC 1$ for $k\in\bN/2$, $\ell-2k$ even, and $p=0$.  Otherwise it is~$0$.

\end{enumerate}
\end{lemma}

\begin{lemma} \label{psifine}
The filtration $\{\Psi^{k-j/2}_{\ell-j}(\lambda,p): j\in\bN\}$ of $\Psi^k_\ell$ is $\sigma_{\lambda,p}$-invariant and compatible with composition.  There is an even equivalence
$$ \Psi^k_\ell/\Psi^{k-1/2}_{\ell-1} \to \bF(p-k)^{\ell\Pi}, \quad
   \omega^p G \o D^{2k}_\ell \mapsto \omega^{p-k} G. $$
\end{lemma}

\begin{lemma} \label{psiCQCS1}
For $k-p\not\in \bN/2$, there is a unique $\ds$-map from $\bF(p-k)^{\ell\Pi}$ to $\Psi^k_\ell(\lambda,p)$ mapping $\omega^{p-k}$ to $\omega^p\o D^{2k}_\ell$.  It is even and takes $\omega^{p-k} G$ to a \psidog\ of symbol $\omega^p G \o D^{2k}_\ell$.
\end{lemma}

In light of Lemma~\ref{psiCQCS1}, it is appropriate to define $\Psi^k_\ell(\lambda,p)$ to be {\em non-resonant\/} if $k-p\not\in\bN/2$.  In order to define the conformal quantization and symbol for \psidog s, let $\prod_{j\in\bN} \bF(p-k+j/2)^{(\ell-j)\Pi}$ denote the set of all formal series with one summand in each factor.  Proposition~\ref{psiCQCS2} below follows from Lemma~\ref{psiCQCS1}.  Then the proof of Theorem~\ref{psimain} follows that of Lemma~\ref{pirs}, Corollary~\ref{uptoscalar}, and Theorem~\ref{main}.

\begin{prop} \label{psiCQCS2}
For $k-p\not\in\bN/2$, the conformal quantization extends uniquely to an even $\ds$-equivalence
$$ \CQ_{\lambda,p}: \prod_{j\in\bN} \bF(p-k+j/2)^{(\ell-j)\Pi}
   \to \Psi^k_\ell(\lambda,p) $$
such that $\CQ_{\lambda,p}(\omega^{p-k+j/2}G)$ has symbol $\omega^p G\o D^{2k-j}_{\ell-j}$.  It maps $\omega^{p-k+j/2}$ to $\omega^p \o D^{2k-j}_{\ell-j}$.  The conformal symbol $\CS_{\lambda,p}$ extends to its inverse.
\end{prop}

\meno{\bf Definition.} {\em
For $k-p\not\in\bN/2$, let $\pi(\lambda,p)$ be the action $\CS_{\lambda,p}\circ\sigma_{\lambda,p}\circ \CQ_{\lambda,p}$ of $\K$ on $\prod_{j\in\bN} \bF(p-k+j/2)^{(\ell-j)\Pi}$.  
Imitate Section~\ref{CS} to regard it as an infinite matrix with entries
$$ \pi_{rs}^m(\lambda,p): \K \to \Hom\bigl(\bF(s)^{m\Pi}, \bF(r)^{(m+2(r-s))\Pi}\bigr)
   = \bH(s,r-s)^{2(r-s)\Pi}, $$
where $m=\ell-2(s+k-p)\in\bZ_2$.  Varying $k$ and $\ell$ defines $\pi_{rs}^m$ for all $(r,s,m)$ such that $r, s\not\in-\bN/2$ and $r-s\in\bZ/2$.\/}

\begin{thm} \label{psimain}
For $r-s\in 1-\bN/2$, $\pi_{rs}^m$ is $\pi_r^{m\Pi}$ if $r=s$ and~$0$ otherwise.  For $r-s\in 1+\bZ^+/2$, there exist scalars $b^m_{rs}(\lambda,p)$ such that
$$ \pi^m_{rs} = b^m_{rs} \beta_{r-s}(s). $$

The $b^m_{rs}$ are polynomial in $(\lambda,p)$ and rational in $s$.  Their poles are contained in the resonant set where $s\in-\bN/2$.  They are identical to the functions $b^m_{rs}$ defined in Section~\ref{CS}.  In particular, for $3/2\le r-s\le 5/2$ they are given by the formulas of Theorem~\ref{main}.
\end{thm}

At this point let us prove that although $\Psi^k_\ell(\lambda,p)$ is resonant for $k-p\in\bN/2$, its quotient by $\Psi^{p-1/2}_{\ell-1-2(k-p)}$ is essentially non-resonant, in the sense that $\ds$ acts on it semisimply.  Moreover, the formula for the action of $\K$ on it with respect to the conformal symbol is the same as in the usual non-resonant case.  The proof is a simple limiting argument.

\begin{cor} \label{URT}
For $k-p\in \bN/2$, there is a unique even $\ds$-equivalence
$$ \CQ_{\lambda,p}:\bigoplus_{j=0}^{2(k-p)} \bF(p-k+j/2)^{(\ell-j)\Pi} \to 
   \Psi^k_\ell(\lambda,p)/\Psi^{p-1/2}_{\ell-1-2(k-p)} $$
such that $\CQ_{\lambda,p}(\omega^{p-k+j/2}G)$ has symbol $\omega^p G\o D^{2k-j}_{\ell-j}$.  It maps $\omega^{p-k+j/2}$ to the coset of $\omega^p \o D^{2k-j}_{\ell-j}$.  Write $\CS_{\lambda,p}$ for its inverse and $\pi(\lambda,p)$ for $\sigma_{\lambda,p}$ transferred to its domain, as usual.  The matrix entries $\pi_{rs}^m(\lambda,p)$ with $s<r\le 0$ are given by the formulas of Theorems~\ref{main} and~\ref{psimain}.
\end{cor}

\meno{\em Proof.\/}
The first two sentences follow from the same $Q_\ds$-eigenspace argument used for Proposition~\ref{CQCS}.  To prove the last, note that the analogous statement for $\Psi^{k+\ep}_\ell(\lambda,p)/\Psi^{p-1/2+\ep}_{\ell-1-2(k-p)}$ with $2\ep\not\in\bZ$ follows from Theorem~\ref{psimain}.  Furthermore, the polynomials in $Q_\ds$ which project to its eigenspaces have no singularities as $\ep$ goes to zero, and so neither do the formulas for the action of $\K$ on their ranges.  The result follows by continuity.  $\Box$

\subsection{Conjugation} \label{COD}
Define $\C:\Psi(\lambda,p) \to \Psi(1/2-p-\lambda,p)$ by
$$ \C(\omega^p G \o D^z_m) :=
   e^{i\pi(z+m)/2} (-1)^{m|G|} \omega^p \o D^z_m G, $$
where in the complex exponential the element $m$ of $\bZ_2$ is taken to be~$0$ or~$1$.  We will see in Section~\ref{S11} that on \dog s, $\C$ is conjugation.  

Define $\Psi^k_\ell(\lambda,p)$ to be {\em special\/} if $c(\lambda,p)=0$, \ie\ $\lambda=1/4-p/2$.  Note that in the special case, $\C$ preserves $\Psi(\lambda,p)$.  The following result is well-known.  In Section~\ref{S11} we outline a conceptual proof; here we give a direct one.

\begin{prop} \label{conjprops}
\begin{enumerate}

\item[(a)] $\C^2$ acts on $\Psi^k_\ell$ by the scalar $e^{i\pi(2k+\ell)}$.

\smallbreak
\item[(b)] $\C$ is a super anti-involution: $\C(T\circ S)=(-1)^{|T||S|} \C(S) \circ \C(T)$.

\smallbreak
\item[(c)] $\C$ is an even equivalence from $\sigma_{\lambda,p}$ to $\sigma_{1/2-p-\lambda,p}$.

\end{enumerate}
\end{prop}

\meno{\it Proof.\/}  
Part~(a) is easy.  For part~(c), check that $\C(\pi_\lambda(X)) = -\pi_{1/2-\lambda}(X)$ for all $X\in\K$ and apply~(b).  To prove~(b), fix $T\in\Psi(\lambda+p,p')$ and $S\in\Psi(\lambda,p)$ such that neither contains any terms $\o D^z_1$.  First we prove the result for $T\circ S$, which comes down to proving $\C(\px^{z'} G \px^z) = \C(G\px^z)\C(\px^{z'})$.  Applying Leibniz' rule reduces this to the identity
$$ G\px^{z'} = \sum_{n=0}^\infty \Bigl(\sum_{j=0}^n (-1)^j 
   {z'\choose j}{z'-j\choose n-j} \Bigr) G^{(n)} \px^{z'-n}. $$

For the general case we must prove the result for $T\pxi\circ S$, $T\circ S\pxi$, and $T\pxi\circ S\pxi$.  The latter two are straightforward but the first is tricky.  Check each choice of $|T|$ and $|S|$ separately.  When both are odd, $\C(S)\C(T)=0$, so the product rule gives $[\pxi,\C(S)]\C(T) = \C(S)[\pxi,\C(T)]$.  $\Box$

\medbreak
We now prove the most important symmetry of the $b_{rs}^m$, Proposition~\ref{CODsym}.  Propositions~\ref{conjprops} and~\ref{CODsym} immediately yield Corollary~\ref{specsplit}, an explicit description of the splitting of $\sigma_{\lambda,p}$ into the eigenspaces of $\C$ in the special case.

For $k-p\not\in\bN/2$, write $\C^{\CS}$ for the endomorphism $\CS_{1/2-p-\lambda,p}\circ\C\circ\CQ_{\lambda,p}$ of $\prod_{j\in\bN} \bF(p-k+j/2)^{(\ell-j)\Pi}$.  By Proposition~\ref{conjprops}c, $\C^{\CS}$ is an even $\K$-equivalence from $\pi(\lambda,p)$ to $\pi(1/2-p-\lambda,p)$.  Define
$$ t^k_\ell(j) := \exp i\pi\bigl(2k-j+([\ell-j] \mod 2)\bigr)/2. $$

\begin{prop} \label{CODsym}
$\C^{\CS}$ preserves $\bF(p-k+j/2)^{(\ell-j)\Pi}$ and acts on it by the scalar $t^k_\ell(j)$.  As a function of $(c,p)$, $b^m_{rs}$ is of parity $(-1)^{r-s}$ in $c$ for $2(r-s)$ even and parity $(-1)^{m+r-s-1/2}$ for $2(r-s)$ odd.  In other words,
$$ b^m_{rs}(1/2-p-\lambda,p) = (-1)^{2\{r-s\}m + \lf r-s \rf}
   \th b^m_{rs}(\lambda,p). $$
\end{prop}

\meno{\em Proof.\/}
The restrictions of $\pi(\lambda,p)$ and $\pi(1/2-p-\lambda,p)$ to $\ds$ are both $\prod_{j\in\bN} \pi_{p-k+j/2}^{(\ell-j)\Pi}$.  Since $\C^{\CS}$ intertwines them, it must act by some scalar on each factor $\bF(p-k+j/2)^{(\ell-j)\Pi}$.  By Proposition~\ref{psiCQCS2}, $\CQ_{\lambda,p}$ maps $\omega^{p-k+j/2}$ to $\omega^p\o D^{2k-j}_{\ell-j}$, so the scalar is that by which $\C$ acts on $\omega^p\o D^{2k-j}_{\ell-j}$.  This yields the first sentence.

For the second, use the fact that $\C^{\CS}$ is a $\K$-map to prove
$$ b^{\ell-j}_{p-k+i/2,p-k+j/2}(1/2-p-\lambda,p) =
   t^k_\ell(i) t^k_\ell(j)^{-1} b^{\ell-j}_{p-k+i/2,p-k+j/2}(\lambda,p). \quad \Box $$

\begin{cor} \label{specsplit}
In the special case, $\C$ is a $\K$-equivalence from $\Psi^k_\ell$ to itself.  It has two eigenvalues, $\pm t^k_\ell(0)$.  Let $\t\Psi^k_\ell$ be its $t^k_\ell(0)$-eigenspace.  Then its $-t^k_\ell(0)$-eigenspace is $\t\Psi^{k-1+\ell/2}_0$, and we have the $\K$-splitting
$$ \Psi^k_\ell(1/4-p/2,\th p) = \t\Psi^k_\ell(1/4-p/2,\th p) 
   \oplus \t\Psi^{k-1+\ell/2}_0(1/4-p/2,\th p). $$

The image of $\t\Psi^k_0$ under $\CS$ is given by the following equation.  The image of $\t\Psi^{k-1/2}_1$ is the same, except without the $\bF(p-k)$ summand.
$$ \CS\bigl(\t\Psi^k_0(1/4-p/2,\th p)\bigr) =
   \bigoplus_{j\in\bN} \bigl(\bF(p-k+2j) \oplus \bF(p-k+2j+1/2)^\Pi\bigr). $$
\end{cor}

\meno{\bf Remarks.}
In fact, the proof of Proposition~\ref{conjprops}b shows that $\C$ is $\VRoo$-invariant.  The true content of Corollary~\ref{specsplit} is that in the special case, the \dog\ modules split under $\VRoo$ along even and odd degree: the finer filtration is irrelevant here.  Thus $\C$ is more properly associated to the full $\VRoo$ action (for which it would, of course, be interesting to carry out the entire \cite{CMZ97} program).

The special case is a useful source of examples, as the modules $\t\Psi^k_\ell$ are very ``tightly knit'': of all the $\K$-modules of infinite length composed of \tdm s, they seem to be the closest to uniserial.

\subsection{The supercircle} \label{S11}
In order to use certain dualities arising from the Berezinian, we now extend our definitions to the supercircle $\Soo$.  The set $\PSoo$ of polynomial functions on $\Soo$ has basis $\{\xi^\ell x^n: \ell=0, 1; n\in\bZ\}$.  Thus $\VRoo$, $\K$, $\bF(\lambda)$, $\bD(\lambda,p)$, and $\Psi^k_\ell(\lambda,p)$ may all be defined over $\Soo$ simply by admitting negative powers of $x$.  We will denote the resulting extensions with a hat: $\h\K$, $\h\bF(\lambda)$, etc.  

The bracket on $\K$ and the formula for $\pi_\lambda$ involve only derivatives, and so depend polynomially on the exponents of $x$ in their arguments.  Therefore the same is true of $\sigma_{\lambda,p}$, $\CS_{\lambda,p}$, $\CQ_{\lambda,p}$, $\pi^\ell_{rs}(\lambda,p)$, and $\beta_p(\lambda)$.  For example, to prove this for $\CS_{\lambda,p}$ and $\CQ_{\lambda,p}$, use the fact that projection to an eigenspace of $Q_\ds$ is given by the action of a polynomial in $Q_\ds$.

It follows that the formulas for all of these functions over $\Roo$ are valid over $\Soo$.  (We do not work over $\Soo$ from the start because many of the lowest weight arguments become muddier but no significant results change.)

We now define the Berezinian and the associated pairing.  Proposition~\ref{Ber} gives its properties; it is well-known and we will not give the proof.  It yields a conceptual explanation of Proposition~\ref{conjprops}.

\meno{\bf Definition.}{\em
The\/ {\em Berezinian\/} $\Ber:\h\bF(1/2)\to\bC^{1|0}$ and the pairing $B:\h\bF(1/2-\lambda) \ot \h\bF(\lambda) \to \bC^{1|0}$ are defined by\/}
$$ \Ber(\omega^{1/2}F) := {\ts \frac{1}{2\pi i} \oint_{S^1} \pxi F dx}, \quad
   B(\omega^{1/2-\lambda}F \ot \omega^\lambda G) := \Ber(\omega^{1/2} FG). $$

\begin{prop} \label{Ber}
The Berezinian is a non-zero odd $\h\K$-map.  Up to a scalar, it is the only non-zero $\h\K$-map from $\h\bF(\lambda)$ to $\bC^{1|0}$ for any $\lambda$.

The pairing $B$ is odd, non-degenerate, and $\h\K$-invariant.  Up to a scalar, it is the only $\h\K$-invariant pairing of $\h\bF(1/2-\lambda)$ with $\h\bF(\lambda)$.  Regarded as a form on $\bigoplus_\lambda \h\bF(\lambda)$, it is symmetric.  On \dog s, the map $\C$ of Proposition~\ref{conjprops} is equal to the map $T\mapsto T^B$ of Lemma~\ref{tauB}.
\end{prop}

Next we define the $\bZ_2$-graded versions of the noncommutative residue and the Adler trace.  Proposition~\ref{Adler} gives their properties and Proposition~\ref{SNCRsym} is the resulting symmetry of the $b^m_{rs}$.  Proposition~\ref{Adler} follows from Proposition~\ref{Ber} essentially as in the non-super case.  We leave the details to the reader (see \cite{FR91} for the case $\lambda=p=0$), except to remark that the $\h\K$-invariance of the super noncommutative residue may be explained as follows: $\SNCR$ is zero on $\h\Psi^k_\ell$ unless $k+1/2$ is in $\bN/2$ and $2k+\ell$ is even, when $\h\Psi^k_\ell=\h\bD^k\oplus\h\Psi^{-1/2}_1$.  In this case it is zero on $\h\bD$ and is the Berezinian of the symbol on $\h\Psi^{-1/2}_1$.

Note that any \psidog\ may be written as $\omega^p \sum_{z\in\bC, m\in\bZ_2} G_{z,m} \o D^z_m$, where the $G_{z,m}$ are zero for $z$ off some set $k-\bN$.

\meno{\bf Definition.} {\em
The\/ {\em super noncommutative residue\/} is the map
$$ \SNCR:\h\Psi(\lambda,0)\to\bC^{1|0}, \quad
   \SNCR\Bigl(\sum_{z, m} G_{z,m} \o D^z_m\Bigr) 
   := \Ber\bigl(\omega^{1/2} G_{-1,1}\bigr). $$

The\/ {\em super Adler trace\/} is the pairing\/}
$$ A:\h\Psi(\lambda+p,-p)\ot\h\Psi(\lambda,p) \to \bC^{1|0}, \quad
   A(T\ot S) := \SNCR(T\circ S). $$

\begin{prop} \label{Adler}
$\SNCR$ is a non-zero even $\h\K$-map.  The super Adler trace is even, non-degenerate, and $\h\K$-invariant.  Regarded as a form on $\bigoplus_{\lambda,p} \h\Psi(\lambda,p)$, it is supersymmetric.
\end{prop}

\begin{prop} \label{SNCRsym}
There is a symmetry across the antidiagonal $r+s=1/2$ relating the entries of the matrix $\pi$ at a given value of $(c,p)$ to those with the same value of~$c$ but opposite~$p$.  More precisely,
$$ b^m_{rs}(\lambda+p,-p) = (-1)^{2\{r-s\}m + \lf r-s \rf}
   \th b^{m+1+2(r-s)}_{1/2-s,1/2-r}(\lambda,p). $$
\end{prop}

\meno{\it Proof.\/}
As stated in Proposition~\ref{Ber}, the odd $\h\K$-invariant pairing $B$ of $\h\bF(1/2-\lambda)$ with $\h\bF(\lambda)$ may be regarded as a supersymmetric form on $\bigoplus_{\lambda\in\bC}\h\bF(\lambda)$.  Taking this point of view, one checks that
\begin{equation} \label{Bformulas}
   (\omega^p)^B = \omega^p, \quad G^B = G, \quad \px^B = -\px, \quad
   \pxi^B = -\pxi, \quad \ep^B = -\ep.
\end{equation}

Throughout this proof, when we write $\ep$ for an endomorphism of either $\h\bF(\lambda)$ or $\h\bF(\lambda)^\Pi$ we shall mean $\ep_{\h\bF(\lambda)}$.  Thus $\ep_{\h\bF(\lambda)^\Pi}=-\ep$.  We will also abbreviate $\h\bF(\lambda)^{m\Pi}$ to $\h\bF(\lambda)^m$.

Using Lemma~\ref{tensorhom} and Section~\ref{Parity}, one finds that for $m\in\bZ_2$, $\ep\ot\ep^m$ (not $\ep\ot_s\ep^m$) is an odd $\h\K$-equivalence from $\h\bF(1/2-\lambda)^{m+1}\ot\h\bF(\lambda)^m$ to $\h\bF(1/2-\lambda)\ot\h\bF(\lambda)$.  Therefore $B\circ(\ep\ot\ep^m)$ is an even $\h\K$-invariant pairing of $\h\bF(1/2-\lambda)^{m+1}$ with $\h\bF(\lambda)^m$.  We will write $\o B$ for it.

Let $\h\bF$ be $\bigoplus_{\lambda\in\bC} \bigl(\h\bF(\lambda) \oplus \h\bF(\lambda)^\Pi\bigr)$.  Regard $\o B$ as an even supersymmetric form on $\h\bF$.  The parity change functor $\Pi$ is an endomorphism of $\h\bF$ in an obvious way: it switches the values of $m$.  Let $M$ be the endomorphism of $\h\bF$ taking the value $m$ on $\bF(\lambda)^m$.  View all endomorphisms of $\bigoplus_\lambda \h\bF(\lambda)$ as endomorphisms of $\h\bF$ commuting with $\Pi$.  (For example, $\ep = (-1)^M \ep_{\h\bF}$.)  We leave the reader to check that
\begin{equation} \label{oBformulas}
   \Pi^{\o B} = (-1)^M\Pi, \quad M^{\o B} = 1-M, \quad
   \ep^{\o B} = -\ep, \and\ T^{\o B} = T^B,
\end{equation}
where $T$ is any endomorphism of $\bigoplus_\lambda \h\bF(\lambda)$.

The strategy is to transfer the Adler trace from \psidog s to $\h\bF$ via $\CQ$.  However, $\CQ$ is not defined in the resonant case, so we must delete from $\h\bF$ all $\h\bF(\lambda)^m$ with $\lambda\in-\bN/2$.  This leaves $\o B$ degenerate on those $\h\bF(\lambda)^m$ with $\lambda\in\bZ^+/2$, so we will delete them also.  Finally, we must allow formal series appropriate to \psidog s.  Thus we arrive at the following definition.

Let $\h\bF_{\reg}$ be the space of formal series of elements in those $\h\bF(\lambda)^m$ with $\lambda\not\in\bZ/2$ ranging over finitely many cosets of $\bN$.  Extend $\o B$ to $\h\bF_{\reg}$ in the natural way.  As usual, write $\pi(\lambda,p)$ for the \r\ of $\h\K$ on $\h\bF_{\reg}$ obtained by transferring the \r\ $\sigma_{\lambda,p}$ on $\h\Psi(\lambda,p)$ via $\CQ_{\lambda,p}$, deleting the terms giving those $\bF(r)^m$ with $r\in\bZ/2$.

We claim that $\o B$ is $-A\circ\CQ_{\lambda,p}$ for all $(\lambda,p)$, and is therefore an invariant pairing of $\pi(\lambda+p,-p)$ with $\pi(\lambda,p)$.  To prove this, note that for any $r\not\in\bZ/2$ the restriction of $A\circ\CQ_{\lambda,p}$ to $\h\bF(1/2-r)^{m+1}\ot \h\bF(r)^m$ is $\ds$-invariant.  Check that up to a scalar there is only one such pairing, so it must be a multiple of $\o B$.  Evaluating at an appropriate spot shows that the multiple is $-1$, independent of $r$.

Thus Proposition~\ref{tauB} gives $\pi(\lambda,p)(X)^{\o B} = -\pi(\lambda+p,-p)(X)$ for all $X$ in $\h\K$.  Therefore
$$ \pi^m_{rs}(\lambda,p)(X)^{\o B} = -\pi^{m+1+2(r-s)}_{1/2-s,1/2-r}(\lambda+p,-p)(X) $$
for $r, s \not\in\bZ/2$.  In light of Theorem~\ref{psimain}, the proof is finished if we prove that 
$$ \beta_{r-s}(s)^{\o B}(X) = 
   -(-1)^{2\{r-s\}M + \lf r-s \rf} \th \beta_{r-s}(1/2-r)(X) $$
(extend to the non-resonant cases with $r,s\in\bZ/2$ by continuity).  In fact, by Lemma~\ref{beta} it is enough to prove this at $X=e_{3/2}$, as by Proposition~\ref{tauB} both sides are $\ds$-relative.  Here it follows from~(\ref{Bformulas}), (\ref{oBformulas}), and $\beta_p(\lambda)(e_{3/2}) = 4\omega^p\o D^{2p-3} \ep^{2p} \Pi^{2p}$.  (The $\Pi^{2p}$ does not appear in Lemma~\ref{beta} because we usually write $1$ for the identity map from $V$ to $V^\Pi$.)  $\Box$

\section{Applications}  \label{Apps}
\setcounter{lemma}{0}

The simplest application of Theorem~\ref{psimain} is the construction of $\Symb_{3/2}$, the unique $\K$-invariant sesquisymbol.  (For comparison, we note that in general there is a unique $\VR$-invariant double symbol over the line, but only an ordinary $\VRM$-invariant symbol over $\bR^m$ for $m>1$; see \cite{LO99}.)  For $k-p\not=-1/2$ or~$0$, define
\begin{equation} \label{sesquidfn}
   \Symb_{3/2}: \Psi^k_\ell(\lambda,p) \to \bF(p-k)^{\ell\Pi} \oplus
   \bF(p-k+1/2)^{(\ell+1)\Pi} \oplus \bF(p-k+1)^{\ell\Pi}
\end{equation}
to be $\CS_{\lambda,p}$ followed by projection along $\prod_{j\in 3+\bN} \bF(p-k+j/2)^{(\ell-j)\Pi}$ (see Proposition~\ref{psiCQCS2}).  For \dog s, the following proposition is essentially a corollary of the uniqueness statement in Proposition~\ref{CQCS} and the first statement of Corollary~\ref{uptoscalar}.  See Proposition~\ref{psiCQCS2}, Theorem~\ref{psimain}, and Corollary~\ref{URT} for the necessary extensions to \psidog s.

\begin{prop} \label{sesqui}
$\Symb_{3/2}$ is the unique $\K$-invariant map as in~(\ref{sesquidfn}) carrying $\omega^p G \o D^{k-j}_{\ell-j}$ to $\omega^{p-k+j/2} G$ modulo higher \tdm s.
\end{prop}

\subsection{Invariant transvectants} \label{ITs}
The aim of this section is to show how our results can be used to reproduce the classification of $\K$-invariant supertransvectants obtained by Leites, Kochetkov, and Weintrob \cite{LKW91}.  (See \cite{GO07} for applications of the $\K$-invariance of $\bJ_{1/2}^{\mu,\nu}$ and $\bJ_1^{\mu,\nu}$.)

To recall the situation for $\bR$, let $F(\lambda)$ be the $\VR$-module of tensor densities of degree~$\lambda$ (as in Section~\ref{FCoho}).  For $\mu+\nu\not\in-\bN$ and $k\in\bN$, there is a projectively invariant \dog\ 
$$ J^{\mu,\nu}_k:F(\mu)\ot F(\nu)\to F(\mu+\nu+k), $$
unique up to a scalar.  These are the {\em transvectants,\/} discovered ``in antiquity'' by P.~Gordan.  Their super analogs are discussed in Section~\ref{TPs}.

The cases in which $J^{\mu,\nu}_k$ is $\VR$-invariant are known.  Still assuming that $\mu+\nu \not\in\bN$, it is always so for $k=0$ or~$1$, and never so for $k\ge 4$.  For $k=2$ it is so \iff\ $\mu\nu=0$, and for $k=3$ it is so \iff\ $\mu=\nu=0$ or~$-2/3$.  All of the invariant cases with $k=2$ or~$3$ are compositions of a $J_1$ with an operator built from the divergence $\nabla:F(0)\to F(1)$ with the exception of $J^{-2/3,-2/3}_3$, the {\em Grozman operator\/} \cite{Gr80}.  This operator has no super analog.

The classification of the $\VR$-invariant $J^{\mu,\nu}_k$ may be deduced from the $\bR$-analog of our main theorem (an expository account of which may be found in \cite{Co05}).  The idea is to construct a $\VR$-injection from $F(\mu)\ot F(\nu)$ to the space of \psidog s of order~$\le -1$ from $F(1-\mu)$ to $F(\nu)$.  

We now carry this out for $\K$.  Recall from Lemma~\ref{Bol}d that $\SBol_{1/2}(0)$ is $\K$-invariant.

\begin{prop} \label{psimoh}
For $\mu+\nu\not\in-\bN/2$, there is an odd $\K$-injection $\bJ$ from $\bF(\mu)\ot\bF(\nu)$ to $\Psi^{-1/2}_1(1/2-\nu,\mu+\nu-1/2)$, unique up to a scalar.  There are non-zero scalars $c_j$ such that $\bJ$ has the form
\begin{equation} \label{J}
   \bJ := \CQ_{1/2-\nu, \mu+\nu-1/2}\circ\bigl(\ts\bigoplus_{j=0}^\infty 
   c_j \ep^{j+1}\circ \bJ^{\mu,\nu}_{j/2} \bigr).
\end{equation}
\end{prop}

\meno{\em Proof.\/}
Note that by Propositions~\ref{Fdecomp} and~\ref{ST}, the map $(\bigoplus_0^\infty \ep^j\circ \bJ^{\mu,\nu}_{j/2})$ from $\bF(\mu)\ot\bF(\nu)$ to $\bigoplus_0^\infty \bF(\mu+\nu+j/2)$ is an even $\ds$-equivalence.  Now work over the supercircle as in Section~\ref{S11}, and recall the Berezinian pairing $B$ from Proposition~\ref{Ber}.  Check that there is an odd non-degenerate $\K$-invariant pairing between $\h\bD(\mu,1/2-\mu-\nu)$ and $\h\bF(\mu)\ot\h\bF(\nu)$, defined by
$$ \la T, \omega^\mu F \ot \omega^\nu G \ra := B\bigl(T(\omega^\mu F) \ot \omega^\nu G\bigr). $$
Since $\h\bD^k$ is $\ds$-equivalent to $\bigoplus_{j=0}^{2k}\h\bF(1/2-\mu-\nu-j/2)$, the annihilator $(\h\bD^k)^\perp$ must be the kernel of $\bigoplus_{j=0}^{2k} \bJ^{\mu,\nu}_{j/2}$.

In Proposition~\ref{Adler} we saw that there is an even non-degenerate $\K$-invariant pairing between $\h\bD(\mu, 1/2-\mu-\nu)$ and $\h\Psi^{-1/2}_1(1/2-\nu, \mu+\nu-1/2)$.  The annihilator of $\h\bD^k$ with respect to this form is $\h\Psi^{-1-k}_{2k}$.  Therefore for each $k$ there is an odd $\K$-equivalence from $\h\bF(\mu)\ot\h\bF(\nu)$ modulo the kernel of $\bigoplus_{j=0}^{2k} \bJ^{\mu,\nu}_{j/2}$ to $\h\Psi^{-1/2}_1/\h\Psi^{-1-k}_{2k}$.  The restriction to the superline of the projective limit of these equivalences is the injection $\bJ$ of the proposition.  It follows from Theorem~\ref{psimain} and Proposition~\ref{ST} that all $\ds$-maps between the two spaces have the form~(\ref{J}).  $\Box$

\begin{cor} \label{NoGrozman} \cite{LKW91}
Assume $\mu+\nu\not\in-\bN/2$.  Then $\bJ^{\mu,\nu}_k$ is $\K$-invariant \iff\ either $k=0$, $1/2$, or $1$; $k=3/2$ and $\mu\nu=0$; or $k=2$ and $\mu=\nu=0$.  

We have the following proportionalities:
$$ \bJ^{0,\nu}_{3/2} \propto \bJ^{1/2,\nu}_1 \circ \bigl(\SBol_{1/2}(0)\ot 1\bigr), \quad
   \bJ^{\mu,0}_{3/2} \propto \bJ^{\mu,1/2}_1 \circ \bigl(1\ot\SBol_{1/2}(0)\bigr), $$
and $\bJ^{0,0}_2 \propto \bJ^{1/2,1/2}_1\circ\bigl(\SBol_{1/2}(0)\ot\SBol_{1/2}(0)\bigr)$.
\end{cor}

\meno{\em Proof.\/}
Applying Proposition~\ref{psimoh}, we find that $\bJ^{\mu,\nu}_k$ is $\K$-invariant \iff\ $\bF(\mu+\nu+k)$ is a quotient of $\Psi^{-1/2}_1(1/2-\nu, \mu+\nu-1/2)$.  By Theorem~\ref{psimain}, this occurs \iff\ $b^{2l+1}_{\mu+\nu+k,\mu+\nu+l}(1/2-\nu,\mu+\nu-1/2)$ is~$0$ for $0\le l\le k-3/2$.  Use the formulas of Theorem~\ref{main} to complete the proof.  $\Box$

\subsection{$\bH(\lambda,p)$-valued cohomology} \label{KCoho}
The computation of the Ext groups between \tdm s and the indecomposable \r s composed of them goes back to Feigin and Fuchs \cite{FF80}.  In a major achievement, they computed $\Ext^n_{\VR}\bigl(F(\lambda),F(\lambda+p)\bigr)$ for all~$n$.  

Let us recall the $n=1$ result in the non-resonant case $p\not=0$ or $1-2\lambda$.  Here all cohomology classes are $\dsl_2$-relative by the $\VR$-analog of Lemma~\ref{IC}, $\dim \Ext^1$ is always either~$0$ or~$1$, and it is~$1$ precisely when $p=2$, $3$, or~$4$; $p=5$ and $\lambda=-4$ or~$0$; or $p=6$ and $2\lambda = -5\pm\sqrt{19}$.

In this section we compute the $\ds$-relative $\Ext^1$ groups and some of the $\Ext^2$ groups over $\K$.  Note in particular that Theorem~\ref{Ext1}f gives the $\K$-analogs of the $p=6$ cocycles of $\VR$.  We begin by proving that the $\dt$- and $\ds$-relative $\Ext$ groups of $\K$ are \dog-valued.  In fact, the proof applies equally well to $\VR$, and we expect also to $\VRM$.

\meno {\bf Notation.}
Recall that for any subalgebra $\dh$ of $\K$, $H^n\bigl(\K,\dh,\bH(\lambda,p)\bigr)$ and $\Ext^n_{\K,\dh}\bigl(\bF(\lambda),\bF(\lambda+p)\bigr)$ are both standard notation for the space of $\dh$-relative $\Hom\bigl(\bF(\lambda),\bF(\lambda+p)\bigr)$-valued $n$-cohomology of $\K$.  For brevity, let us denote this space and its dimension as follows:
$$ \Ext^n_{\K,\dh}(\lambda,p) := \Ext^n_{\K,\dh}\bigl(\bF(\lambda),\bF(\lambda+p)\bigr),
   \ \ \dim^n_{\K,\dh}(\lambda,p) := \dim \Ext^n_{\K,\dh}(\lambda,p). $$

We will focus on $\ds$-relative cohomology, and it will be convenient to use the abbreviations $C^n_{\lambda,p}$, $Z^n_{\lambda,p}$, and $B^n_{\lambda,p}$ for the $\ds$-relative $n$-cochains, cocycles, and coboundaries of $\bH(\lambda,p)$.  However, we will do so only in the proofs.

\begin{lemma} \label{ExtN}
For all $\lambda$, $p$, and $n$, $C^n\bigl(\K,\dt,\bH(\lambda,p)\bigr)$ is $\bD^{p-n/2}(\lambda,p)$-valued for $2p\in n+\bN$, and~$0$ otherwise.  Similarly, $C^n\bigl(\K,\ds,\bH(\lambda,p)\bigr)$ is $\bD^{p-3n/2}(\lambda,p)$-valued for $2p\in 3n+\bN$, and~$0$ otherwise.

In particular, all $\dt$- and $\ds$-relative cohomology of $\bH(\lambda,p)$ is $\bD(\lambda,p)$-valued, $\Ext^n_{\K,\dt}(\lambda,p) = 0$ for $2p\not\in n+\bN$, and $\Ext^n_{\K,\ds}(\lambda,p) = 0$ for $2p\not\in 3n+\bN$.
\end{lemma}

\meno{\it Proof.\/}
By~(\ref{Kmodsubalg}) and Lemma~\ref{SymF}, $\Lambda^n_s(\K/\dt)$ is $\dt$-equivalent to a sum of $\bF(\mu)^{2\mu\Pi}$ with $2\mu\in n+\bN$, and $\Lambda^n_s(\K/\ds)$ is $\ds$-equivalent to a sum of $\bF(\mu)^{2\mu\Pi}$ with $2\mu\in 3n+\bN$.  Apply Lemma~\ref{good}e.  $\Box$

\begin{thm} \label{Ext1}
Excepting the case $p=0$ and the special case $\lambda = \frac{1}{4} - \frac{p}{2}$, the inclusion map from $\Ext^1_{\K,\ds}(\lambda,p)$ to $\Ext^1_\K(\lambda,p)$ is an isomorphism.  The space $\Ext^1_{\K,\ds}(\lambda,p)$ is of parity $(-1)^{2p}$ for all $(\lambda,p)$.  When $p=2$, $5/2$, or~$3$ and in addition $\lambda=0$ or~$-p+1/2$, it consists of non-trivial cup products.

\begin{enumerate}

\item[(a)]  For $p=0$ or $(\lambda,p) = (0,\, \frac{1}{2})$, $\dim^0_\K(\lambda,p) = 1$.  Otherwise it is $0$.

\smallbreak
\item[(b)]  For $p\not\in\{\frac{3}{2},\, 2,\, \frac{5}{2},\, 3,\, 4\}$, $\dim^1_{\K,\ds}(\lambda,p) = 0$ for all $\lambda$.

\smallbreak
\item[(c)]  For $p=\frac{3}{2}$ or $\frac{5}{2}$, $\dim^1_{\K,\ds}(\lambda,p)=1$ unless $\lambda = \frac{1}{4} - \frac{p}{2}$, when it is~$0$.

\smallbreak
\item[(d)]  For $p=2$, $\dim^1_{\K,\ds}(\lambda,p) = 1$ for all $\lambda$.

\smallbreak
\item[(e)]  For $p=3$, $\dim^1_{\K,\ds}(\lambda,p) = 1$ for $\lambda = -\frac{5}{2}$ or~$0$, and~$0$ otherwise.

\smallbreak
\item[(f)]  For $p=4$, $\dim^1_{\K,\ds}(\lambda,p) = 1$ for $\lambda= \frac{-7\pm\sqrt{33}}{4}$, and~$0$ otherwise.

\end{enumerate}
\end{thm}

\meno{\em Proof.\/}
The first and second sentences follow from Proposition~\ref{ghost} and Lemmas~\ref{IC} and~\ref{beta}.  Next, we claim that for $p\in 1/2+\bN$ and $2q\in 3+\bN$,
\begin{equation} \label{nontrivcups} \begin{array}{rcl}
   \SBol_p\bigl(\frac{1}{4} - \frac{p}{2}\bigr) \cup 
   \o\beta_q\bigl(\frac{1}{4} - \frac{p}{2} - q\bigr) &=&
   \o\beta_{p+q}\bigl(\frac{1}{4} - \frac{p}{2} - q\bigr), \\[4pt]
   \o\beta_q\bigl(\frac{1}{4} + \frac{p}{2}\bigr) \cup
   \SBol_p\bigl(\frac{1}{4} - \frac{p}{2}\bigr) &=&
   - \o\beta_{p+q}\bigl(\frac{1}{4} - \frac{p}{2}\bigr).
\end{array} \end{equation}
To verify this, note that in each equation both sides are $\ds$-invariant, and apply Lemmas~\ref{Bol} and~\ref{srelC1}.  To prove the third sentence, use~(\ref{nontrivcups}) and parts~(a), (c), and~(d) of the theorem to see that the 1-classes in question are cup products of the 0-class $[\SBol_{1/2}(0)]$ with the appropriate 1-class (see below).

Recall $c:=\lambda+p/2-1/4$ from Section~\ref{CS}.  By Lemma~\ref{Bol}c, $C^0_{\lambda,p}$ is $\bC 1$ if $p=0$, $\bC\SBol_p(\lambda)$ if $p\in 1/2+\bN$ and $c=0$, and~$0$ otherwise.  By Lemma~\ref{Bol}d, $Z^0_{\lambda,p}$ is $\bC 1$ if $p=0$, $\bC\SBol_{1/2}(0)$ if $(\lambda,p)=(0,1/2)$, and~$0$ otherwise.  This gives part~(a), and we also see that $B^1_{\lambda,p}$ is $\bC\partial\SBol_p(\lambda)$ if $p\in 3/2+\bN$ and $c=0$, and~$0$ otherwise.

By Lemma~\ref{beta}, $C^1_{\lambda,p} = \bC\o\beta_p(\lambda)$ for $2p\in 3+\bN$, and~$0$ otherwise.  By Lemma~\ref{ExtN}, $C^2_{\lambda,p}=0$ for $p\not\in 3+\bN/2$.  This proves~(c) and~(d).  To proceed further, we must decide when $\partial\o\beta_p$ is zero.

Recall~(\ref{wedge2s}): $\Lambda^2_s(\K/\ds) \scong \bF(3) \oplus \bF(9/2)^\Pi \oplus \bF(5) \oplus\cdots$.  Since $\partial\o\beta_p$ is $\ds$-relative, it annihilates some $\ds$-copy of $\bF(\mu)^{2\mu\Pi}$ \iff\ it annihilates the \lwv\ of weight $\mu$.  Moreover, it must map this \lwv\ to a \lwv\ of weight $\mu$, \ie\ a multiple of $\omega^p \o D^{2(p-\mu)}$.  In particular, for $p=3$ or~$4$, $\partial\o\beta_p$ is~$0$ \iff\ it annihilates the $\bF(3)$.

The \lwv\ of the $\bF(3)$ is $\bigwedge^2 e_{3/2}$.  Check that if $V$ is any $\K$-module and $\phi\in C^1(\K,\ds, V)$ has parity $|\phi|$, then 
$$ \partial\phi({\ts\bigwedge}^2 e_{3/2}) = (-1)^{|\phi|} (2e_{3/2}-e_{1/2}^3) \phi(e_{3/2})
   = 2(-1)^{|\phi|}\t s_{3/2} \phi(e_{3/2}), $$
where $\t s_{3/2}$ is the step algebra element given in Lemma~\ref{345}.  Hence by~(\ref{obeta}), $\partial \o\beta_p(\lambda)(\bigwedge^2 e_{3/2})=0$ \iff\ $\sigma_{\lambda,p}(\t s_{3/2}) (\omega^p \o D^{2p-3})=0$.

The following formulas can be verified directly:
\begin{equation} \begin{array}{rcl} \label{coboundary}
   8\sigma_{\lambda,p}(\t s_{3/2}) \bigl(\omega^p \o D^{2p-3}\bigr)
   &=& -(p-2)(16c^2 - 8p -1)\omega^p \o D^{2p-6}, \\[4pt]
   4\sigma_{\lambda,p}(\t s_{3/2}) \bigl(\omega^p \o D^{2p-3}\bigr)
   &=& -(2p-3)(2p-5) c \omega^p \o D^{2p-6},
\end{array} \end{equation}
for $p$ integral and half-integral, respectively.  However, it is more efficient to apply Section~\ref{Proofs}: by~(\ref{sigmats}), $\sigma_{\lambda,p}(\t s_{3/2})(\omega^p\o D^{2p-3}_m)=0$ \iff\ $P^m_{3,3/2}(\lambda,p)=0$, which occurs \iff\ $b^m_{3,3/2}(\lambda,p)=0$.  Now use Theorem~\ref{main}.  

This proves~(e) and~(f), as well as $\partial\o\beta_p\not=0$ for $p\in 5+\bN$ and $16c^2 \not= 8p+1$, and for $p\in 7/2 + \bN$ and $c\not=0$.  Since $\o\beta_p\propto\partial\SBol_p$ for $p\in 3/2+\bN$ and $c=0$, to complete the proof of~(b) it only remains to prove that $\partial\o\beta_p\not=0$ for $p\in 5+\bN$ and $16c^2 = 8p+1$, the cases in which $\partial\o\beta_p$ annihilates the $\bF(3)$.  

We claim that in these cases $\partial\o\beta_p$ also annihilates the $\bF(9/2)^\Pi$.  For if not, its symbol would be a non-zero element of $Z^2\bigl(\K,\ds,\bF(9/2)^\Pi\bigr)$, contradicting Proposition~\ref{F012coho}f.  Therefore it has image in $\bD^{2p-10}$.  Since we have not simplified the formula for $P^1_{5,3/2}$, we will verify directly that here $\partial\o\beta_p(\bigwedge^2 e_{5/2}) \not= 0$.  Check that it is $2\bigl(\sigma_{\lambda,p}(e_{5/2})\o\beta_p(e_{5/2}) - \o\beta_p(e_5)\bigr)$, and then use~(\ref{Brsfact1}) and~(\ref{obeta}) to prove that this is a non-zero multiple of $\sigma_{\lambda,p}(e_{5/2} e_{1/2}^2 - e_{1/2}^7)(\omega \o D^{2p-3})$.

Since the latter expression lies in $\bD^{2p-10}_5$, it is a multiple of $\omega^p \o D^{2p-10}$.  Looking ahead to Corollary~\ref{Zen} and applying $Z$, we find that the multiple is
$$ \ts \Bigl(\frac{c}{3} - \frac{1}{4}\Bigr) \Bigl[{p-2 \atop 3}\Bigr] \Bigl(36c-12p-15 - 
   \frac{1}{2} \Bigl[{2c-5/2 \atop 3}\Bigr]\Bigr). $$
If $8p=16c^2-1$, this expression becomes a multiple of $x^3+24x^2+188x+384$, where $x=4c-7$.  This cubic is irreducible over $\bQ$, so it is never zero for $p\in 5+\bN$, where $c$ is quadratic over $\bQ$.  $\Box$

\begin{prop} \label{Ext2}
Consider $\dim^2_{\K,\ds}(\lambda,p)$.  For $p\not\in 3+\bN/2$, it is~$0$.  For $p=7/2$ and~$9/2$, it is~$1$ in the special case $c=0$, and~$0$ otherwise.  For $p=3$ and~$4$, it is~$1$ for $16c^2=8p+1$, and~$0$ otherwise.  For $p=5$, it is~$1$ for all~$\lambda$.  For $p=11/2$, it is~$2$ for $c=0$, and~$1$ otherwise.
\end{prop}

\meno{\it Proof.\/}
The second sentence is from Lemma~\ref{ExtN}.  By~(\ref{wedge2s}), $\dim C^2_{\lambda,p}$ is $1$ for $p=3$, $7/2$, and~$4$; $2$~for $p=9/2$; and~$3$ for $p=5$ and~$11/2$.  By~(\ref{wedge3s}), $\dim C^3_{\lambda,p}$ is~$0$ for $p<9/2$ and~$1$ for $p=9/2$, $5$, and~$11/2$.  Moreover, $C^3_{\lambda,p} = B^3_{\lambda,p}$ for these $p$-values: by Proposition~\ref{F012coho}f, the 2-cochain with values in the copy of $\bF(9/2)^\Pi$ is not a cocycle even at the symbol level.  Since $C^1_{\lambda,p} = \bC\o\beta_p$, Theorem~\ref{Ext1} gives the necessary $B^2_{\lambda,p}$.  $\Box$

\meno {\bf Remarks.}
We will see in Section~\ref{Resonant} that when $c=0$, the non-trivial $\ds$-relative 2-classes at $p=7/2$ and~$9/2$, as well as one of the two 2-classes at $p=11/2$, are in fact the coboundaries of the $\dt$-relative 1-cochains $\o\alpha_p$ discussed there.  It also seems likely that $\dim^2_{\K,\ds}(\lambda,6)=1$ for all $\lambda$.  The situation appears to be parallel to that for $\VR$ \cite{FF80}, leading us to the following provisional extension of Conjecture~1:

\meno {\bf Conjecture~2.}  Assume $p\in\bN/2$, and define $G^\pm(n)=n^2\pm n/2$.  For $p\ge G^+(n+1)$ or $p < G^-(n)$, $\dim^n_{\K,\dt}(\lambda,p) = 0$.  For $G^-(n+1) \le p < G^+(n+1)$, $\dim^n_{\K,\dt}(\lambda,p) = 0$ unless $p - G^-(n+1)$ is integral and $\lambda$ is one of two particular values $\lambda^\pm_n(p)$ lying in some quadratic extension of $\bQ$, when it is~$1$.  For $G^+(n) \le p < G^-(n+1)$, $\dim^n_{\K,\dt}(\lambda,p) = 1$.  For $G^-(n) \le p < G^+(n)$, $\dim^n_{\K,\dt}(\lambda,p) = 0$ unless $p - G^-(n)$ is integral and $\lambda = \lambda^\pm_{n-1}(p)$, when it is~$1$.  Finally, for $n\ge 2$ all $\dt$-relative cohomology classes are $\ds$-relative.

\medbreak
Let us mention further that the modules $\bD(\lambda,p)$ exhibiting cohomology have only one continuous parameter (namely, $\lambda$), while those studied by Feigin and Fuchs have two.  The attempt to obtain a second parameter by studying the modules $\Psi(\lambda,p)$ does not immediately succeed: in general these modules have no cohomology because their composition series extends infinitely far ``to the right''.  However, this can easily be remedied by passing to quotients, suggesting the following:

\meno{\bf Question.}  Fix $\ell\in\bZ_2$ and $\lambda$, $p$, $q\in \bC$ such that $p-q \in \bZ/2$.  What is
$$ H^\bullet\Bigl(\K,\bigl\{\cup_{k\in\bN} \Psi^{q+k}_\ell(\lambda,p)\bigr\}/ 
   \Psi^{q-1/2}_{\ell+1}(\lambda,p)\Bigr)? $$

\medbreak
We hope that this question will lead to a deeper analogy with the ``parabola picture'' of Feigin and Fuchs.  For example, here~(\ref{coboundary}) holds for $2(p-q) + \ell$ even and odd, respectively.  In particular, for $2(p-q) + \ell$ even and $p-q = 3$, $7/2$, $4$, or $9/2$, the above $H^1$ is $\bC$ for $16c^2 = 8p+1$ and $0$ otherwise.

\subsection{Uniserial modules} \label{Uniserial}
The $\ds$-split extensions of \tdm s are associated to $\Ext^1_{\K,\ds}(\lambda,p)$ and $\Ext^2_{\K,\ds}(\lambda,p)$ (see Section~\ref{LengthN}).  Two questions arise.  First, what are the \ind\ extensions?  And second, which of them can be realized in some natural setting, for example as \sq s of a differential operator or tensor product module?  In the case of $\VR$, to our knowledge O.~Mathieu was the first to ask the latter question for \dog\ modules.  For tensor products the idea is in \cite{FF80}.  

The length~2 extensions of \tdm s of $\VR$ are given in \cite{FF80}, and those of length~3 having distinct Casimir eigenvalues were computed in \cite{Co01}.  It is possible to realize many uniserial examples as \sq s of \dog\ modules \cite{Co05}.

Here we realize various uniserial $\K$-modules as \sq s of \psidog\ modules.  For $\lambda\in\bC$ and $p_1,\ldots,p_n\in\bZ/2$, define a {\em $(\lambda; p_1,\ldots,p_n)$ extension \/} to be a $\K$-extension of
$$ \bF(\lambda)\to\bF(\lambda+p_1)^{2p_1\Pi}\to\cdots\to
   \bF(\lambda+p_1+\cdots+p_n)^{2(p_1+\cdots+p_n)\Pi}. $$

\begin{prop} \label{length2K}
There are no \ind\ $\ds$-split extensions of\/ $\bF(\lambda)\to\bF(\lambda+p)^{2(p+1)\Pi}$ for any $(\lambda,p)$.  There are no \ind\ $\ds$-split $(\lambda; p)$ extensions unless $\dim^1_{\K,\ds}(\lambda,p) = 1$ (see Theorem~\ref{Ext1}), in which case there is a unique such extension.  It is uniserial unless $\lambda=0$ and $p\in\{2, 5/2, 3\}$, when it has $\bC\oplus\bF(p)^{2p\Pi}$ as a submodule.  In these cases its quotient by $\bC$ is the uniserial extension of\/ $\bF(1/2)^\Pi\to\bF(p)^{2p\Pi}$.

For $p=3/2$ and $-2\lambda\not\in\{0, 1, 2\}$, $p=2$ and $-2\lambda\not\in\{0, 1, 2, 3\}$, or $p=5/2$ and $-2\lambda\not\in\{0, \ldots, 4\}$, the above extension is a \sq\ of a \psidog\ module (up to a possibly odd equivalence).
\end{prop}

\meno{\it Proof.\/}
First note that $\bF(\lambda)$ is $\K$-\irr\ unless $\lambda=0$, when it is an \ind\ extension of $\bF(1/2)^\Pi \to \bC$.  Therefore the first sentence is a corollary of Theorem~\ref{Ext1} and Lemma~\ref{length2}, unless $\lambda=0$.  By Proposition~\ref{F012coho}c, there are no non-trivial $\ds$-split extensions of $\bC\to\bF(\mu)^{m\Pi}$ unless $\mu=3/2$ and $m=1$, in which case there is one.  We leave the rest of the first paragraph to the reader.

In order to prove the second paragraph, for $i_1<i_2<\cdots<i_a$ in $\bZ^+/2$ and $k-q\not\in \bN/2$ we define a subspace $\Psi^k_\ell(\nu,q)_{i_1,\ldots,i_a}$ of $\Psi^k_\ell(\nu,q)$ by
$$ \Psi^k_\ell(\nu,q)_{i_1,\ldots,i_a}\ :=\ \CQ_{\nu,q}
   \bigoplus_{j\in\bN\backslash\{2i_1,2i_2,\ldots,2i_a\}}
   \bF(q-k+j/2)^{(\ell+j)\Pi}. $$
For $k-q\in i_a+\bN/2$ we modify the definition to fit Corollary~\ref{URT}: let $\Psi^k_\ell(\nu,q)_{i_1,\ldots,i_a}$ be the pull-back to $\Psi^k_\ell(\nu,q)$ of the subspace of $\Psi^k_\ell/\Psi^{q-1/2}_{\ell-1-2(k-q)}$ obtained by applying $\CQ_{\nu,q}$ to the sum of those $\bF(q-k+j/2)^{(\ell+j)\Pi}$ with $0\le j\le 2p$ and $j\not\in\{2i_1,\ldots,2i_a\}$.

By Theorem~\ref{psimain} and Corollary~\ref{URT}, $\Psi^k_\ell(\nu,q)_{i_1,\ldots,i_a}$ is $\K$-invariant whenever $\{i_1,\ldots,i_a\}$ is a subset of $\{1/2,1\}$.  Suppressing $(\nu,q)$, we find that
\begin{eqnarray*}
   &&\bigl((\Psi^{q-\lambda}_\ell)_{\frac{1}{2},1} / 
   (\Psi^{q-\lambda-2}_\ell)\bigr)^{\ell\Pi}
   \mbox{\rm\ is a $(\lambda; \frac{3}{2})$ extension,} \\[4pt]
   &&\bigl((\Psi^{q-\lambda}_\ell)_{\frac{1}{2},1} / 
   (\Psi^{q-\lambda-3/2}_{\ell+1})_{\frac{1}{2}}\bigr)^{\ell\Pi}
   \mbox{\rm\ is a $(\lambda; 2)$ extension,} \\[4pt]
   &&\bigl((\Psi^{q-\lambda}_\ell)_{\frac{1}{2},1} / 
   (\Psi^{q-\lambda-3/2}_{\ell+1})_1\bigr)^{\ell\Pi}
   \mbox{\rm\ is a $(\lambda; \frac{5}{2})$ extension.}
\end{eqnarray*}
These extensions are \ind\ \iff\ $b^\ell_{\lambda+p,\lambda}(\nu,q)\not=0$ for $p=3/2$, $2$, and~$5/2$, respectively, which holds for generic values of $(\nu,q)$.  $\Box$

\meno{\bf Remark.}
The $\bigl((-7\pm\sqrt{33})/4;\th 4\bigr)$ extensions are probably also \sq s of \psidog\ modules.  As evidence, restate Corollary~\ref{specsplit} as follows:
\begin{equation} \label{specsplit2}
   \Psi^k_0({\ts \frac{1}{4}-\frac{q}{2}},\th q)_{1,\th 
   \frac{3}{2},\th 3,\th \frac{7}{2},\ldots}\ \and\ 
   \Psi^k_1({\ts \frac{1}{4}-\frac{q}{2}},\th q)_{\frac{1}{2},
   \th 1,\th \frac{5}{2},\th 3,\ldots}
\end{equation}
are $\K$-invariant.  The reader may check that if $4\nu=1-2q$ and both $b_{\lambda+2,\lambda}^0(\nu,q)$ and $b_{\lambda+4,\lambda+5/2}^1(\nu,q)$ are zero, then (suppressing $(\nu,q)$)
$$ (\Psi^{q-\lambda}_0)_{\frac{1}{2},\th 1,\th \frac{3}{2},\th 2} / 
   (\Psi^{q-\lambda-5/2}_1)_{\frac{3}{2}} 
   \mbox{\rm\ is a $(\lambda; 4)$ extension.} $$
Applying Theorem~\ref{main}, we find that this occurs only when $4\lambda=-7\pm\sqrt{33}$ and $2q=-\lambda-3$.  Thus these cases resolve the question affirmatively provided that they are in fact uniserial, \ie\ provided that they satisfy $b^0_{\lambda+4,\lambda}\not=0$.  We expect that this is so, but we have not verified it.

\begin{prop} \label{length3K}
There exists a uniserial $\ds$-split $(\lambda; p_1,p_2)$ extension \iff\ there exist uniserial $\ds$-split $(\lambda; p_1)$ and $(\lambda+p_1; p_2)$ extensions and $\beta_{p_2}(\lambda+p_1)\cup\beta_{p_1}(\lambda)$ is a 2-coboundary.  When such an extension exists, it is unique and we will call it $\bF(\lambda; p_1, p_2)$.

\begin{enumerate}

\item[(a)]  $\bF(\lambda;\th 3/2,\th 3/2)$ exists \iff\ $-2\lambda\not\in\{0, 1, 4, 5\}$.  It is a \psidog\ \sq\ for $-2\lambda\not\in\{0,\ldots, 5\}$.

\smallbreak
\item[(b)]  $\bF(\lambda;\th 2,\th 3/2)$ exists \iff\ $-2\lambda\not\in\{0, 3, 5\}$.  It is a \psidog\ \sq\ for $-2\lambda\not\in\{0,\ldots,6\}$.

\smallbreak
\item[(c)]  $\bF(\lambda;\th 3/2,\th 2)$ exists \iff\ $-2\lambda\not\in\{1, 3\}$.  It is a \psidog\ \sq\ for $-2\lambda\not\in\{0,\ldots, 6\}$.

\smallbreak
\item[(d)]  $\bF(\lambda;\th 5/2,\th 3/2)$ exists \iff\ $-4\lambda\not\in\{0, 4, 12, 7\pm\sqrt{33}\}$.  It is a \psidog\ \sq\ for $-2\lambda\not\in\{0,\ldots,7\}$ and $4\lambda\not=-7\pm\sqrt{33}$.

\smallbreak
\item[(e)]  $\bF(\lambda;\th 3/2,\th 5/2)$ exists \iff\ $-4\lambda\not\in\{2, 6, 10, 7\pm\sqrt{33}\}$.  It is a \psidog\ \sq\ for $-2\lambda\not\in\{0,\ldots,7\}$ and $4\lambda\not=-7\pm\sqrt{33}$.

\smallbreak
\item[(f)]  $\bF(\lambda;\th 2,\th 2)$ exists \iff\ $-4\lambda\not\in\{0, 8, 7\pm\sqrt{33}\}$.  It is a \psidog\ \sq\ for $-2\lambda\not\in\{0,\ldots,7\}$ and $4\lambda\not=-7\pm\sqrt{33}$.

\smallbreak
\item[(g)]  $\bF(\lambda;\th 5/2,\th 2)$ exists \iff\ $-4\lambda\not\in\{0, 4, 8, 10, 13\}$.  It is a \psidog\ \sq\ for $-2\lambda\not\in\{0,\ldots,8\}$.

\smallbreak
\item[(h)]  $\bF(\lambda;\th 2,\th 5/2)$ exists \iff\ $-4\lambda\not\in\{0, 3, 8, 12\}$.  It is a \psidog\ \sq\ for $-2\lambda\not\in\{0,\ldots,8\}$.

\smallbreak
\item[(i)]  $\bF(-4;\th 3/2,\th 3)$, $\bF(-9/2;\th 5/2,\th 5/2)$, and $\bF(-9/2;\th 2,\th 3)$ exist.

\end{enumerate}
\end{prop}

\meno{\em Proof.\/}
For the first paragraph, use Lemma~\ref{length3}, Theorem~\ref{Ext1}, and Proposition~\ref{length2K}.  For the first sentences of~(a) through~(h), we must avoid those values of $\lambda$ for which either $\dim^1_{\K,\ds}(\lambda,p_1)$ or $\dim^1_{\K,\ds}(\lambda+p_1,p_2)$ is zero, or a $(0;p)$ extension with $p=2$, $5/2$, or $3$ occurs as a \sq\ (as these are not uniserial), or $[\beta_{p_2}\cup\beta_{p_1}]\not=0$ (in $\ds$-relative cohomology).  

Clearly $[\beta_{p_2}\cup\beta_{p_1}]=0$ \iff\ $[\o\beta_{p_2}\cup\o\beta_{p_1}]=0$.  By direct calculation, $\o\beta_{p_2}\cup\o\beta_{p_1}(\bigwedge^2 e_{3/2})$ is never zero.  As we saw in the proof of Proposition~\ref{Ext2}, $B^2_{\lambda,p}$ is spanned by $\partial\o\beta_p$, 
and $\dim Z^2_{\lambda,p} = 1$ for $p=3$, $7/2$, $4$, and~$9/2$.  Therefore $[\o\beta_{p_2}\cup\o\beta_{p_1}]=0$ \iff\ $\partial\o\beta_{p_1+p_2}\not=0$.  The first sentences of~(a)-(h) now follow from Theorem~\ref{Ext1} and Propositions~\ref{Ext2} and~\ref{length2K}.

For the second sentences of~(a)-(h) we use the modules $\Psi^k_\ell(\nu,q)_{i_1,\ldots,i_a}$ from the proof of Proposition~\ref{length2K}.  The parameters $(\nu,q)$ will be suppressed.  The cases $-2\lambda=0,\ldots, 2(p_1+p_2)-1$ are missing because they lead to resonant \psidog\ \sq s.

For~(a), note that (up to a possibly odd equivalence)
$$ (\Psi^{q-\lambda}_\ell)_{\frac{1}{2},1} / 
   (\Psi^{q-\lambda-2}_\ell)_1
   \mbox{\rm\ is a $(\lambda; \frac{3}{2}, \frac{3}{2})$ extension.} $$
It is uniserial for $b^\ell_{\lambda+3/2,\lambda} \not= 0 \not= b^{\ell+1}_{\lambda+3,\lambda+3/2}$, which is true for most $(\nu,q)$.

For~(b), check that if $b^\ell_{\lambda+3/2,\lambda}=0$ then
$$ (\Psi^{q-\lambda}_\ell)_{\frac{1}{2},1, \frac{3}{2}} / 
   (\Psi^{q-\lambda-5/2}_{\ell+1})_1
   \mbox{\rm\ is a $(\lambda; 2, \frac{3}{2})$ extension.} $$
It is uniserial for $b^\ell_{\lambda+2,\lambda} \not= 0 \not= b^\ell_{\lambda+7/2,\lambda+2}$, which is the case for most values of $(\nu,q)$ on the curve $b^\ell_{\lambda+3/2,\lambda}=0$.

Case~(c) is essentially dual to~(b).  More precisely, working over $\Soo$ as in Section~\ref{S11}, the dual of a $(\lambda; p_1,\ldots, p_n)$ extension is a $(1/2-\lambda-\sum_i p_i; p_n,\ldots,p_1)$ extension (with reversed parity if $\sum_i 2p_i$ is even).  Moreover, \sq s of \psidog\ modules are dual-closed: by Proposition~\ref{Adler},
$$ \h\Psi^k_\ell(\nu,q)/ \h\Psi^{k-p-1/2}_{\ell+2p+1}(\nu,q) \and
   \h\Psi^{p-1/2-k}_{\ell+2p+1}(\nu+q,-q)/ \h\Psi^{-1-k}_\ell(\nu+q,-q) $$
are dual.  The reader may use Proposition~\ref{SNCRsym} to convert these statements to a proof that~(b), (d), and~(g) imply~(c), (e), and~(h), respectively.

For~(d), if $b^\ell_{\lambda+3/2,\lambda} = b^\ell_{\lambda+4,\lambda+2} =0$, then
$$ (\Psi^{q-\lambda}_\ell)_{\frac{1}{2},1, \frac{3}{2}} / 
   (\Psi^{q-\lambda-2}_\ell)_{\frac{1}{2}, 2}
   \mbox{\rm\ is a $(\lambda; \frac{5}{2}, \frac{3}{2})$ extension.} $$
For uniseriality we need $b^\ell_{\lambda+5/2,\lambda} \not= 0 \not= b^{\ell+1}_{\lambda+4,\lambda+5/2}$.  The two equalities are met for $\ell=0$, $c(\nu,q)=0$, and $-4q=2\lambda+5$.  At these values the necessary inequalities fail only for $-6\lambda=5$, $13$, and~$17$.  To fill these gaps we switch to the case $\ell=1$, where we need $16c^2=4q\lambda+2q+1$ and $2q(2\lambda-5)=(2\lambda+3)(2\lambda+7)$.  Here the necessary inequalities fail only for $\lambda = (5k+8)/(2k-16)$, $k\in\{0,1,2,3,4\}$, different gaps.

Proceed similarly for~(f), using the two spaces
$$ (\Psi^{q-\lambda}_0)_{\frac{1}{2}, 1, \frac{3}{2}} / 
   (\Psi^{q-\lambda-5/2}_1)_{\frac{3}{2}}
   \and (\Psi^{q-\lambda}_0)_{\frac{1}{2}, 1, \frac{3}{2}, \frac{5}{2}} / 
   (\Psi^{q-\lambda-3}_0)_1. $$
They are $(\lambda; 2, 2)$ extensions for $b^0_{\lambda+3/2,\lambda} = b^1_{\lambda+4,\lambda+5/2} =0$, and for $b^0_{\lambda+3/2,\lambda} = b^0_{\lambda+5/2,\lambda} =0$, respectively.  In both cases they are uniserial for $b^0_{\lambda+2,\lambda} \not= 0 \not= b^0_{\lambda+4,\lambda+2}$.  Taking $c(\nu,q)=0$ and working out the arithmetic gives the result.

For~(g) we use~(\ref{specsplit2}): for $c(\nu,q)=0$, the spaces
$$ (\Psi^{q-\lambda}_0)_{\frac{1}{2}, 1, \frac{3}{2}, 3, \frac{7}{2}} / 
   (\Psi^{q-\lambda-2}_0)_{\frac{1}{2}, 1, \frac{3}{2}, \frac{5}{2}}
   \and (\Psi^{q-\lambda}_0)_{\frac{1}{2}, 1, \frac{3}{2},2, 3, \frac{7}{2}} / 
   (\Psi^{q-\lambda-4}_0)_{\frac{1}{2}} $$
are $(\lambda; 5/2, 2)$ extensions for $b^0_{\lambda+9/2,\lambda+2}=0$ and for $b^0_{\lambda+2,\lambda}=0$, respectively.  In both cases they are uniserial for $b^0_{\lambda+5/2,\lambda} \not= 0 \not= b^1_{\lambda+9/2,\lambda+5/2}$.  Arithmetic gives the desired \sq s.

For~(i), note that any $(\lambda; p_1,\ldots,p_n)$ extension with $\lambda=-\sum_i p_i$ contains $\bF(0)$ and hence $\bC$ as submodules.  Division by $\bC$ gives a $(\lambda, p_1, \ldots, p_{n-1}, p_n+1/2)$ extension.  $\Box$

\meno
Next we construct some uniserial modules of length~4 as \sq s of \psidog\ modules.  Here various special values of $\lambda$ arise.  As far as we could determine, Proposition~\ref{length4K} gives all length~4 $\ds$-split uniserial \sq s.  In particular, we could find no $(\lambda;p_1,p_2,p_3)$ \sq s with $(2p_1,2p_2,2p_3)$ equal to any of $(3,4,3)$, $(5,3,3)$, $(3,5,3)$, $(4,3,4)$, $(5,4,3)$, $(5,3,5)$, or $(4,5,4)$.  We were also unable to find any uniserial \sq s of length~$\ge 5$.

\begin{prop} \label{length4K}
In the following cases there is a unique uniserial $\ds$-split $(\lambda,p_1,p_2,p_3)$ extension, which arises as a \sq\ of a \psidog\ module.

\begin{enumerate}

\item[(a)]
$(\lambda; \frac{3}{2}, \frac{3}{2}, \frac{3}{2})$ for $-2\lambda\not\in\{0, 1,\ldots, 8\}$.
\smallbreak

\item[(b)]
$(1; 2, \frac{3}{2}, \frac{3}{2})$ and the dual case $(-\frac{11}{2}; \frac{3}{2}, \frac{3}{2}, 2)$.
\smallbreak

\item[(c)]
$(\frac{-9\pm\sqrt{57}}{4}; 2, 2, \frac{3}{2})$ and the dual case $(\frac{-11\pm\sqrt{57}}{4}; \frac{3}{2}, 2, 2)$.
\smallbreak

\item[(d)]
$(\frac{-11\pm\sqrt{89}}{4}; \frac{5}{2}, \frac{3}{2}, 2)$ and the dual case $(\frac{-11\pm\sqrt{89}}{4}; 2, \frac{3}{2}, \frac{5}{2})$.
\smallbreak

\item[(e)]
$(\frac{-11\pm\sqrt{73}}{4}; 2, \frac{5}{2}, \frac{3}{2})$ and the dual case $(\frac{-11\pm\sqrt{73}}{4}; \frac{3}{2}, \frac{5}{2}, 2)$.
\smallbreak

\item[(f)]
The dual cases $(-7; 2, 2, 2)$ and $(3/2; 2, 2, 2)$.
\smallbreak

\item[(g)]
$(\frac{-11\pm\sqrt{105}}{4}; \frac{5}{2}, 2, 2)$ and the dual case $(\frac{-13\pm\sqrt{105}}{4}; 2, 2, \frac{5}{2})$.
\smallbreak

\end{enumerate}
\end{prop}

\meno{\em Proof.\/}
In all cases, uniqueness follows from the cup equation (see Section~\ref{LengthN}) and Theorem~\ref{Ext1}.  Since \psidog\ \sq s are dual closed (see the proof of Proposition~\ref{length3K}), it is only necessary to prove one of the two dual statements in each case.  Throughout we use the fact that if $\Psi^{q-\lambda}_\ell(\nu,q)$ has a $(\lambda; p_1, p_2, p_3)$ \sq, it is uniserial \iff\ 
$$ b^\ell_{\lambda+p_1,\lambda} \not=0, \quad 
   b^{\ell+2p_1}_{\lambda+p_1+p_2,\lambda+p_1} \not=0, \quad
   b^{\ell+2p_1+2p_2}_{\lambda+p_1+p_2+p_3,\lambda+p_1+p_2} \not=0. $$

For~(a), note that if $b^\ell_{\lambda+2,\lambda}$ and $b^{\ell+1}_{\lambda+9/2,\lambda+5/2}$ are zero then
$$ (\Psi^{q-\lambda}_\ell)_{\frac{1}{2}, 1, 2} / 
   (\Psi^{q-\lambda-5/2}_{\ell+1})_{\frac{1}{2}, 2}
   \mbox{\rm\ is a $(\lambda; \frac{3}{2}, \frac{3}{2}, \frac{3}{2})$ extension.} $$
The values $-2\lambda\in\{0,1,\ldots,8\}$ involve resonant \sq s.  Taking $\ell=0$ and applying the formulas of Theorem~\ref{main}, we find that for all other $\lambda$ except $-k/3$, $k\in\{4,5,7,8\}$, we can pick $q$ and $c$ so that both the \sq\ and uniseriality conditions are satisfied (the arithmetic is elementary, but delicate).  Taking $\ell=1$, we can satisfy the conditions for all $\lambda$ except $-4k/7$, $k\in\{1,2,\ldots,6\}$.  Thus there are no gaps.

For~(b) we must take $c(\nu,q)\not=0$.  For~(c)-(g) we take $c=0$, so that certain $b^\ell_{rs}$ are zero by Corollary~\ref{specsplit} (see also~(\ref{specsplit2})).  We list only the appropriate family of \sq s in which to look for satisfactory values of $\lambda$, along with the conditions necessary for the \sq\ to exist (in addition to $c=0$ for~(c)-(g)).

\begin{enumerate}

\item[(b)]
$(\Psi^{q-\lambda}_0)_{\frac{1}{2}, 1, 2, \frac{5}{2}} / (\Psi^{q-\lambda-7/2}_1)_{\frac{3}{2}}$, $b^0_{\lambda+2,\lambda}=b^0_{\lambda+5/2,\lambda}=b^1_{\lambda+5,\lambda+7/2}=0$.
\smallbreak

\item[(c)]
$(\Psi^{q-\lambda}_1)_{\frac{1}{2}, 1, \frac{3}{2}, \frac{5}{2}, 3} / (\Psi^{q-\lambda-7/2}_0)_{\frac{1}{2}, 2}$, $b^1_{\lambda+3/2,\lambda}=b^0_{\lambda+11/2,\lambda+7/2}=0$.
\smallbreak

\item[(d)]
$(\Psi^{q-\lambda}_0)_{\frac{1}{2}, 1, \frac{3}{2}, 2, 3, \frac{7}{2}} / (\Psi^{q-\lambda-9/2}_1)_{\frac{3}{2}}$, $b^0_{\lambda+2,\lambda}=b^1_{\lambda+6,\lambda+9/2}=0$.
\smallbreak

\item[(e)]
$(\Psi^{q-\lambda}_0)_{\frac{1}{2}, 1, \frac{3}{2}, \frac{5}{2}, 3, \frac{7}{2}} / (\Psi^{q-\lambda-4}_0)_{\frac{1}{2}, 2}$, $b^0_{\lambda+5/2,\lambda}=b^0_{\lambda+6,\lambda+4}=0$.
\smallbreak

\item[(f)]
$(\Psi^{q-\lambda}_0)_{\frac{1}{2}, 1, \frac{3}{2}, \frac{5}{2}, 3, \frac{7}{2}} / (\Psi^{q-\lambda-9/2}_1)_{\frac{3}{2}}$, $b^0_{\lambda+5/2,\lambda}=b^1_{\lambda+6,\lambda+9/2}=0$.
\smallbreak

\item[(g)]
$(\Psi^{q-\lambda}_0)_{\frac{1}{2}, 1, \frac{3}{2}, 2, 3, \frac{7}{2}} / (\Psi^{q-\lambda-4}_0)_{\frac{1}{2}, \frac{5}{2}}$, $b^0_{\lambda+2,\lambda}=b^0_{\lambda+13/2,\lambda+4}=0$.  $\Box$
\smallbreak

\end{enumerate}

\subsection{Equivalences and symmetries} \label{Equivs}
The study of the equivalence classes and symmetries of modules of \dog s on tensor densities was initiated by Duval and Ovsienko \cite{DO97}, who treated the modules of \dog s of order~2 over arbitrary manifolds.  Subsequently the cases of modules of higher order \dog s for manifolds of dimension $m>1$ and $m=1$ were treated in \cite{LMT96} and \cite{GO96}, respectively.  Then Lecomte and Ovsienko \cite{LO99} considered the more general case of modules of subquotients of \psdog s.

The 1-dimensional case is special: there the \dog\ modules have an invariant double symbol and consequently more equivalences.  Here we give an informal discussion of the case of subquotients of the finer filtration of the \dog\ modules of $\K$, which was specifically formulated at the end of \cite{GMO07}.  We have obtained precise results which will be given in a forthcoming paper.

For $\nu$, $q$, $\lambda\in\bC$ and $j\in\bZ^+$, define
\begin{eqnarray*}
   && \SQ_{\nu,q}^\ell(\lambda, j) := 
   \Psi^{q-\lambda}_\ell(\nu,q) / \Psi^{q-\lambda-j/2}_{\ell+j}(\nu,q),
   \mbox{\rm\ an extension of}\\
   &&\bF(\lambda)^\ell\to\bF(\lambda+1/2)^{\ell+1}
   \to\cdots\to \bF(\lambda+(j-1)/2)^{\ell+j-1}
\end{eqnarray*}
(see Section~\ref{LengthN} for extensions).  Here we restrict to non-resonant \sq s: the case that $2\lambda\in-\bN$ and $2\lambda+j-1\in\bZ^+$ is excluded (see Theorem~\ref{psimain} and Corollary~\ref{URT}), so the extension is $\ds$-split.

Duval and Ovsienko asked two questions.  First, when are $\SQ_{\nu,q}^\ell(\lambda,j)$ and $\SQ_{\nu',q'}^{\ell'}(\lambda',j')$ equivalent $\K$-modules?  Clearly equivalence requires $\lambda'=\lambda$ and $j'=j$.  We allow odd equivalences, so $\ell'$ need not be $\ell$.  Each choice of $\ell'$ leads to conditions on $\nu$ and $q$.

Second, what is the $\K$-endomorphism ring $\End_\K$ of $\SQ_{\nu,q}^\ell(\lambda,j)$?  In the non-resonant case any equivalence must drop to scalars on the composition series modules $\bF(\lambda+n/2)^{\ell+n}$.  Therefore $\End_\K$ is simply $\bC^e$, where $e$ is the number of indecomposable summands of $\SQ_{\nu,q}^\ell(\lambda,j)$.

At $j\le 3$, Proposition~\ref{sesqui} shows that $\SQ_{\nu,q}^\ell(\lambda,j)$ splits as $\bigoplus_{n=0}^{j-1} \bF(\lambda+n/2)^{\ell+n}$, so all $(\nu,q)$ give equivalent modules with $\End_\K = \bC^j$.

At $j=4$, $\SQ_{\nu,q}^\ell(\lambda,4)$ belongs to one of two equivalence classes: for $b^\ell_{\lambda+3/2,\lambda}=0$ it is $\bigoplus_0^3 \bF(\lambda+n/2)^{\ell+n}$, while for $b^\ell_{\lambda+3/2,\lambda}\not=0$ it is equivalent to the direct sum of $\bigoplus_1^2 \bF(\lambda+n/2)^{\ell+n}$ and the \ind\ $(\lambda;3/2)$ extension of Proposition~\ref{length2K}.  The first case is exceptional and has $\End_\K = \bC^4$, while the second is generic and has $\End_\K = \bC^3$.

At $j=5$, the two \sq s $\SQ_{\nu,q}^\ell(\lambda,5)$ and $\SQ_{\nu',q'}^{\ell'}(\lambda,5)$ are equivalent \iff\ $b^\ell_{\lambda+m/2,\lambda+n/2}(\nu,q)$ and $b^{\ell'}_{\lambda+m/2,\lambda+n/2}(\nu',q')$ are either both zero or both non-zero for $(m,n)=(3,0)$, $(4,1)$, or $(4,0)$.  Thus there are eight possible equivalences classes, all of which do occur.  For generic values of $(\nu,q)$ all three coefficients are non-zero, the module is the direct sum of $\bF(\lambda+1)^\ell$ and an \ind\ module, and $\End_\K = \bC^2$.

At $j=6$ a new phenomenon occurs.  To our knowledge, in all settings in which these questions have been considered to date the class of modules with a given $(\lambda,j)$ have been either almost all equivalent or essentially all inequivalent.  Here the equivalence classes are curves in $(\nu,q)$ space.  

To be more precise, recall from Section~\ref{COD} that the two modules $\Psi(\nu,q)$ and $\Psi(1/2-q-\nu,q)$ are always equivalent, by conjugation of \dog s.  These two modules have opposite values of $c$, so it is natural to regard the equivalence classes of the modules $\SQ^\ell_{\nu,q}(\lambda,j)$ as subsets of $(c^2,q)$ space.  As such, for $j=6$ the equivalence classes form the pencil of conics through four fixed points depending only on $\lambda$ and $\ell$ (for $\ell=1$, one of the four points turns out to be double).  The modules are generically indecomposable, with $\End_\K = \bC$.

To prove this, one first checks that the condition from the case $j=5$ is still necessary for equivalence, now for
$$ (m,n)\in\{(3,0), (4,1), (5,2), (4,0), (5,1), (5,0)\}. $$
But here this is not sufficient, there is an additional condition:
$$ \bigl( b^\ell_{\lambda+5/2,\lambda} b^{\ell+1}_{\lambda+2,\lambda+1/2} \bigr) /
   \bigl( b^{\ell+1}_{\lambda+5/2,\lambda+1/2} b^\ell_{\lambda+2,\lambda} \bigr) $$
is an invariant.  Its level curves form the pencil of conics.

We have not yet analyzed the case $j=7$.  It may be that these modules are equivalent only to their conjugates.  Exceptional cases would be interesting.

\section{Proofs}  \label{Proofs}
\setcounter{lemma}{0}

In this section we prove Theorems~\ref{main} and~\ref{psimain}.  The proof of Corollary~\ref{uptoscalar} covers those $\pi^m_{rs}$ in Theorem~\ref{psimain} with $r-s\in 1-\bN/2$.  Note that $b^m_{rs}(\lambda,p)$ is a function of the three complex variables $\lambda$, $p$, and $s$ and the two discrete variables $r-s\in 1+\bZ^+/2$ and $m\in\bZ_2$.  It is clearly rational; all the constructions in this paper are.  It is defined for all $(\lambda,p)$ for $s\not\in-\bN/2$, so it is polynomial in $(\lambda,p)$.  The set on which it is defined for \dog s is Zariski dense, so any rational formula valid on \dog s is necessarily the unique such formula valid on all \psidog s.  

Thus we come down to proving the formulas in Theorem~\ref{main}.  We use the strategy of \cite{Co01}.  Consider the equation
$$ \pi^m_{rs}(\lambda,p)(e_{r-s})(\omega^s) = 
   b^m_{rs}(\lambda,p) \beta_{r-s}(s)(e_{r-s})(\omega^s). $$
Clearly there are scalars $B_{rs}$ and $P^m_{rs}(\lambda,p)$ such that
$$ \pi^m_{rs}(e_{r-s})(\omega^s) = P^m_{rs} \omega^r, \quad
   \beta_{r-s}(e_{r-s})(\omega^s) = B_{rs} \omega^r. $$
We will compute them independently and obtain $b^m_{rs}$ as $P^m_{rs}/B_{rs}$.

Let us begin with $B_{rs}$ (the easier of the two).  Check that
\begin{equation} \label{Brsfact2}
   \pi_\lambda(e_{1/2})^n (\omega^\lambda) =
   (2\lambda)(2\lambda+1)\cdots (2\lambda+{\ts \lf\frac{n-1}{2}}\rf)
   \xi^{2\{n/2\}} x^{\lf n/2\rf} \omega^\lambda.
\end{equation}

By Lemma~\ref{beta}, $\beta_{r-s}(s)\bigl(\ad(e_{1/2})^{2r-2s-3}(e_{3/2})\bigr)\omega^s$ is
\begin{equation} \label{betars}
   4\bigl(\sigma_{s,r-s}(e_{1/2})^{2r-2s-3}(\omega^{r-s}\o D^{2r-2s-3})\bigr)(\omega^s).
\end{equation}
By~(\ref{saction}), for any $i>j$ we have $\o D^i \pi_s(e_{1/2})^j (\omega^s) =0$.  Therefore, if we expand the power of $\sigma_{s,r-s}(e_{1/2})$ in~(\ref{betars}) using the fact that $\sigma_{s,r-s}(e_{1/2})$ is $L\bigl(\pi_r(e_{1/2})\bigr) - R\bigl(\pi_s(e_{1/2})\bigr)$ (see Section~\ref{Hom}), we find that those terms containing a positive power of $L\bigl(\pi_r(e_{1/2})\bigr)$ give zero.  Taking into account signs, this reduces~(\ref{betars}) to
\begin{displaymath}
   4(-1)^{\lf r-s+1/2 \rf} \omega^{r-s} \o D^{2r-2s-3} 
   \circ \pi_s(e_{1/2})^{2r-2s-3} (\omega^s).
\end{displaymath}
Applying~(\ref{Brsfact1}) and~(\ref{Brsfact2}), we arrive at
\begin{equation} \label{Brs}
   B_{rs} = 2^{1+2\{r-s\}} (2s)(2s+1)\cdots (2s-2+\lf r-s \rf).
\end{equation}
Note that $B_{s+3/2,s}=4$, $B_{s+2,s}=4s$, and $B_{s+5/2,s}=8s$.

Now we turn to $P^m_{rs}(\lambda,p)$.  An elementary weight argument shows that 
$$ \pi(\lambda,p)(e_{r-s}) \bigl(\Pi^m(\omega^s)\bigr) =
   \sum_{j=0}^{2(r-s)} q_j \xi^{2\{j/2\}} x^{\lf j/2 \rf} 
   \Pi^{m+2(r-s)-j}(\omega^{r+j/2}) $$
for some scalars $q_j$, where by definition $P^m_{rs}(\lambda,p) = q_0$.

Recall that $\pi(\lambda,p)=\bigoplus_{\lambda,m}\pi^m_\lambda$ on $\ds$, and note that the image of $\pi_\lambda(e_{1/2})$ in $\bF(\lambda)$ is the space of $\omega^\lambda F$ such that $F(0,0)=0$.  (By $F(0,0)$ we mean the constant term of $F$.)  It follows that
$$ \pi(\lambda,p)(e_{r-s})\bigl(\Pi^m(\omega^s)\bigr) \equiv
   P^m_{rs} \Pi^{m+2(r-s)}(\omega^r)\ \mod\ 
   \Im\bigl(\pi(\lambda,p)(e_{1/2})\bigr). $$
Transferring this via $\CS_{\lambda,p}$, we find
\begin{equation} \label{proof1}
   \sigma_{\lambda,p}(e_{r-s})(\omega^p \o D^{2(p-s)}_m) \equiv
   P^m_{rs} \omega^p \o D^{2(p-r)}_{m+2(r-s)}\ \mod\ 
   \Im\bigl(\sigma_{\lambda,p}(e_{1/2})\bigr).
\end{equation}
We can easily calculate the left side of this equation, but we cannot compute $P^m_{rs}$ from it because the image of $\sigma_{\lambda,p}(e_{1/2})$ is too complicated.  Instead we use the step algebra element $\t s_{r-s}$ from Section~\ref{EPSA}.

By Lemma~\ref{step}, $\sigma_{\lambda,p}(\t s_{r-s})(\omega^p \o D^{2(p-s)}_m)$ is a \lwv\ of weight~$r$, hence by Lemma~\ref{Bol}a a multiple of $\omega^p \o D^{2(p-r)}_{m+2(r-s)}$ (assuming that it is defined, and we will see that it is).  On the other hand, $\t s_{r-s} \equiv e_{r-s}$ modulo $e_{1/2}\odU(\K)$, so $\sigma_{\lambda,p}(\t s_{r-s})(\omega^p \o D^{2(p-s)}_m)$ is congruent to~(\ref{proof1}) modulo the image of $\sigma_{\lambda,p}(e_{1/2})$.

Since $\omega^r$ is not in the image of $\pi_r(e_{1/2})$ and $\CS_{\lambda,p}$ is an $\ds$-equivalence, $\omega^p \o D^{2(p-r)}_{m+2(r-s)}$ is not in the image of $\sigma_{\lambda,p}(e_{1/2})$.  Thus
\begin{equation} \label{sigmats}
   \sigma_{\lambda,p}(\t s_{r-s})(\omega^p \o D^{2(p-s)}_m) =
   P^m_{rs} \omega^p \o D^{2(p-r)}_{m+2(r-s)}.
\end{equation}
The proof is concluded by computing the left hand side of~(\ref{sigmats}) essentially directly.  In order to do this, let us define the ``evaluation at zero'' operator $Z$ on \psidog s by
$$ Z(\omega^p G \o D^z_m) := \omega^p G(0,0) \o D^z_m. $$

\begin{lemma} \label{Z}
\begin{enumerate}

\item[(a)]  $Z(T\circ S)= Z\bigl(Z(T)\circ S\bigr)$.

\smallbreak
\item[(b)]  If $Z\bigl(\pi_{\lambda+p}(X)\bigr)=0$, then $Z\sigma_{\lambda,p}(X) = Z\sigma_{\lambda,p}(X) Z$.

\smallbreak
\item[(c)]  If $f(x)=0$, then $Z\bigl(\pi_\lambda(fD)\bigr)=0$ and $Z \sigma_{\lambda,p}(fD)(\omega^p G \o D^z_m) =$
$$ -\omega^p G(0,0) 
   Z\Bigl(\pxi\circ\o D^z_m\circ f + (-1)^m(2\lambda-1)\o D^z_m\circ\xi f' \Bigr). $$

\smallbreak
\item[(d)]  If $n>0$, then $Z\bigl(\pi_\lambda(e_n)\bigr) = 0$ for all $\lambda$.

\smallbreak
\item[(e)]  For $n<2\mu$, $Z \sigma_{\lambda,p}(e_{1/2}^n e_{\mu-n/2}) = \bigl(Z \sigma_{\lambda,p}(e_{1/2})\bigr)^n \bigl(Z \sigma_{\lambda,p}(e_{\mu-n/2})\bigr)$.

\end{enumerate}
\end{lemma}

This lemma is easily proven using the definitions in Sections~\ref{Flambda} and~\ref{DOGs}; for part~(c), use~(\ref{sigac1}) and~(\ref{sigac2}).  Some computation with the definition of $e_k$ and $\pi_\lambda$ leads to the following corollary:

\begin{cor} \label{Zen}
For $n\in\bZ^+$, $k$ and $\mu$ in $\bZ^+/2$, $n\le 2\mu$, and $m=0$ or~$1$,

\begin{enumerate}

\item[(a)]  $Z \sigma_{\lambda,p}(e_{1/2}^n) (\omega^p \o D^z_m) = E^{z,m}_n\omega^p \o D^{z-n}_{m-n}$, where
$$ E^{z,m}_n = \lb {(z-m)/2 \atop \lf n/2 \rf + 2(1-m)\{n/2\}}\rb
   \lb {(z+m)/2 -1+2\lambda \atop \lf n/2 \rf + 2m\{n/2\}}\rb. $$

\smallbreak
\item[(b)]  $Z \sigma_{\lambda,p}(e_k) (\omega^p \o D^z_m) = F^{z,m}_k \omega^p \o D^{z-2k}_{m-2k}$, where
\begin{eqnarray*}
   F^{z,m}_n &=& -(-1)^n \lb {(z-m)/2 \atop n} \rb \Bigl(\frac{z}{2} - 
   \frac{n}{2^m} + (n+1)\lambda \Bigr), \\[6pt]
   F^{z,m}_{n+1/2} &=& (-1)^n \lb {(z-m)/2 \atop n} \rb \Bigl(\frac{z+m}{2} - 
   n + 2m(n+1)\lambda \Bigr).
\end{eqnarray*}

\smallbreak
\item[(c)]  $Z \sigma_{\lambda,p}(e_{1/2}^n e_{\mu-n/2})(\omega^p \o D^z_m) = G^{z,m}_{\mu,n} \omega^p \o D^{z-2\mu}_{m-2\mu}$, where
$$ G^{z,m}_{\mu,n} = E^{z-2\mu+n,m-2\mu+n}_n\th F^{z,m}_{\mu-n/2}. $$

\end{enumerate}
\end{cor}

Now note that~(\ref{sigmats}) remains unchanged if we apply $Z$ to both sides.  Therefore, using Corollary~\ref{Zen} and~(\ref{ts}) of Section~\ref{EPSA}, we may rewrite~(\ref{sigmats}) as
\begin{equation} \label{Pmrs}
   P^m_{rs} = \sum_{n=0}^{2(r-s-1)} \t s_{r-s,n}(r) G^{2(p-s),m}_{r-s,n}.
\end{equation}
The rational functions $\t s_{r-s,n}(e_0)$ are given in Lemma~\ref{345}, although they are fully simplified only at $r-s=3/2$, $2$, and~$5/2$.  By~(\ref{Psdfn}), they can only have poles at $r$ in resonant cases.  Long but straightforward computation yields the simplified formulas of Theorem~(\ref{main}).  This completes the proof.

\section{The resonant and special cases} \label{Resonant}
\setcounter{lemma}{0}

Here we give the non-$\ds$-relative 1-cohomology of the $\bD(\lambda,p)$.  These results can be found in \cite{BBBBK}, but that article states an incorrect version of our Theorem~\ref{Ext1} (the 1-cocycles at $p=2$, $3$, and~$4$ are missed).  We also discuss the ``resonant conformal symbol'' in the sense of \cite{CS04}.

Recall the $\bF(0)$-valued 1-cocycle $\theta$ from Proposition~\ref{F012coho}.  Clearly $\theta$ may be regarded as an even $\bD^0(\lambda,0)$-valued 1-cocycle $\theta_0(\lambda)$ for all $\lambda$.

\begin{lemma} \label{urelDcoho}
Up to a scalar (and modulo $\dt$-relative classes), the only non-$\dt$-relative 1-cohomology classes of the modules $\bH(\lambda,p)$ are $[\theta_0(\lambda)]$ for $p=0$, and $[\SBol_{1/2}(0)\cup \theta_0(0)]$ for $(\lambda,p)=(0,1/2)$.
\end{lemma}

\meno {\em Proof.\/}
Suppose that $[\phi]$ is a non-$\dt$-relative 1-class.  By Lemmas~\ref{FCL0} and~\ref{FCL1}, there is a $\du$-relative $\dt$-invariant representative $\phi$ such that $\phi(e_0)$ is non-zero and $\K$-invariant.  By Lemma~\ref{Bol}d, either $p=0$ and $\phi(e_0)\propto 1$, or $(\lambda,p) = (0,1/2)$ and $\phi(e_0)\propto\SBol_{1/2}$.  The rest is easy.  $\Box$

\medbreak
Now recall the $\bF(1/2)$- and $\bF(3/2)$-valued 1-cocycles $\alpha$ and $\beta$ from Proposition~\ref{F012coho}.  For $2p\in 3+\bN$, Lemma~\ref{beta} lifts $\beta$ to the $\bD^{p-3/2}(\lambda,p)$-valued 1-cochain $\o\beta_p$ via $\CQ_{\lambda,p}$.  The $\dt$- but non-$\ds$-relative 1-classes arise from a more subtle lift of $\alpha$: for $p\in \frac{1}{2} +\bN$ and $c\not=0$, define
\begin{equation} \label{oalphap}
   \o\alpha_p(\lambda) := \ts \frac{1}{c}\Bigl(\partial\SBol_p(\lambda)
   + \frac{1}{16}(2p-1)(2p+1-16c^2)\o\beta_p\Bigr).
\end{equation}

\begin{lemma} \label{specialalpha}
$\o\alpha_p(\lambda)$ is an odd $\bD^{p-1/2}(\lambda, p)$-valued $\dt$-relative 1-cochain mapping $e_{1/2}$ to $-2\omega^p \o D^{2p-1}$.  Its symbol, viewed as an $\bF(1/2)$-valued 1-cochain, is $\alpha$.  Conjugation gives $\C\circ\o\alpha_p(\lambda) = (-1)^{p-1/2} \o\alpha_p(-\lambda-p+1/2)$.  

There is a continuous extension of the domain of $\o\alpha_p$ to $c=0$, and at $c=0$, $\o\alpha_p$ is not a coboundary for any $p$.  The 2-cochain $\partial\o\alpha_p$ is $\ds$-relative and $\bD^{p-3}$-valued, and its $\bF(3)$-symbol is $\frac{3}{8} {p-1/2 \choose 3} \bigl(p + \frac{1}{2} - 8c^2 \bigr)\th \beta\cup\beta$.
\end{lemma}

\meno{\it Proof.\/}
The $\dt$-relativity follows from Lemmas~\ref{Bol}b and~\ref{beta}, and the value at $e_{1/2}$ is given by~(\ref{SBoleoh}).  The symbol statement then follows from Lemma~\ref{Bol}b again.  The conjugation parity is straightforward.

For the continuous extension, note that $\partial\SBol_p$ is $\bD^{p-1/2}$-valued, so we may apply $\CS_{\lambda,p}$ to it.  Observe that in Theorem~\ref{main}, $b^1_{3/2, 0}$ is defined, even though $s=0$ is resonant.  We claim that $\CS_{\lambda,p}\bigl(\partial\SBol_p + b^1_{3/2,0} \o\beta_p\bigr)$ has no $\bF(3/2)$ term.  This may be proven by a (long) direct calculation, or by taking appropriate limits using Theorem~\ref{main}.  On the other hand, $\partial\SBol_p\propto\o\beta_p$ at $c=0$, so $\partial\SBol_p + b^1_{3/2,0} \o\beta_p$ must be zero there.  Since it is polynomial in $c$, it is divisible by $c$.

To prove that $\o\alpha_p$ is not a coboundary at $c=0$, use the facts that $[\alpha]\not=0$, and at $c=0$, $\partial\SBol_p \propto \o\beta_p$.  For the final sentence, apply $\partial$ to~(\ref{oalphap}) and use~(\ref{coboundary}) along with $\beta\cup\beta(\bigwedge^2 e_{3/2}) = 32 \omega^3$.  $\Box$

\begin{prop} \label{trelDcoho}
Up to a scalar, the only $\dt$-relative but non-$\ds$-relative 1-cohomology classes of the modules $\bH(\lambda,p)$ are $[\o\alpha_p(\lambda)]$ for $(\lambda,p)=(0, \frac{1}{2})$, $(-\frac{1}{2}, \frac{3}{2})$, and $(-1, \frac{5}{2})$.
\end{prop}

\meno {\it Proof.\/}
By Lemmas~\ref{IC} and~\ref{ExtN}, a $\dt$- but non-$\ds$-relative 1-class $[\phi]$ can occur only for $p\in\bZ^+/2$ and $c=0$.  Moreover, we may assume that $\phi$ is $\bD^{p-1/2}$-valued, $\phi(e_{1/2}) \propto \omega^p \o D^{2p-1}$, and $\phi$'s $\bF(1/2)$-symbol is proportional to $\alpha$.  

For $p\in\bZ^+$, (\ref{SBoleoh})~gives $\partial\SBol_p(e_{1/2}) = p\omega^p \o D^{2p-1}$, so $[\phi]$ is $\ds$-relative.  

For $p=1/2$, $3/2$, or~$5/2$ and $c=0$, Lemma~\ref{specialalpha} implies that $\o\alpha_p$ is a $\dt$-relative 1-cocycle whose class is non-$\ds$-relative.  If $\phi$ is a another such 1-cocycle, then $\phi(e_{1/2}) \propto \o\alpha_p(e_{1/2})$, so their classes are proportional modulo $\ds$-relative classes, of which there are none by Theorem~\ref{Ext1}.

For $p\in 7/2 +\bN$ and $c=0$, Lemma~\ref{specialalpha} implies that $\partial\o\alpha_p \not= 0$.  Here $\o\beta_p \propto \partial\SBol_p$, so it will suffice to prove that any $\dt$-relative $\bD^{p-1/2}$-valued 1-cochain $\psi$ of symbol $\alpha$ with an $\ds$-relative coboundary is of the form $\o\alpha_p + t\o\beta_p$ for some constant $t$.  

Since $\psi - \o\alpha_p$ is $\dt$-relative and $\bD^{p-1}$-valued, its $\bF(1)$-symbol is a multiple of $\SBol_2(-1) \circ \bX^{-1}$ (see the proof of Proposition~\ref{F012coho}).  But this does not have an $\ds$-relative coboundary, so in fact $\psi - \o\alpha_p$ must be $\bD^{p-3/2}$-valued.  Therefore it annihilates $\ds$ and has $\ds$-relative coboundary, so it is itself $\ds$-relative and therefore a multiple of $\o\beta_p$.  $\Box$

\medbreak
We remark that $\o\alpha_p$ is not the conformal lift $\CQ_{\lambda,p}\circ\ep^{2p+1}_{\bF(1/2)}\circ\alpha$ of $\alpha$, except at $p=1/2$ and~$3/2$.

The program of \cite{CS04} can be carried over to the resonant $\K$-modules $\Psi^k_\ell(\lambda,p)$ with $k-p\in\bN/2$, where $Q_\ds$ has doubled eigenvalues.  For example, consider the resonant special case in which $c=0$.  When $2(k-p)\not\equiv\ell$ modulo~$2$, the conformal symbol $\CS_{\lambda,p}$ exists but is not unique.  However, there is a unique choice of $\CS_{\lambda,p}$ which commutes with the conjugation $\C$.  The formulas for its $b^m_{rs}$ are those of the non-resonant case, but certain terms in the numerator and the denominator which are both zero at $c=0$ must be reduced in a specific way (one takes the resonant limit of the non-resonant special case).  The 1-cochains $\o\alpha_p$ with $c=0$ take values in these modules.

In the general resonant case, $\Psi^k_\ell$ is not completely reducible under $\ds$.  However, there is a $\dt$-invariant choice of $\CS_{\lambda,p}$ which is as close to $\ds$-invariant as possible, in the sense that the resulting matrix entries are zero on $\ds$ except on the {\em antidiagonal,\/} where $r+s=1/2$.  These antidiagonal entries are multiples of $\alpha_p := R(\ep)\circ\o\alpha_p$, and the multiples can be computed effectively.  

The entries with $s\le 0< r$ and $|r+s-1/2|>1$ are zero on $\ds$ but not necessarily $\ds$-relative.  They are difficult to compute, excepting those on the first and second sub and super antidiagonals where $r+s$ is $-1/2$, $0$, $1$, or $3/2$: these are given by the formulas of Theorem~\ref{main}, just as are the entries in the non-resonant zones $0<s<r$ and $s<r\le 0$.

\meno {\bf Acknowledgments.}  I am very grateful to Valentin Ovsienko for generously sharing with me his progress on the supersymmetric line, and I thank Paula Cohen for encouraging me to work on it.  Ian Musson gave helpful advice, and the referee's suggestions greatly improved both the organization and the content of the paper.

\section{Appendix: Representations of Lie superalgebras}  \label{LSA}
\setcounter{lemma}{0}

In this appendix we state several elementary but delicate properties of \lsa s and their \r s which streamline the calculations in the body of the paper.  All proofs are standard and are left to the reader.  

Fix a \lsa\ $\dg$ and \r s $\pi$ and $\sigma$ of $\dg$ on superspaces $V$ and $W$.  Occasionally we will need additional \r s $(\phi, U)$, $(\phi', U')$, $(\pi', V')$, and $(\sigma', W')$ of $\dg$.  We continue to write $V_\even$ and $V_\odd$ for the even and odd homogeneous subspaces of $V$, $V^\dg$ for the subspace of $\dg$-invariant vectors, and $\ep_V$ (or simply $\ep$) for the parity endomorphism:
$$ \ep_V:V \to V, \quad \ep_V(v):= (-1)^{|v|}v $$

Note that $\ep_\dg$ is an involution of $\dg$.  Write $\pi^\ep$ for the twist of $\pi$ by $\ep$: the \r\ on $V$ defined by $\pi^\ep := \pi\circ\ep_\dg$.  The map $\ep:V\to V$ is an even $\dg$-equivalence from $\pi$ to $\pi^\ep$.

\subsection{The representations $\Hom(V,W)$ and $V^*$} \label{Hom}
The space $\Hom(V,W)$ is naturally both a superspace and a \r\ of $\dg$, under the action $\Hom(\pi,\sigma)$ defined by
\begin{equation} \label{superhom}
   \Hom(\pi,\sigma)(X)(T) := \sigma(X)\circ T - 
   \ep_{\Hom(V,W)}^{|X|}(T)\circ \pi(X).
\end{equation}
We sometimes write $T^X$ for $\Hom(\pi,\sigma)(X)(T)$.  Here $\dg$ acts by superderivations: $(TS)^X = T^X S + (-1)^{|T||X|}TS^X$.

The elements of $\Hom(V,W)^\dg$ are called {\em intertwining maps,\/} or {\em $\dg$-maps.\/}  Bijective intertwining maps are called {\em $\dg$-equivalences.\/}  Note that intertwining maps $T$ do not satisfy $\sigma(X)\circ T = T\circ \pi(X)$ unless they are even.  In general they satisfy the ``rule of signs'': $\sigma(X)\circ T = (-1)^{|T||X|} T\circ\pi(X)$.  The space of intertwining maps is $\bZ_2$-graded, \ie\ equal to the sum of its even and odd parts.

Let $\pi_\triv$ be the even trivial \r\ $X\mapsto 0$ on $\bC^{1|0}$.  The {\em dual\/} $\pi^*$ of $\pi$ is defined to be the \r\ $\Hom(\pi, \pi_\triv)$ on $V^* = \Hom(V,\bC^{1|0})$.  Given any map $T:V\to W$, write $T^t:W^*\to V^*$ for its {\em transpose,\/} defined by $T^t(\mu)(v):=\mu(T(v))$.  The following facts are useful (note in particular that $(\pi^*)^*$ is not $\pi$):
\begin{displaymath}
   \pi^*(X) = -\pi(X)^t\circ \ep_{V^*}^{|X|}, \quad
   \ep_V^t = \ep_{V^*}, \quad
   (\pi^*)^* = \pi^\ep.
\end{displaymath}

\begin{lemma} \label{tau}
The map $\tau:\Hom(V,W)\to \Hom(W^*,V^*)$ defined by $\tau(T):= T^t\circ\ep_{W^*}^{|T|}$ is an even $\dg$-equivalence.
\end{lemma}

The next lemma is another illustration of the rule of signs.  Define maps
\begin{eqnarray*}
   &&L:\Hom(V,W)\to \Hom\bigl(\Hom(U,V),\Hom(U,W)\bigr), \\
   &&R:\Hom(U,V)\to \Hom\bigl(\Hom(V,W),\Hom(U,W)\bigr)
\end{eqnarray*}
by $L(T)(S):=T\circ S$ and $R(S)(T):=(-1)^{|T||S|}T\circ S$.

\begin{lemma} \label{RL}
The maps $L$ and $R$ are even $\dg$-maps.  In particular, if $T:V\to W$ and $S:U\to V$ are $\dg$-maps then so are $L(T)$ and $R(S)$.  Where defined, $L(T)\circ L(S) = L(T\circ S)$ and $R(S)\circ R(T) = (-1)^{|T||S|} R(T\circ S)$.
\end{lemma}

\subsection{Tensor products}
Given any linear maps $\nu:V\to V'$ and $\omega:W\to W'$, define
\begin{equation} \label{ots}
   \nu\ot_s\omega:= (\nu\circ\ep_V^{|\omega|})\ot\omega:V\ot W\to V'\ot W'.
\end{equation}
Note that $(\nu\ot_s\omega)(v\ot w)=(-1)^{|\omega||v|}\nu(v)\ot\omega(w)$, in accordance with the rule of signs.

The tensor product $V\ot W$ is naturally both a superspace and a \r\ of $\dg$.  Its parity is $|v\ot w| := |v|+|w|$, and its action is
$$ (\pi\ot\sigma)(X) := \pi(X)\ot_s 1 + 1\ot_s \sigma(X). $$
In the category of super \r s the tensor product operation is associative but supercommutative rather than commutative, in the sense that the $\dg$-equivalence from $\pi\ot\sigma$ to $\sigma\ot\pi$ is not $v\ot w\mapsto w\ot v$ but instead the map $s$ predicted by the rule of signs:
\begin{equation} \label{exchange}
   s: V\ot W \to W\ot V, \quad s(v\ot w) := (-1)^{|v||w|}w\ot v.
\end{equation}

The next lemma is not hard to prove directly but it can be tricky to find.  It can be systematically deduced using the following two facts.  First, we have the usual even $\dg$-equivalence from $W\ot V^*$ to $\Hom(V,W)$ defined by $(w\ot\lambda) (v) := w\lambda(v)$.

Second, the obvious ``identity'' map from $V^*\ot W^*$ to $(V\ot W)^*$ is not a $\dg$-equivalence, but the exchange map $\lambda\ot\mu\mapsto\mu\ot\lambda$ {\em is\/} a $\dg$-equivalence from $V^*\ot W^*$ to $(W\ot V)^*$ (there is no sign term because dualization reverses order).

\begin{lemma} \label{tensorhom}
The rule $\nu\ot\omega\mapsto \nu\ot_s\omega$ defines an even $\dg$-equivalence from $\Hom(V, V')\ot\Hom(W,W')$ to $Hom(V\ot W, V'\ot W')$.  In particular, if $\nu:V\to V'$ and $\omega:W\to W'$ are $\dg$-maps then $\nu\ot_s\omega$ is also a $\dg$-map.
\end{lemma}

The next lemma shows that composition behaves as expected.  We remark that it is equivalent to the fact that $\dg$ acts by superderivations on $\End(U\oplus V\oplus W)$ and preserves its lower triangular block subalgebras.

\begin{lemma}
The rule $T\ot S\mapsto T\circ S$ defines an even $\dg$-map
$$ \Comp:\Hom(V,W)\ot \Hom(U,V)\to \Hom(U,W). $$
\end{lemma}

\subsection{Tensor algebras}
Write $\bigotimes V := \bigoplus_0^\infty \bigotimes^n V$ for the tensor algebra on $V$, and $\bigotimes\pi:=\bigoplus_0^\infty \bigotimes^n\pi$ for the natural action of $\dg$ on it by superderivations.  The algebras $\S_s V$ of supersymmetric tensors and $\Lambda_s V$ of superalternating tensors are the quotients of $\bigotimes V$ by the ideals $\I_{\S_s}V$ and $\I_{\Lambda_s}V$ generated by the images of the endomorphisms $1-s$ and $1+s$ of $V\ot V$, respectively (see~(\ref{exchange})).  These are graded $\dg$-invariant ideals, so we have gradings 
$$ \S_s V = {\ts\bigoplus_{0}^\infty} \S_s^n V\ \and\ 
   \Lambda_s V = {\ts\bigoplus_0^\infty} \Lambda_s^n V $$
with associated $\dg$-actions $\S_s\pi = \bigoplus_0^\infty \S_s^n\pi$ and $\Lambda_s\pi = \bigoplus_0^\infty \Lambda_s^n \pi$.  As vector spaces, 
$$ \S_s V = \S V_\even\ot\Lambda V_\odd\ \and\ 
   \Lambda_s V = \Lambda V_\even \ot \S V_\odd, $$
where $\S$ and $\Lambda$ indicate the ordinary symmetric and alternating algebras.

Now note that the operation $\ot_s$ defined in~(\ref{ots}) is associative: if $\mu:U\to U'$ is a third linear map, then $\mu\ot_s(\nu\ot_s\omega)=(\mu\ot_s\nu)\ot_s\omega$, so $\mu\ot_s\nu\ot_s\omega$ is unambiguous.  Given any map $T:V\to W$, define a map $\bigotimes_s T:\bigotimes V\to \bigotimes W$ by $\bigotimes_s T:=\bigoplus_0^\infty\bigotimes_s^n T$, where 
$$ {\ts\bigotimes_s^n} T := T\ot_s T\ot_s\cdots\ot_s T 
   \mbox{\rm\ \ ($n$ copies of $T$)}. $$

\begin{lemma} \label{otsT}
If $T:V\to W$ is a $\dg$-map, then so are $\bigotimes_s^n T$ and $\bigotimes_s T$.  

If $T$ is an even $\dg$-map, then $\bigotimes_s T$ is a $\dg$-algebra homomorphism which carries $\I_{\S_s}V$ to $\I_{\S_s}W$ and $\I_{\Lambda_s}V$ to $\I_{\Lambda_s}W$.   Hence it factors through to $\dg$-algebra homomorphisms from $\S_sV$ to $\S_sW$ and from $\Lambda_sV$ to $\Lambda_sW$, which we also denote by $\bigotimes_s T$.  

If $T$ is an odd $\dg$-map, then $\bigotimes_s T$ carries $\I_{\S_s}V$ to $\I_{\Lambda_s}W$ and $\I_{\Lambda_s}V$ to $\I_{\S_s}W$.  Hence it factors through to $\dg$-maps from $\S_sV$ to $\Lambda_sW$ and from $\Lambda_sV$ to $\S_sW$, which we again denote by $\bigotimes_s T$.
\end{lemma}

\subsection{Parity reversal} \label{Parity}
At the level of superspaces, the parity-reversing functor $V\mapsto V^\Pi$ is simply tensoring with the odd trivial \r\ $\bC^{0|1}$.  At the level of \r s we must decide which side to tensor on: by convention, $V^\Pi := V\ot\bC^{0|1}$.  Thus $V^\Pi$ carries the \r\ $\pi^\Pi$ defined by $\pi^\Pi(X):=\pi(X)$ for all $X\in\dg$.  In particular, $(\pi^\Pi, V^\Pi)$ and $(\pi, V)$ are the same as \r s of an algebra on a vector space, but not as \r s of a \lsa\ on a superspace.  Note that $\ep_{V^\Pi} = -\ep_V$, and both $1:V\to V^\Pi$ and $\ep_V:V\to V^\Pi$ are odd bijections.  In fact, $1$ is not a $\dg$-equivalence but $\ep_V$ is.  Here we describe the interaction between parity reversal and Hom, tensor products, and tensor algebras.

As superspaces, $\Hom(V,W^\Pi)$ and $\Hom(V^\Pi,W)$ are equal to $\Hom(V,W)^\Pi$, while $\Hom(V^\Pi,W^\Pi)$ is equal to $\Hom(V,W)$.  The $\dg$-actions $\Hom(\pi,\sigma^\Pi)$ and $\Hom(\pi^\Pi,\sigma)$ are equal, as are $\Hom(\pi,\sigma)$ and $\Hom(\pi^\Pi,\sigma^\Pi)$.  However, due to the change in the sign of the $\ep$ in~(\ref{superhom}) these two pairs of actions are not the same.  Note that $\Hom(\pi,\sigma)^\Pi$ has the parity of the former pair but the $\dg$-action of the latter pair.  The next three lemmas are corollaries of Lemmas~\ref{RL}, \ref{tensorhom}, and~\ref{otsT}, respectively.

\begin{lemma} \label{leprep}
The maps $L(\ep_W)$ and $R(\ep_V)$ from $\Hom(V,W)$ to $\Hom(V^\Pi,W)$ are equal and are odd $\dg$-equivalences.  In particular, if $T:V\to W$ is a $\dg$-map then so is $\ep_W\circ T:V\to W^\Pi$, which may also be viewed as the $\dg$-map $(-1)^{|T|}T\circ\ep_V:V^\Pi\to W$.
\end{lemma}

\begin{lemma} \label{otep}
The \r s $V\ot W^\Pi$ and $(V\ot W)^\Pi$ are equal.  The map $\ep_V\ot 1$ is a $\dg$ equivalence from $V\ot W$ to $V^\Pi\ot W$.
\end{lemma}

\begin{lemma} \label{taep}
The map $\bigotimes_s\ep:\bigotimes V\to \bigotimes(V^\Pi)$ is a $\dg$-equivalence which factors through to $\dg$-equivalences from $\S_sV$ to $\Lambda_s(V^\Pi)$ and from $\Lambda_s V$ to $\S_s(V^\Pi)$.  Its homogeneous components $\bigotimes_s^n\ep$ have parity $(-1)^n$.
\end{lemma}

\subsection{Invariant pairings}
An {\em invariant pairing\/} of $V$ with $W$ is an intertwining map $B:V\ot W\to\bC^{1|0}$, \ie\ a map $B:V\ot W\to\bC$ such that
$$ B((Xv)\ot w) + (-1)^{|X||v|}B(v\ot(Xw)) =0. $$
The situation is supersymmetric with respect to $V$ and $W$, in the sense that $\t B:=B\circ s$ is an invariant pairing of $W$ with $V$ \iff\ $B$ is an invariant pairing of $V$ with $W$.

An invariant pairing of $V$ with itself is called an {\em invariant form\/} on $V$.  An invariant form $B$ on $V$ is {\em supersymmetric\/} if it is equal to $\t B$.  Odd supersymmetric forms are symmetric, while even supersymmetric forms are symmetric on $V_\even$ and skewsymmetric on $V_\odd$.  It can be shown that if $V$ admits a non-degenerate even symmetric invariant form, then $\dg_\odd$ acts by zero.  Note that if $B$ is any invariant pairing of $V$ with $W$, then $B\oplus\t B$ is a supersymmetric invariant form on $V\oplus W$, non-degenerate \iff\ $B$ is, and of the same parity as $B$.

Fix a non-degenerate invariant pairing $B$ of $V$ with $W$.  Given $v_0\in V$ and $T\in\End(W)$, define $v_0^B\in W^*$ and $T^B\in \End(V)$ by
\begin{displaymath}
   v_0^B(w):=B(v_0\ot w), \quad B\bigl(T^B(v)\ot w\bigr):=B\bigl(\ep^{|T|}(v)\ot T(w)\bigr).
\end{displaymath}

\begin{lemma} \label{tauB}
The map $v\mapsto v^B$ is a $\dg$-equivalence from $V$ to $W^*$ of parity $|B|$.  The map $T\mapsto T^B$ from $\End(W)$ to $\End(V)$ satisfies $(T\circ S)^B=(-1)^{|T||S|}S^B\circ T^B$ and $\sigma(X)^B = -\pi(X)$.  It is an even $\dg$-equivalence such that $(T^B)^{\t B}=T$ and $(\ep_W)^B=(-1)^{|B|}\ep_V$.
\end{lemma}

\subsection{\bf The universal enveloping algebra} \label{UEA}
The \uea\ $\dU(\dg)$ of $\dg$ is the quotient of $\bigotimes\dg$ by the ideal generated by all elements of the form $X\ot Y-(-1)^{|X||Y|}Y\ot X - [X,Y]$.  We will frequently write simply $\dU$ when $\dg$ is determined by the context.  Let $\dU^n$ be the degree filtration.  

We first recall the super \PBW\ theorem.  Given an element $\alpha$ of the symmetric group $S_n$, define the {\em superpermutation endomorphism\/} $\SP(\alpha)$ of $\bigotimes^n\dg$ as follows: if $X_1,\ldots, X_n$ are homogeneous elements of $\dg$, then 
$$ \SP(\alpha)(X_1\ot\cdots\ot X_n) :=
   \pm X_{\alpha^{-1}(1)}\ot\cdots\ot X_{\alpha^{-1}(n)}, $$
the sign being chosen so that the two sides are equal in $\S_s\dg$.  One checks that $\SP$ is a \r\ of $S_n$ on $\bigotimes^n\dg$ commuting with the adjoint action of $\dg$, and that the formula
$$ \Sym := {\ts\frac{1}{n!}}\sum_{\alpha\in S_n} \SP(\alpha) $$
defines a projection on $\bigotimes\dg$ with kernel $\I_{\S_s\dg}$.  The super \PBW\ theorem states that this projection drops to a $\dg$-equivalence $\Sym$ from $\S_s\dg$ to $\dU(\dg)$, the {\em symmetrizer map.\/}

Let us use $X$ and $Y$ for elements of $\dg$ and $\Omega$ and $\Theta$ for elements of $\dU$.  There are some useful automorphisms of $\dU$.  First, the involution $\ep_\dg$ of $\dg$ extends to an algebra involution of $\dU$.  This involution has a square root $\iota$, an anti-automorphism of $\dU$ of order~4 defined by
$$ \iota(X):=-(-i)^{|X|}X, \quad
   \iota(\Omega\Theta):=\iota(\Theta)\iota(\Omega). $$
One checks that $\iota^2=\ep$.

There is also a super anti-involution $\Omega\mapsto\Omega^T$, defined by
$$ X^T := -X, \quad (\Omega\Theta)^T := 
   (-1)^{|\Omega||\Theta|}\Theta^T\Omega^T. $$
It is related to dualization: $\pi^*(\Omega) = \pi(\Omega^T)^t\circ\ep_{V^*}^{|\Omega|}$.

Although we will not need this, we remark that $\iota$ and $T$ commute and their composite $j(\Omega):=\iota(\Omega^T)$ has square $\ep$ and is a superautomorphism: $j(\Omega\Theta)=(-1)^{|\Omega||\Theta|}j(\Omega)j(\Theta)$.  It has the curious formula
$$ j(X_1\cdots X_n) = (-i)^{\bigl(\sum_k|X_k|\bigr)^2} X_1\cdots X_n. $$

\subsection{The adjoint prime action and the super and ghost centers}
To our knowledge the ``adjoint prime'' \r\ $\ad'$ of $\dg$ on $\dU^\Pi$ was first defined in \cite{ABF97}:
$$ \ad'(X)(\Omega) := X\Omega -(-1)^{|X||\Omega+1|}\Omega X. $$
Direct computation shows that this is indeed a \r.  It does not preserve the degree filtration, but it increases the degree by at most $\dim(\dg_\odd)$: the \PBW\ theorem gives
$$ \ad'(\dU)(\dU^n) \subseteq \ad'(\Lambda\dg_\odd)(\dU^n)
   \subseteq \dU^{n+\dim(\dg_\odd)}. $$

As is pointed out in \cite{Gor00}, $\ad'$ is best viewed as the adjoint action on the bimodule $\dU^\Pi$.  Let us remark that the adjoint \r\ $\ad$ of $\dg$ on itself defines an algebra map $\ad:\dU\to\End(\dg)$ which intertwines the natural $\dg$ actions.  This same map intertwines the $\ad'$ action on $\dU^\Pi$ with $\Hom(\dg,\dg^\Pi)$, so $\ad'$ may be thought of as a pull-back of $\Hom(\ad,\ad^\Pi)$.  However, although $T\mapsto\ep\circ T$ is an odd $\dg$-equivalence from $\Hom(\dg,\dg)$ to $\Hom(\dg,\dg^\Pi)$, the $\ad$ and $\ad'$ actions on $\dU$ and $\dU^\Pi$ are in general not equivalent.

The {\em supercenter\/} of $\dU(\dg)$ is $\dZ(\dg):=\dU(\dg)^\dg$.  In \cite{Gor00} one finds the following definitions: the {\em anticenter\/} $\dZ'(\dg)$ and {\em ghost center\/} $\t\dZ(\dg)$ are  
$$ \dZ'(\dg):=\bigl(\dU(\dg)^\Pi\bigr)^{\ad'(\dg)}, \quad \t\dZ(\dg):= \dZ(\dg)+\dZ'(\dg). $$

The intersection $\dZ\cap\dZ'$ consists of those elements of $\dU^{\dg_\even}$ which annihilate $\dg_\odd$.  When this intersection is zero, $\t\dZ$ is $\bZ_2\times\bZ_2$-graded: the $(0,0)$-space is $\dZ_\even$ and the other three spaces are $\dZ_\odd$, $\dZ'_\even$, and $\dZ'_\odd$.  Note that
\begin{eqnarray*}
   \dZ_\even + \dZ'_\odd &=&
   \bigl\{\Omega\in\dU: X\Omega=\Omega X\ \forall\ X\in\dg\bigr\},\\
   \dZ'_\even + \dZ_\odd &=&
   \bigl\{\Omega\in\dU: X\Omega=(-1)^{|X|}\Omega X\ \forall\ X\in\dg\bigr\}.
\end{eqnarray*}

\def\eightit{\it} \def\bib{\bf}

\end{document}